\documentclass{article}
\usepackage{graphicx} 
\usepackage{amsmath}
\usepackage{amssymb}
\usepackage{amsthm}
\usepackage{bbm}
\usepackage{dirtytalk}
\usepackage{amsfonts}
\usepackage{xcolor}
\usepackage{enumitem}
\usepackage{url}
\usepackage{hyperref} 
\usepackage[a4paper, total={6in, 9in}]{geometry}
\theoremstyle{plain}
\newtheorem{corollary}{Corollary}

\newtheorem{theorem}{Theorem}
\newtheorem{proposition}{Proposition}
\newtheorem{lemma}{Lemma}
\newtheorem{definition}{Definition}
\newtheorem{example}{Example}
\newtheorem{remark}{Remark}
\newtheorem{assumption}{Assumption}

\title{Small noise fluctuations and large deviations of conservative SPDEs with Dirichlet boundary conditions}
\author{
Shyam Popat\thanks{Mathematical Institute, University of Oxford, Oxford OX2 6GG, UK (popat@maths.ox.ac.uk)}}
\date{\today}

\begin{document}

\maketitle

\small
\begin{center}
ABSTRACT. We establish a central limit theorem and large deviations principle that characterises small noise fluctuations of the generalised Dean--Kawasaki stochastic PDE. 
The fluctuations agree to first order with fluctuations of certain interacting particle systems, such as the zero range process, about their hydrodynamic limits.
Our main contribution is that we are able to consider stochastic PDEs on general $C^2$ bounded domains with Dirichlet boundary conditions. 
On the level of particles, the boundary condition corresponds to absorption or injection of particles at the boundary.
\end{center}
\vspace{10pt}
\textbf{MSC2020 subject classification:} \hspace{5pt}60F05, 60F10, 60H15, 82B21, 82B31, 60K35, 35Q70\\
\textbf{Keywords:}\hspace{5pt} Central limit theorem, Large deviations principle, Interacting particle system, Stochastic partial differential equation, Dean--Kawasaki equation, Dirichlet boundary conditions
\normalsize
\section{Introduction}
In this paper, we study the small noise fluctuations of solutions $\rho^{\epsilon,K}$ to the non-linear stochastic PDE for every $K\in\mathbb{N}$ and $\epsilon\in(0,1)$, 
\begin{equation}\label{generalised Dean-Kawasaki Equation Stratonovich with epsilon}
   \begin{cases}
    \partial_t\rho^{\epsilon,K}=\Delta\Phi(\rho^{\epsilon,K})-\sqrt{\epsilon}\nabla\cdot(\sigma(\rho^{\epsilon,K})\circ\dot{\xi}^K)-\nabla\cdot \nu(\rho^{\epsilon,K}), & \text{on} \hspace{5pt} U\times(0,T],\\
\Phi(\rho^{\epsilon,K})=\bar{f},&\text{on}\hspace{5pt} \partial U\times[0,T],\\
    \rho^{\epsilon,K}(\cdot,t=0)=\rho_0 ,&\text{on}\hspace{5pt} U\times\{t=0\},
   \end{cases} 
\end{equation}
where $U\subset\mathbb{R}^d$ denotes a $C^2$-bounded domain and $T>0$ is a fixed terminal time.
    For $K\in\mathbb{N}$, $\dot\xi^{K}$ denotes truncations of $\mathbb{R}^d$-valued space-time white noise defined in Definition \ref{definition truncated noise} below, and $\circ$ represents the fact that we are looking at Stratonovich noise.
The non-linearities $\Phi,\sigma$ and $\nu$ cover a wide range of relevant examples from the point of view of fluctuations of interacting particle systems, including the full range of fast diffusion and porous medium equations $\Phi(\xi)=\xi^m$ for every $m\in(0,\infty)$ and diffusion coefficients including the degenerate square root $\sigma(\xi)=\sqrt{\xi}$.
The non-negative initial condition $\rho_0$ is assumed to be $L^1(U)$-integrable.
The irregularity of solutions to \eqref{generalised Dean-Kawasaki Equation Stratonovich with epsilon} means that we can only express the boundary condition $\bar{f}$ indirectly via $\Phi(\rho^{\epsilon,K})$.\\
In Section \ref{sec: application to particle system} we motivate the study of generalised Dean-Kawasaki type equations \eqref{generalised Dean-Kawasaki Equation Stratonovich with epsilon} with an example of how they are useful in the theory of fluctuating hydrodynamics.
Briefly, the SPDEs correctly capture to first order the fluctuations of the zero range process particle system about its hydrodynamic limit. 
They also correctly describe the rare events exhibited by the particle system in the sense that the rate functions of the SPDE and the empirical measure of the particle system coincide.
In Section \ref{sec: background literature} we provide a background of the solution theory for equations of the form \eqref{generalised Dean-Kawasaki Equation Stratonovich with epsilon}, as well as relevant results in the study of central limit theorems and large deviation principles.
In the literature equations such as \eqref{generalised Dean-Kawasaki Equation Stratonovich with epsilon} are referred to as \say{conservative SPDEs} due to the fact the equation conserves mass (conserves $L^1(U)$-norm) on the torus.
However, on the bounded domain with Dirichlet boundary data, mass is not preserved so the description is slightly misleading.\\
The goal of the present work is to study fluctuations of solutions $\rho^{\epsilon,K}$ about the zero noise ($\epsilon\to0, K\to\infty$), hydrodynamic limit
\begin{equation}\label{heat equation hydrodynamic limit}
   \begin{cases}
    \partial_t\bar\rho=\Delta\Phi(\bar\rho)-\nabla\cdot \nu(\bar\rho), & \text{on} \hspace{5pt} U\times(0,T],\\
\Phi(\bar\rho)=\bar{f},&\text{on}\hspace{5pt} \partial U\times[0,T],\\
    \bar\rho(\cdot,t=0)=\rho_0 ,&\text{on}\hspace{5pt} U\times\{t=0\}.
   \end{cases} 
\end{equation}
In order to consider weak solutions of \eqref{generalised Dean-Kawasaki Equation Stratonovich with epsilon} we need the diffusion coefficient $\sigma$ to be sufficiently regular, and in order to apply It\^o's formula to obtain energy estimates we will need to add a regularisation to the equation, see Definition \ref{defn: regularised equation} below.
The first result we will prove is two quantitative law of large number estimates for the regularised equation $\rho^{n,\epsilon,K}$ in the form of $L^p(U\times[0,T])$-estimates. 
\begin{proposition}[Propositions \ref{ppn: energy estimate for difference of regularised equation and hydrodynamic limit},  \ref{prop: second energy estimate for rho n epsilon k} below]\label{ppn: energy estimates in introducition}
    Let $\rho^{n,\epsilon,K}$ be the weak solution of the regularised equation \eqref{equation for rho n epsilon k}, and let $\bar\rho$ denote solution to the regularised hydrodynamic limit equation \eqref{eq: regularised hydrodynamic limit}.
    Under Assumptions \ref{assumption on noise}, \ref{assumption of constant boundary data and constant initial condition}, \ref{asm: new assumptions} and \ref{asm: assumptions for well posedenss of equation}, for every $\epsilon\in(0,1)$ and $K\in\mathbb{N}$, there exists a constant $c\in(0,\infty)$ independent of $T, \alpha, \epsilon$, $n$ and $K$ such that for every $p\geq 2$, we have the $p$-independent estimate
    \begin{multline*}
    \frac{1}{Tp(p-1)}\mathbb{E}\left[\int_0^T\int_U \left(\rho^{n,\epsilon,K}-\bar\rho\right)^p\right]+\frac{1}{p(p-1)}\sup_{t\in[0,T]}\mathbb{E}\left[\int_U \left(\rho^{n,\epsilon,K}-\bar\rho\right)^p\right]\\
    \leq c \epsilon T\left(\|\nabla\cdot F^K_2\|_{L^\infty(U)}+\|F^K_3\|_{L^\infty(U)}\right),
\end{multline*}
and the $p$-dependent estimate
\begin{multline*}
    \frac{1}{Tp(p-1)}\mathbb{E}\left[\int_0^T\int_U \left(\rho^{n,\epsilon,K}-\bar\rho\right)^p\right]+\frac{1}{p(p-1)}\sup_{t\in[0,T]}\mathbb{E}\left[\int_U \left(\rho^{n,\epsilon,K}-\bar\rho\right)^p\right]\\
    \leq  c \epsilon^{p/2} T\left(\|\nabla\cdot F^K_2\|^{p/2}_{L^\infty(U)}+\|F^K_3\|^{p/2}_{L^\infty(U)}\right)\left(1+\epsilon T\left(\|\nabla\cdot F^K_2\|_{L^\infty(U)}+\|F^K_3\|
    _{L^\infty(U)}\right)\right),
\end{multline*}
where $F_i^K$, $i=1,2,3$ relate to the spatial component of the noise $\xi^K$, see Definition \ref{definition truncated noise}.
\end{proposition}
To quantify the rate at which $\rho^{\epsilon,K}$ converges to $\bar\rho$, we define
\begin{equation}\label{eq: for v epsilon k}
    v^{\epsilon,K}:=\epsilon^{-1/2}(\rho^{\epsilon,K}-\bar\rho).
\end{equation}
We will show, under an appropriate joint scaling regime $\epsilon\to 0, K\to\infty$, that solutions $v^{\epsilon,K}$ converge in probability in a space of distributions to the linearised\footnote{The equation is linear in $v$.} stochastic PDE
\begin{equation}\label{eq: linearised spde v}
   \begin{cases}
    \partial_t v=\Delta(\Phi'(\bar\rho)v))-\nabla\cdot(\sigma(\bar\rho)\dot\xi+\nu'(\bar\rho)v), & \text{on} \hspace{5pt} U\times(0,T],\\
v=0,&\text{on}\hspace{5pt} \partial U\times[0,T],\\
    v(\cdot,t=0)=0 ,&\text{on}\hspace{5pt} U\times\{t=0\},
   \end{cases} 
\end{equation}
where $\dot\xi$ is space-time white noise and $\bar\rho$ is the solution to the hydrodynamic limit equation \eqref{heat equation hydrodynamic limit}.
Due to the space-time white noise, solutions to equation \eqref{eq: linearised spde v} live in the space of distributions, we will show they live in $H^{-s}(U)$ for every $s>d/2$, see Proposition \ref{proposition existence of strong solutions to OU process} below.
The linearisation present in \eqref{eq: linearised spde v} can be formally seen by observing that $v^{\epsilon,K}$ solves
\begin{align}\label{formal computation for convergence of v epsilon to v}
    \partial_t v^{\epsilon,K}&=\Delta(\epsilon^{-1/2}(\Phi(\rho^{\epsilon,K})-\Phi(\bar\rho))-\nabla\cdot(\epsilon^{-1/2}(\nu(\rho^{\epsilon,K})-\nu(\bar\rho))-\nabla\cdot \left(\sigma(\rho^{\epsilon,K})\circ \dot\xi^{K}\right)\nonumber\\
    &=\Delta\left(v^\epsilon \frac{\Phi(\rho^{\epsilon,K})-\Phi(\bar\rho)}{\rho^{\epsilon,K}-\bar\rho}\right)-\nabla\cdot\left(v^\epsilon \frac{\nu(\rho^{\epsilon,K})-\nu(\bar\rho)}{\rho^{\epsilon,K}-\bar\rho}\right)-\nabla\cdot\left(\sigma(\rho^{\epsilon,K})\circ \dot\xi^{K}\right),
\end{align}
with zero initial condition and boundary data.
Owing to the law of large numbers result of Proposition \ref{ppn: energy estimates in introducition}, $\rho^{\epsilon,K}\to \bar\rho$ and also $\dot\xi^K\to\dot\xi$, which gives the formal convergence of $v^{\epsilon,K}$ to $v$. 
One of the main results of the present work will be to make the above ansatz rigorous.
As mentioned above, solutions $\rho^{\epsilon,K}$ to the singular equation \eqref{generalised Dean-Kawasaki Equation Stratonovich with epsilon} can only be considered in a stochastic kinetic sense. 
However, the non-linear nature of this definition makes it incompatible with convergence in the space of distributions.
As a result, we begin by proving the below central limit theorem for the regularised equation.
\begin{theorem}[Theorem \ref{theorem CLT for approximating equation} below]\label{thm: clt for regularised equation in intro}
    Fix $\alpha\in(0,1)$, $K\in\mathbb{N}$ and $\epsilon\in(0,1)$.
    Let $v^{n,\epsilon,K}$ denote the weak solution to the regularised equation \eqref{eq: definition of v n epsilon k with regularised equations}, and let $v$ be the solution of the regularised, linearised SPDE \eqref{eq: regularised OU process} in the sense of Definition \ref{definition strong solution to ou process}.
    Under Assumptions \ref{assumption on noise}, \ref{assumption of constant boundary data and constant initial condition}, \ref{asm: new assumptions} and \ref{asm: assumptions for well posedenss of equation},  there exists constants $\beta, c\in(0,\infty)$ independent of $n$, $\epsilon$, $K$ and $\alpha$ such that
     \begin{align*}
        \mathbb{E}&\|v^{n,\epsilon,K}-v\|^2_{L^2([0,T];H^{-s}(U))}\leq c\left(1+\epsilon T\left(\|\nabla\cdot F^K_2\|_{L^\infty(U)}+\|F^K_3\|
    _{L^\infty(U)}\right)\right)\nonumber\\
        &\times\left[\epsilon T^2\left(\|\nabla\cdot F^K_2\|^{2}_{L^\infty(U)}+\|F^K_3\|^{2}_{L^\infty(U)}\right)+\epsilon^{\beta+1} T^2\left(\|\nabla\cdot F^K_2\|^{\beta+2}_{L^\infty(U)}+\|F^K_3\|^{\beta+2}_{L^\infty(U)}\right)\right. \nonumber\\
        &\left.+\epsilon T\left(\|F_1^K\|^2_{L^\infty(U)}\|\sigma_n'\|^2_{L^\infty(U)}+\|F_2^K\|^2_{L^\infty(U;\mathbb{R}^d)}\right)\right] +c\mathbb{E}\|\sigma_n(\rho^{n,\epsilon,K})-\sigma(\bar\rho)\|_{L^2(U\times[0,T])}+cT\,\mathcal{T}_s(K),
    \end{align*}
    where $F_i^K$, $i=1,2,3$ relate to the spatial component of the noise $\xi^K$ and $\mathcal{T}_s(K)$ denotes a tail sum of norms of spatial components of the noise, defined in Theorem \ref{theorem CLT for approximating equation}, that decays to zero as $K\to\infty$.
\end{theorem}
The proof of Theorem \ref{theorem CLT for approximating equation} below illustrates how the ansatz \eqref{formal computation for convergence of v epsilon to v} is made rigorous.
Furthermore, in the above bound, the presence of $\|\sigma'_n\|^2_{L^\infty(U)}$ prevents us from directly taking the regularisation limit $n\to \infty$.
In the case of the singular square root diffusion, this term diverges when the solution approaches its zero set.
To get around this, one needs to prove the following $L^\infty(U\times[0,T])-$estimate.
\begin{theorem}[Theorem \ref{thm: l infinity estimate for singular equation} below]
Let $\epsilon\in(0,1)$, $K\in\mathbb{N}$ and under Assumption \ref{assumption of constant boundary data and constant initial condition}, let $\rho_0\equiv M>0$.
Assume also that Assumptions \ref{assumption on noise}, \ref{asm: new assumptions} and \ref{asm: assumptions for well posedenss of equation} hold.
Then, for stochastic kinetic solutions $\rho^{\epsilon,K}$ of \eqref{generalised Dean-Kawasaki Equation Stratonovich with epsilon}, there exist constants $c,\gamma\in(0,\infty)$ independent of $\epsilon$ such that
\[\mathbb{E}\|(\rho^{\epsilon,K}-M)_-\|_{L^\infty(U\times[0,T])}\leq c\left(\epsilon\left(\|F_3^K\|_{L^\infty(U)}+\|\nabla\cdot F_2^K\|_{L^\infty(U)}\right)\right)^\gamma. \]
\end{theorem}
Formally the above estimate says that if our initial condition satisfies $\rho_0\geq M$ for some $M>0$, then along appropriate scaling limits the solutions of equation \eqref{generalised Dean-Kawasaki Equation Stratonovich with epsilon} satisfy $\rho^{\epsilon,K}>M/2$ uniformly in space and time with high probability.
Using this and the pathwise uniqueness of stochastic kinetic solutions (Theorem 3.6 of \cite{popat2025well}), in Theorem \ref{thm: CLT for singular equation} we are able to prove a central limit theorem in probability for the singular equation.
Based on the above estimates, in Remark \ref{rmk: choice of joint scaling} we identify that an appropriate joint scaling regime is one which ensures that as $\epsilon\to0, K\to\infty$,
 \begin{equation}\label{eq: joint scaling intro}
     \epsilon\left(\|F_1^K\|^2_{L^\infty(U)}+\|F_2^K\|^2_{L^\infty(U;\mathbb{R}^d)}+\|F_3^K\|^2_{L^\infty(U)}+\|\nabla\cdot F_2^K\|^2_{L^\infty(U)}\right)\to0.
 \end{equation}
 In the case that we consider noise as in Example \ref{example of noise: eigenvalues of laplacian}, the scaling limit is explicit
 \[\epsilon K^{d+2}\to0.\]
Hence we are able to show that solutions $\rho^{\epsilon,K}$ of equation \eqref{generalised Dean-Kawasaki Equation Stratonovich with epsilon} converge to solutions $\bar\rho$  of the hydrodynamic limit equation \eqref{heat equation hydrodynamic limit}, and the central limit theorem helps us to identify the rate of convergence to first order.
Even though we expect the convergence of $\rho^{\epsilon,K}$ to $\bar\rho$, we are interested in quantifying the (un)likelihood that we see a profile very different to $\bar\rho$ after a long time. 
Whilst unlikely, the occurrence of a rare event can have catastrophic consequences, for instance mechanical failure of a machine or a population going extinct, see Section \ref{sec: motivation of boundary data} for further examples which are relevant on the bounded domain.
The proof of the uniform large deviations principle is based on the variational representation of infinite dimensional Brownian motion presented in \cite{budhiraja2008large}.
It involves the study of the below SPDE, which is equation \eqref{generalised Dean-Kawasaki Equation Stratonovich with epsilon} with an additional vector valued control term $g\in L^2(U\times[0,T];\mathbb{R}^d)$,
\begin{equation}\label{eq: controlled SPDE intro}
    \partial_t\rho^{\epsilon,g}=\Delta\Phi(\rho^{\epsilon,g})-\sqrt{\epsilon}\nabla\cdot(\sigma(\rho^{\epsilon,g})\circ\dot{\xi}^K)-\nabla\cdot(\sigma(\rho^{\epsilon,g})P_K g)-\nabla\cdot \nu(\rho^{\epsilon,g}),
\end{equation}
where $P_K g$ denotes the projection of the control onto the first $K$ modes of the noise.
The proof relies on three key ingredients.
The first is for every $\epsilon\in(0,1)$ and boundary condition $\bar{f}\in H^1(\partial U)$, the existence of a measurable solution map $\mathcal{G}^{\bar{f},\epsilon}: Ent_\Phi(U)\times C([0,T];(\mathbb{R}^d)^\infty)\to L^1(U\times[0,T])$ that takes the initial data and noise to the (stochastic kinetic) solution of equation  \eqref{eq: controlled SPDE intro}, and the existence of a similar map $\mathcal{G}^{\bar{f},0}$ for the limiting (zero noise) parabolic-hyperbolic PDE.
\begin{proposition}[Propositions \ref{prop: existence of solution map for skeleton equation},\ref{prop: existence of solution map for controlled SPDE} below]
    Suppose Assumptions \ref{asm: assumptions for well posedenss of equation} and \ref{asm: assumptions on boundary data for well-posedness} are satisfied, and fix the boundary data $\bar{f}\in H^1(\partial U)$.
    Then for every $\epsilon\in(0,1)$, initial condition $\rho_0\in Ent_\Phi(U)$, predictable control $g\in L^2(U\times[0,T];\mathbb{R}^d)$ and $g_k(s):=\langle g(\cdot,s),f_k\rangle_{L^2(U)}$, there exists measurable maps
    \[\mathcal{G}^{\bar{f},\epsilon},\mathcal{G}^{\bar{f},0}:Ent_\Phi(U)\times C([0,T];(\mathbb{R}^d)^\infty)\to L^1(U\times[0,T])\]
    such that
    \[\rho^{\epsilon,g}(\bar{f},\rho_0)=\mathcal{G}^{\bar{f},\epsilon}\left(\rho_0,\sqrt{\epsilon}(B_k)_{k\in\mathbb{N}}+\left(\int_0^\cdot g_k^\epsilon(s)\,ds\right)_{k\in\mathbb{N}}\right), \,\,\rho(\bar{f},\rho_0)=\mathcal{G}^{\bar{f},0}\left(\rho_0,\left(\int_0^\cdot g_k^\epsilon(s)\,ds\right)_{k\in\mathbb{N}}\right),\]
    where $\rho^{\epsilon,g}(\bar{f},\rho_0)$ denotes the unique stochastic solution of \eqref{eq: controlled SPDE intro} with initial data $\rho_0$ and boundary data $\Phi(\rho^{\epsilon,g})|_{\partial U}=\bar{f}$ in the sense of Definition 2.8 of \cite{popat2025well} and $\rho(\bar{f},\rho_0)$ denotes the unique weak solution of \eqref{eq: skeleton equation} with initial data $\rho_0$ and boundary data $\Phi(\rho)|_{\partial U}=\bar{f}$ in the sense of Definition \ref{def: weak solution of skeleton equation} below.
\end{proposition}
If we had the existence of a measurable map that also took the boundary data as an input, we would have been able to extend the large deviations principle to uniform subsets of bounded initial data and boundary data.
Through careful reasoning, in Section \ref{sec: lack of uniformity of solution map with respect to initial condition} we illustrate why such a map does not exist, so the large deviations principle is only uniform with respect to the initial data and not both the initial and boundary data.\\
The second condition is the collapse of the controlled SPDE \eqref{eq: controlled SPDE intro} to the limiting parabolic-hyperbolic PDE under the joint scaling \eqref{eq: joint scaling intro}.
\begin{theorem}[Theorem \ref{thm: convergence of solutions when initial data and controls converge} below]
   Let the coefficients of \eqref{eq: regularised controlled SPDE} satisfy Assumptions \ref{assumption on noise}, \ref{asm: assumptions for well posedenss of equation}, \ref{asm: assumptions on boundary data for well-posedness} and inequality \eqref{eq: assumption on quotient Phi / Phi'}. 
    Let $\{g^\epsilon\}_{\epsilon\in(0,1)}$ and $g$ be uniformly bounded $L^2(\Omega; L^2(U\times[0,T];\mathbb{R}^d))$-valued, predictable processes satisfying $g^\epsilon\to g$ weakly in $L^2(U\times[0,T];\mathbb{R}^d)$ as $\epsilon\to0$ and let $\{\rho_0^\epsilon\}_{\epsilon\in(0,1)},\rho_0\in Ent_\Phi(U)$ be uniformly bounded and such that 
    $\rho_0^\epsilon\to\rho_0$ weakly in $L^1(U)$ as $\epsilon\to0$.
    Let $\{\epsilon, K(\epsilon)\}_{\epsilon\in(0,1)}$ be such that the joint scaling \eqref{eq: joint scaling intro} is satisfied. 
    For every $\epsilon\in(0,1)$ and fixed boundary data $\bar{f}\in H^1(\partial U)$, denote by $\rho^{\epsilon,g^\epsilon}(\bar{f},\rho^\epsilon_0)$ the stochastic kinetic solutions of the controlled SPDE \eqref{eq: controlled SPDE intro} with control $g^\epsilon$, boundary data $\bar{f}$, and initial data $\rho_0^\epsilon$.
    Then we have that
    \[\rho^{\epsilon,g^\epsilon}(\bar{f},\rho_0^\epsilon)\to\rho(\bar{f},\rho_0) \quad\text{weakly in}\quad L^1(U\times[0,T]) \quad \text{as} \quad \epsilon\to0,\]
    where recall that we used $\rho$ to denote the unique weak solution of the limiting PDE \eqref{eq: skeleton equation}.
\end{theorem}
The third main condition for the large deviations is the lower-semicontinuity of the proposed rate function, proved in Proposition \ref{prop: rate function is lower semicontinuous}, which follows from the weak lower-semicontinuity of the $L^2(U\times[0,T];\mathbb{R}^d)-$norm.
This allows us to prove the large deviations principle.
\begin{theorem}[Theorem \ref{thm: LDP theorem from BDM} below]\label{thm: intro ldp theorem}
     Let $\bar{f}\in H^1(\partial U)$ and  $\rho_0\in Ent_\Phi(U)$ denote the boundary data and initial condition respectively.
     For arbitrary $\rho\in L^1(U\times[0,T])$, define 
    \begin{multline*}
        I_{\bar{f},\rho_0}(\rho):=\frac{1}{2}\inf_{g\in L^2(U\times[0,T];\mathbb{R}^d)}\left\{\|g\|^2_{L^2(U\times[0,T];\mathbb{R}^d)}:\partial_t\rho=\Delta\Phi(\rho)-\nabla\cdot(\sigma(\rho)g+ \nu(\rho))\right.\nonumber\\
        \left.: \Phi(\rho)|_{\partial U}=\bar{f}\,\, \text{and} \,\, \rho(\cdot,0)=\rho_0\right\}.
    \end{multline*}
    Then under Assumptions \eqref{assumption on noise}, \ref{asm: assumptions for well posedenss of equation} and \ref{asm: assumptions on boundary data for well-posedness}, for every $\rho\in L^1(U\times[0,T])$,
     $\rho\mapsto I_{\bar{f},\rho_0}(\rho)$ is a rate function on $L^1(U\times[0,T])$, the family $\{I_{\bar{f},\rho_0}(\cdot): \bar{f}\in H^1(\partial U),\, \rho_0\in Ent_\Phi(U)\}$ of rate functions has compact level sets on compacts, and solutions $\{\rho^{\epsilon,K(\epsilon)}\}_{\epsilon\in(0,1)}$ of \eqref{generalised Dean-Kawasaki Equation Stratonovich with epsilon} satisfy the large deviations principle on $L^1(U\times[0,T])$ with rate function $I_{\bar{f},\rho_0}$.
\end{theorem}

\subsection{Application to particle systems}\label{sec: application to particle system}
To motivate how equations such as \eqref{generalised Dean-Kawasaki Equation Stratonovich with epsilon} arise in the study of fluctuating hydrodynamics, briefly consider the concrete example of the symmetric zero range process $\eta^N$ on the torus $\mathbb{T}^d_N:=(\mathbb{Z}^d\setminus N\mathbb{Z}^d)$.
 That is, the process satisfying that if there are $k$ particles at a site $x\in \mathbb{T}_N^d$, then wait a random exponential time (with density proportional to $k\exp(-kt)$, for $t\in(0,\infty)$), then jump to a neighbour uniformly at random.
 Jumps at different sites $x,y\in\mathbb{T}_N^d$ occur independently. 
 Let 
 \[\mu^N_t:=\frac{1}{N^d}\sum_{x\in(\mathbb{Z}^d/N\mathbb{Z}^d)}\delta_{x/N}\cdot\eta^N_{N^2t}(x)\] denote the parabolically rescaled empirical measure for the particle system, which is a discrete measure that tracks the scaled position of particles after a scaled amount of time $t$. 
It was shown in \cite{ferrari1987local} that $\mu^N$ converges to the deterministic measure $\bar\rho\,dx$ where $\bar\rho:\mathbb{T}^d\times[0,T]\to\mathbb{R}$ is the unique solution to the heat equation $\partial_t\bar\rho=\frac{1}{2}\Delta\bar\rho$, in the sense that for every continuous $f:\mathbb{T}^d\times[0,T]\to\mathbb{R}$ and $\delta\in(0,1)$, we have
\[\lim_{N\to\infty}\mathbb{P}\left[|\langle f,\mu^N\rangle-\langle f,\bar\rho\rangle|>\delta\right]=0.\]
This formally illustrates that we expect particles to uniformly spread out on the torus, which can be seen by the fact that sites with a larger number of particles jump more frequently on average.\\
Suppose more generally that if there are $k$ particles at a site, we instead wait an exponential time with density $g(k)\exp(-g(k)t)$ for some non-decreasing function function $g:\mathbb{N}_0\to\mathbb{R}_{\geq 0}$ satisfying $g(0)=0$ and $g(k)>0$ for $k>0$.
In this case $\bar\rho$ instead solves $\partial_t\bar\rho=\frac{1}{2}\Delta\Phi(\bar\rho)$ for the mean local jump rate $\Phi$, see \cite{kipnis1998scaling}.
Furthermore, nonequilibrium central limit fluctuations was shown in Chapter 11 of \cite{kipnis1998scaling} and in \cite{dirr2016entropic}, where the measures
\[m^N:=N^{d/2}(\mu^N-\bar\rho\,dx)\]
are shown to converge in the space of distributions as $N\to\infty$ to the linearised stochastic PDE
\begin{equation}\label{eq: equation linearised DK}
    \partial_t m=\Delta(\Phi'(\bar\rho)m)-\nabla\cdot(\Phi^{1/2}(\bar\rho)\dot{\xi}),
\end{equation}
where $\dot\xi$ is $d$-dimensional space-time white noise and $\bar\rho$ is the solution to the heat equation $\partial_t\bar\rho=\frac{1}{2}\Delta\Phi(\bar\rho)$ above.
Simulating particle systems to observe central limit convergence and rare event behaviour is expensive when there are a large number of particles. 
Instead, if we instead consider the stochastic PDE
    \begin{equation*}
        \partial_t\rho^N=\Delta\Phi(\rho^N)-\frac{1}{\sqrt{N}}\nabla\cdot(\Phi^{1/2}(\rho^N)\circ\dot\xi^N)
    \end{equation*}
    for suitable approximations of space-time white noise $\dot\xi^N$, see Definition \ref{definition truncated noise} below, then the central limit theorem result of Theorem \ref{thm: CLT for singular equation} tells us that even on a bounded domain, as $N\to\infty$, fluctuations of the function\footnote{For the SPDE, instead of the scaling $\sqrt{N}$ we could instead look at the scaling $N^{d/2}$ as in the particle system.
    This would simply amount to a different scaling regime, with $\epsilon$ replaced by $\epsilon^d$ in equation \eqref{eq: joint scaling intro}.} $\sqrt{N}(\rho^N-\bar\rho)$ converge to the solution of equation \eqref{eq: equation linearised DK}.\\
    That is to say, solutions of equations such as \eqref{generalised Dean-Kawasaki Equation Stratonovich with epsilon} correctly describe the non-equilibrium fluctuations of particle systems.\\
     Let us now discuss the large deviations. 
     Based on the central limit theorem $m^N\to m$, a natural first guess is to consider the measure arising from the first order expansion  $\bar{\rho}^N:=\bar\rho\,dx+\frac{1}{\sqrt{N}}m$.
    Schilder's theorem and the contraction principle prove that $\bar{\rho}^N$ satisfies a large deviations principle with rate function
        \small
       \begin{align*}
        \bar{I}_{\rho_0}(\rho)&:=\frac{1}{2}\inf_{g\in L^2(U\times[0,T];\mathbb{R}^d)}\left\{\|g\|^2_{L^2(U\times[0,T];\mathbb{R}^d)}:\partial_t(\rho-\bar\rho)=\Delta\left(\Phi'(\bar{\rho})(\rho-\bar{\rho})\right)-\nabla\cdot(\Phi^{1/2}(\bar\rho)g): \rho(\cdot,0)=\rho_0\right\}.
    \end{align*}
        \normalsize
    However, it was shown in \cite{benois1995large} that the large deviations of the particle system $\mu^N$ are governed by the rate function
        \begin{align*}
        I_{\rho_0}(\rho)&:=\frac{1}{2}\inf_{g\in L^2(U\times[0,T];\mathbb{R}^d)}\left\{\|g\|^2_{L^2(U\times[0,T];\mathbb{R}^d)}:\partial_t\rho=\Delta\Phi(\rho)-\nabla\cdot(\Phi^{1/2}(\rho)g):\, \rho(\cdot,0)=\rho_0\right\}.
    \end{align*}
This does not coincide with the rate function $\bar{I}_{\rho_0}$ from the first order expansion.
On the other hand, in Theorem \ref{thm: LDP theorem from BDM} below we are able to show that even on a bounded domain, the truncated SPDE $\rho^N$ correctly captures the same large deviations rate function as the particle system $\mu^N$.

\subsection{Background literature}\label{sec: background literature}
 Well-posedness of equations such as \eqref{generalised Dean-Kawasaki Equation Stratonovich with epsilon} is a tricky problem in itself. 
Formal derivations of stochastic PDEs such as \eqref{generalised Dean-Kawasaki Equation Stratonovich with epsilon} were derived in physics literature in the nineties by Dean \cite{dean1996langevin} and Kawasaki \cite{kawasaki1994stochastic} when deriving an equation satisfied by the empirical measure of an interacting particle system with Langevin dynamics.
In their derivations, the SPDE appears with space-time white noise instead of truncated or correlated noise.
The irregularity of space-time white noise means that equations such as \eqref{generalised Dean-Kawasaki Equation Stratonovich with epsilon} with space-time white noise are supercritical in the language of regularity structures, so are not renormalisable using Hairer's regularity structures \cite{hairer2014theory} or Gubinelli, Imkeller and Perkowski’s paracontrolled distributions \cite{gubinelli2015paracontrolled}, even in dimension $d=1$.
Furthermore, Konarovskyi, Lehmann and von Renesse show in Theorem 2.2 of \cite{konarovskyi2019dean} and Theorem 1 of \cite{konarovskyi2020dean} that the only martingale solutions of such equations coincide with equations satisfied by empirical measures of particle systems, in the sense that they agree as semi-martingales.\\
With these negative results, attention in the community turned to proving well-posedness for regularised versions of the Dean--Kawasaki equation. 
We describe a non-exhaustive list of research in this direction, noting from Section \ref{sec: application to particle system} that a key case of interest is the case when $\Phi$ is the identity and so the diffusion coefficient is the critical square root $\sigma(\cdot)=\Phi^{1/2}(\cdot)=\sqrt{\cdot}$.\\
The works of Fehrman and Gess \cite{fehrman2019well,fehrman2021path} proved path-by-path well-posedness of equations such as \eqref{generalised Dean-Kawasaki Equation Stratonovich with epsilon} on the torus using techniques from rough path theory \cite{lyons1998differential}, motivated also by the theory of stochastic viscosity solutions, see the works of Lions and Souganidis \cite{lions1998fully1, lions2000fully}. 
The techniques were extended to a bounded domain with homogeneous Dirichlet boundary conditions by Clini \cite{clini2023porous}.
However, the works required the diffusion coefficient to be six times continuously differentiable, so fell far short of the square root diffusion.\\
A construction of probabilistic solutions via an entropy formulation was considered by Dareiotis and Gess \cite{dareiotis2020nonlinear}, based on the earlier work of Dareiotis, Gerencs\'er and Gess \cite{dareiotis2019entropy}.
However, the techniques require higher order integrability of the initial data than in the $L^1(U)$-theory of kinetic solutions, and require $C^{3}_b$-regularity of the diffusion coefficient, so again fall short of the square root diffusion.\\
The existence of probabilistic strong solutions in the case of truncated noise and mollified square root was provided by Djurdjevac, Kremp and Perkowski \cite{djurdjevac2022weak}.
The proof follows a variational approach for a transformed equation and energy estimates for the Galerkin projected system. 
A recent extension of the result was given by Djurdjevac, Ji and Perkowski \cite{djurdjevac2025weak}.\\
Finally, an approach proving local well-posedness of was proved via paracontrolled theory by Martini and Mayorcas \cite{martini2022additive}.
Motivated by the theory of linear fluctuating hydrodynamics, the authors replace the critical square root diffusion by the square root of the solution to the corresponding zero-noise PDE.\\
In all of the above, the square root diffusion was out of reach.
The first work which encapsulated the critical square root diffusion was proved on the torus by Fehrman and Gess \cite{fehrman2024well}, in which the space-time white noise was replaced by infinite dimensional Stratonovich noise which is white in time and sufficiently regular in space.
From the point of view of particle systems, the spatial correlation is justified due to the natural correlation length given by the grid size\footnote{If particles evolve on a lattice with grid size $\epsilon>0$, then we would expect the noise to be constant on blocks of size $\epsilon/2$, not uncorrelated.}.
The authors consider well-posedness of stochastic kinetic solutions, a solution theory first introduced in the partial differential equation setting by Lions, Perthame and Tadmor in the nineties \cite{lions1994kinetic} when studying multidimensional scalar conservation laws.
Relevant extensions to degenerate parabolic-hyperbolic PDEs were given by Chen and Pertame \cite{chen2003well}, Bendahmane and
Karlsen \cite{bendahmane2004renormalized}, and Karlsen and Riseboro \cite{karlsen2000uniqueness}.
In the critical square root diffusion case, kinetic solution theory is necessary to deal with singular terms that appear in the It\^o-to-Stratonovich conversion of \eqref{generalised Dean-Kawasaki Equation Stratonovich with epsilon}, see \eqref{eq: ito form of dk equation} below.
The solution theory is based on the equation's kinetic formulation which renormalises the solution away from it's zero set.
The derivation of the kinetic equation for the controlled equation \eqref{eq: controlled SPDE intro} is presented in the Appendix, see Section \ref{sec: kinetic formulation for controlled SPDE}.\\
The well-posedness in the same setting was extended to a $C^2$, bounded domain with Dirichlet boundary conditions by the author in \cite{popat2025well}, and the boundary data $\bar{f}$ considered encompassed all non-negative constant functions including zero and all smooth functions bounded away from zero.
This encompasses the well-posedness of equation \eqref{generalised Dean-Kawasaki Equation Stratonovich with epsilon} for every $\epsilon\in(0,1)$ and $K\in\mathbb{N}$ considered in the present work.
A recent extension by Fehrman and Gess \cite{fehrman2024conservative} prove well-posedness of the equation on the whole space.
Motivated by the asymmetric zero range process, Wang, Wu and Zhang \cite{wang2022dean} extend the well-posedness of stochastic kinetic solutions to equations incorporating a a non-local convolution
\begin{equation}\label{eq: dk with convolution}
    \partial_t\rho=\Delta\rho-\sqrt{\epsilon}\nabla\cdot(\sqrt{\rho}\circ \dot{\xi}^K)-\nabla\cdot(\rho V\ast \rho).
\end{equation}
where $V$ is an interaction kernel satisfying the Ladyzhenskaya-Prodi-Serrin condition \cite{ladyzhenskaya1967uniqueness,prodi1959teorema,serrin1961interior} and $\ast$ represents a spatial convolution.\\
Let us now discuss literature in the direction of a central limit theorem.
Central limit theorems for stochastic heat equations with Lipschitz continuous diffusion coefficients were studied by by Huang, Nualart and Viitasaari \cite{huang2020central}, and Chen, Khoshnevisan, Nualart and Pu \cite{chen2022central}.
A central limit theorem for the stochastic wave equation with Lipschitz continuous diffusion coefficients was obtained by  Delgado-Vences, Nualart and Zheng \cite{delgado2020central}.\\
A central limit theorem for the generalised Dean--Kawasaki equation on the torus was proved by Clini and Fehrman \cite{clini2023central}.
Due to the more technical estimates used in the the present work, we are able to handle the full range of fast diffusion and porus medium non-linearities $\Phi(\xi)=\xi^m$ for every $m\in(0,\infty)$.
In \cite{clini2023central} the authors can only handle the non-linearity for $m$ sufficiently small due to an assumed bound on the growth of $|\Phi''(\xi)|$, see Assumption 2.8(ii) of \cite{clini2023central}.\\
A central limit theorem for SPDEs capturing fluctuations of the symmetric simple exclusion process was given by Dirr, Fehrman and Gess \cite{dirr2020conservative}.
The authors study fluctuations of
\[\partial_t\rho^\epsilon=\Delta\rho^\epsilon-\sqrt{\epsilon}\nabla\cdot\left(\sqrt{\rho^\epsilon(1-\rho^\epsilon)}\circ\dot\xi^K\right).\]
The well-posedness of the above equation poses the additional difficulty that there are multiple points of singularity $\{\rho^\epsilon\approx 0\}$ and $\{\rho^\epsilon\approx 1\}$, which on the kinetic solution level requires a renormalisation of the solution away from both of these sets.
On the other hand, due to the fact that solutions are $[0,1]$-valued, estimates in \cite{dirr2020conservative} are often more simple than in the present work.\\
Large deviations for reaction diffusion equations with additive white in time, coloured in space noise was shown by Cerrai and Freidlin \cite{cerrai2011approximation}.
Large deviations for stochastic porus media equations have been shown using exponential estimates and a generalized contraction principle by R\"ockner, Wang, and Wu \cite{rockner2006large} and more recently using the weak approach to large deviations by Zhang \cite{zhang2020large}.
The latter work more closely resembles the present work as kinetic solutions are considered, and furthermore the approach to large deviations is also based on the weak approach of Budhiraja, Dupuis, and Maroulas \cite{budhiraja2008large}.\\
The weak approach to large deviations has also developed in the context of Markov chains and random walk models by Dupuis and Ellis \cite{dupuis2011weak}, and for several systems including infinite dimensional Brownian noise in the more recent book by Budhiraja and Dupuis \cite{budhiraja2019analysis}.
Motivated by the same works, Fehrman and Gess \cite{fehrman2023non} prove a large deviations principle for the Dean--Kawasaki equation
\[\partial_t\rho=\Delta\Phi(\rho)-\sqrt{\epsilon}\nabla\cdot(\Phi^{1/2}(\rho)\circ \dot{\xi}^K)\] on the torus, capturing the large deviation behaviour of the zero range process.
In follow up work, Wu and Zhang \cite{wu2022large} prove the large deviations principle for the equation \eqref{eq: dk with convolution} with non-local transport, again on the torus.
For the simple exclusion process the corresponding large deviation results was given by Dirr, Fehrman and Gess \cite{dirr2020conservative}.\\
The weak approach for large deviations was also used to study the large deviations of reaction–diffusion equations with arbitrary polynomial non-linearity by Cerrai and Debussche \cite{cerrai2019large}, for stochastic Landau–Lifshitz equation on a bounded interval by Brze\'zniak, Goldys, and Jegaraj \cite{brzezniak2017large}, and for first-order scalar conservation laws perturbed by small multiplicative noise by Dong, Wu, Zhang and Zhang \cite{dong2020large}.\\
Compared with the above works on the torus, there are several technical points on the bounded domain that need to be addressed.
For example, from Definition \ref{definition truncated noise} and Assumption \ref{assumption on noise} below that we need to define the noise in a precise way to ensure that it has the desired behaviour at the boundary.
Furthermore, energy estimates are significantly complicated due to the presence of boundary data and the lack of conservation of mass.
Since the domain $U\subset \mathbb{R}^d$ is quite general, we often have to use techniques from PDE theory such as the Sobolev extension theorem (see Chapter 5.4 of \cite{evans2022partial}) to prove higher order spatial regularity, see Remark \ref{rmk: extension operator}.

\subsection{Motivation of boundary data}\label{sec: motivation of boundary data}
We noted previously that a similar result in a less general framework was proved on the torus in \cite{clini2023central}. 
For mathematicians and physicists, the torus is often a good starting point to prove properties of an equation due to its nice properties like boundedness, connectedness and periodic boundary conditions.
In practice, when particle systems are used to model real life phenomena, it is more natural and appropriate to consider the system as evolving on a bounded domain. 
In this case, homogeneous Neumann boundary conditions are used to model reflection of particles at the boundary, and Dirichlet boundary conditions model absorption or injection of particles.\\
For the central limit theorem results, in Assumption \ref{assumption of constant boundary data and constant initial condition} below we assume that the Dirichlet boundary data is a random positive constant.
If the density of particles around the boundary is lower than the boundary value, then the boundary data corresponds to injection of particles.
On the other hand, if the density of particles around the boundary is higher than the boundary value, then the boundary condition corresponds to absorption of particles at the boundary. 
Below we give several examples motivating physical systems that can be modelled by particle systems with non-zero Dirichlet boundary conditions.
\begin{itemize}
    \item The most basic example would be to consider the heat equation ($\sigma=\nu=0$ in \eqref{generalised Dean-Kawasaki Equation Stratonovich with epsilon}) which models the heat of a liquid in a container.
    The Dirichlet boundary condition considered here corresponds to external fixed heat (thermal reservoir) being applied to the system at the boundary of the container.
    \item The symmetric zero range process can be used to model vehicle traffic jams or pedestrian movement, see \cite{kaupuvzs2005zero}, where individuals are moving between different road sections.
    The boundary data corresponds to individuals leaving or entering the area of interest.
    \item In economics, the zero range process can be used to model wealth flow, see \cite{pinasco2018game,cardoso2021wealth}. The sites represent individuals and the particles represent wealth moving from one agent to another.
    The boundary data can model wealth entering or exiting the system, though most wealth models assume conservation of wealth.
      \item In biochemical reactions, the boundary condition corresponds to a steady influx of reactants being added at the boundary, corresponding to a fixed concentration of a certain reactant there.
    \item In epidemiology and population modelling, when studying the evolution of a population in a particular space such as island or in a country, the non-zero Dirichlet boundary condition models immigration and emigration of individuals into or out of a specific area.\\
    In the context of modelling disease spread, the boundary condition models new infected individuals entering the population, and already infected individuals leaving the population. 
\end{itemize}
Rather than just strictly positive boundary data, is natural to wonder whether it is possible to consider the homogeneous Dirichlet boundary data $\Phi(\rho^{\epsilon,K})|_{\partial U}=0$ in \eqref{generalised Dean-Kawasaki Equation Stratonovich with epsilon}, which corresponds to the case of pure absorption of particles at the boundary.
We recall that the solutions $\rho^{\epsilon,K}$ of \eqref{generalised Dean-Kawasaki Equation Stratonovich with epsilon} have a singularity on the zero set of the solution $\{\rho^{\epsilon,K}\approx 0\}$ due to the It\^o-to-Stratonovich conversion of \eqref{generalised Dean-Kawasaki Equation Stratonovich with epsilon}.
This is overcome in the $L^\infty(U\times[0,T])-$estimate of Theorem \ref{thm: l infinity estimate for singular equation} below, that formally says that if the solution starts at least distance $M>0$ away from zero, then it remains at least $M/2$ distance away from zero at all times $t\in[0,T]$ and points in $U$ with high probability. 
Such an estimate would not hold if we considered homogeneous Dirichlet boundary data $\Phi(\rho^{\epsilon,K})|_{\partial U}=0$.
Indeed, in the fast diffusion case $\Phi(\xi)=\xi^m$ for $m\in(0,1)$, Section 5 of \cite{bonforte2024cauchy} gives that the solution extinguishes in finite time almost surely, even in the much more simple regime $\sigma=\nu=0$.
\subsection{Organisation of the paper}
In Section \ref{sec: setup} we set up the problem.
In Definition \ref{definition truncated noise} we define the truncated noise considered in equation \eqref{generalised Dean-Kawasaki Equation Stratonovich with epsilon}, and state the various assumptions that are required for the spatial component of the noise.
An important example of the noise coefficients is given in Example \ref{example of noise: eigenvalues of laplacian}.
A simplifying assumption about the boundary data and initial condition of equation \eqref{generalised Dean-Kawasaki Equation Stratonovich with epsilon} is then given in Assumption \ref{assumption of constant boundary data and constant initial condition}.
The assumption is required for the central limit theorem but is not needed for the large deviations principle.\\
In Section \ref{sec: quantitative law of large numbers}, in Definition \ref{defn: regularised equation} we define a regularised version of \eqref{generalised Dean-Kawasaki Equation Stratonovich with epsilon} for which we can define weak solutions.
Quantitative law of large number estimates are given in Proposition \ref{ppn: energy estimate for difference of regularised equation and hydrodynamic limit} and \ref{prop: second energy estimate for rho n epsilon k}, providing $p$-independent and $p$-dependent estimates for $L^p(U\times[0,T])$ and $L^\infty([0,T];L^p(U))$-norms of the difference $(\rho^{n,\epsilon,K}-\bar\rho)$ for $p\geq 2$.
After defining strong solutions to the linearised stochastic PDE $v$ as in \eqref{eq: linearised spde v}, in Definition \ref{definition strong solution to ou process}, we prove a central limit theorem in for the regularised equation in Theorem \ref{theorem CLT for approximating equation}, quantifying convergence of $v^{n,\epsilon,K}:=\epsilon^{-1/2}(\rho^{n,\epsilon,K}-\bar\rho)$ to $v$ in the space of distributions.\\
In Section \ref{sec: clt for singular equation} we aim to extend the central limit theorem to the singular equation.
The result is based on the $L^\infty(U\times(0,T))$-estimate of Theorem \ref{thm: l infinity estimate for singular equation}.
The central limit theorem in probability is then given in Theorem \ref{thm: CLT for singular equation}.
Finally, Remark \ref{rmk: choice of joint scaling} provides an explicit example of a joint scaling $\epsilon\to0, K \to\infty$ under which the results of this work hold.\\
Section \ref{sec: ldp} is devoted to the proof of the large deviations principle.
In Assumption \ref{asm: Assumption for large deviations principle} and Theorem \ref{thm: LDP theorem from BDM} we state sufficient conditions for a uniform large deviations principle to hold.\\
In Section \ref{sec: Weak solutions of the skeleton equation, existence of solution map for controlled SPDE} we begin with Definition \ref{def: weak solution of skeleton equation} that defines what it means to be a weak solution to the limiting parabolic-hyperbolic PDE that appears in the rate function.
We then prove the well-posedness of weak solutions of the limiting PDE, and equivalence of weak and kinetic solutions for the PDE, see Theorem \ref{thm: well-posedness of skeleton equation assumed}.
We conclude the section  with Propositions \ref{prop: existence of solution map for skeleton equation} and \ref{prop: existence of solution map for controlled SPDE}, where we prove the existence of solution maps for the limiting PDE and the controlled SPDE which are needed for Theorem \ref{thm: LDP theorem from BDM}.\\
The uniformity of the large deviations principle is with respect to bounded subsets of the initial data, and in Section \ref{sec: lack of uniformity of solution map with respect to initial condition} we illustrate why the uniformity can not be extended to bounded subsets of the boundary data.\\
In Section \ref{sec: Energy estimates for the regularised controlled SPDE}, in Propositions \ref{prop: l2 estimate from 4.14 of popat}, \ref{prop: higher order spatial regularity for regularised controlled SPDE} and \ref{prop: higher order time regularity for regularised controlled SPDE} we prove various estimates for a regularised version of the controlled SPDE. 
This enables us in  Corollary \ref{corr: tightness of laws of controlled SPDE} to prove tightness of the laws of the controlled SPDE.\\
In Section \ref{sec: checking conditions for BDM} we prove the two remaining conditions required to apply the large deviations Theorem \ref{thm: LDP theorem from BDM}.
In Theorem \ref{thm: convergence of solutions when initial data and controls converge} we first prove that solutions of the controlled SPDE weakly converge to solutions of the limiting PDE under the joint scaling of Remark \ref{rmk: choice of joint scaling}, if we have weak convergence of the initial data and the controls. 
The second condition, showing that the rate function is lower semi-continuous is given in Proposition \ref{prop: rate function is lower semicontinuous}.\\
In the Appendix, in Section \ref{subsec: Appendix example of noise} we give properties of one possible choice of noise coefficients, the eigenfunctions of the Laplacian with homogeneous Dirichlet boundary conditions. 
In Section \ref{sec: assumptions for well-posedness} we give some underlying Assumptions from \cite{popat2025well} under which well-posedness of stochastic kinetic solutions of \eqref{generalised Dean-Kawasaki Equation Stratonovich with epsilon} holds.
The assumptions are split into assumptions for the coefficients (Assumption \ref{asm: assumptions for well posedenss of equation}) and Assumptions on the boundary data (Assumption \ref{asm: assumptions on boundary data for well-posedness}).
In Section \ref{sec: Correct space for initial data for the large deviations principle} we give a formal estimate for the limiting PDE that illustrates that the correct space for the initial data is the space of functions with finite entropy.
We conclude in Section \ref{sec: kinetic formulation for controlled SPDE} by providing a derivation the kinetic formulation for the controlled SPDE.

\section{Setup}\label{sec: setup}
To begin we define the truncated noise considered in equation \eqref{generalised Dean-Kawasaki Equation Stratonovich with epsilon}.
\begin{definition}[Truncated noise]\label{definition truncated noise}
  Let $(\Omega,\mathcal{F},(\mathcal{F}_t)_{t\geq0},\mathbb{P})$ be a filtered probability space with adapted, independent $d$-dimensional Brownian motions $(B^k)_{k\in\mathbb{N}}$ with values in the space $C([0,\infty);(\mathbb{R}^d)^\infty)$ equipped with the metric topology of co-ordinate wise convergence.
  Let $\{f_k\}_{k\in\mathbb{N}}$ be a sequence of continuously differentiable real valued functions on the domain $U\subset \mathbb{R}^d$.
  For every $K\in\mathbb{N}$, define the truncated noise
\[\xi^K: U\times[0,T]\to\mathbb{R}^d,\hspace{20pt}\xi^K(x,t):=\sum_{k=1}^K f_k(x) B^k_t.\]
\end{definition}
For convenience, for every $K\in\mathbb{N}$ define the three quantities related to the spatial component of the noise,
\[ F^K_1:U\to\mathbb{R}
\hspace{10pt} \text{defined by} \hspace{10pt}F^K_1(x):=\sum_{k=1}^K f_k^2(x);\]
\[ F^K_2:U\to\mathbb{R}^d
\hspace{10pt} \text{defined by} \hspace{10pt}F^K_2(x):=\sum_{k=1}^K f_k(x)\nabla f_k(x)=\frac{1}{2}\sum_{k=1}^K \nabla f_k^2(x);
\]
\[ F^K_3:U\to\mathbb{R}
\hspace{10pt} \text{defined by} \hspace{10pt}F^K_3(x):=\sum_{k=1}^K |\nabla f_k(x)|^2.\] 
We need the below assumption on the noise.
\begin{assumption}\label{assumption on noise}
    For every $K\in\mathbb{N}$ and $i=1,2,3$, we assume that
    \begin{enumerate}
        \item The sums $F_i^K$ are continuous on $U$, and also that $\nabla\cdot F_2^K$ is bounded in $U$.
        \item For every $s>d/2$ we have
    \[\sum_{k=1}^\infty \|f_k\|^2_{H^{-s}(U)}<\infty.\]
    Repeatedly, when we want to bound sums of $H^{-s+1}(U)$-norms of the noise coefficients, we will require $s>\frac{d+2}{2}$, which is equivalent to the above assumption.
    \item For every $k\in\mathbb{N}$, $f_k$ satisfy homogeneous Dirichlet boundary conditions, that is
    $f_k(\cdot)|_{\partial U}=0$.
    \item For $j,k\in\mathbb{N}$, if $j\neq k$, then $f_k$ and $f_j$ are orthogonal in $L^2(U)$.
    \item The noise is probabilistically stationary in critical $\sigma(\cdot)=\sqrt{\cdot}$ case. 
    That is to say that $F_1$ is constant on $U$, by which it follows that $\nabla\cdot F_2\geq0$ on $U$.
    \end{enumerate}
\end{assumption}

\begin{remark}\label{rmk: limit of noise coefficients is unbounded}
  Since the noise coefficients $f_k$ are assumed to be orthogonal in $L^2(U)$, we do not expect the limits $\lim_{K\to\infty} F^K_i$ to exist pointwise for $i=1,2,3$.
    For example, on the torus $U=\mathbb{T}^d$ (cf. \cite[Asm.~3.1]{dirr2020conservative}), define the spectral approximation of $\mathbb{R}^d$-valued space-time white noise by 
    \[\sum_{k\in\mathbb{Z}^d: |k|\leq K}\sqrt{2}(sin(k\cdot x) dB^k_t+cos(k\cdot x)dW_t^k).\]
    A direct computation shows that $F^K_1=\#\{k\in\mathbb{Z}^d:|k|\leq K\}$ and $F^K_3=\sum_{|k|\leq K}|k|^2$, neither of which are bounded in the $K\to\infty$ limit.
       \end{remark}
  \begin{remark}\label{rmk: bound for divergence of F2}
    By direct computation we have $\nabla\cdot F_2^K= F_3^K+\sum_{k=1}^K f_k\Delta f_k$.
Therefore, by the triangle inequality and the Cauchy-Schwarz inequality, it holds that
\[\|\nabla\cdot F_2^K\|_{L^\infty(U)}\leq \|F_3^K\|_{L^\infty(U)}+\|F_1^K\|_{L^\infty(U)}^{1/2}\left\|\sum_{k=1}^K(\Delta f_k)^2\right\|^{1/2}_{L^\infty(U)}.\]
Hence in the first point of Assumption \ref{assumption on noise}, the boundedness of $\nabla\cdot F_2^K$ is essentially an assumption on boundedness of the sum 
\[\sum_{k=1}^K(\Delta f_k)^2.\]
By equation \eqref{eq: bound for divergence of F2 in laplacian eigenfunction example} in the Appendix, we see that for the choice of noise in Example \ref{example of noise: eigenvalues of laplacian} below, $\nabla\cdot F_2^K$ can be bounded above by $F_3^K$ for every $K\in\mathbb{N}$.
  \end{remark}     

\begin{remark}
    For every $K\in\mathbb{N}$, $F_2^K$ is a vector in $\mathbb{R}^d$. 
    By the fact that all norms on $\mathbb{R}^d$ are equivalent, there is no ambiguity when we write the norm  $\|F_2^K\|_{L^\infty(U;\mathbb{R}^d)}$ in the below estimates.
\end{remark}
The canonical example on a bounded domain that one should have in mind are the eigenfunctions of the Laplacian with homogeneous Dirichlet boundary conditions.
\begin{example}[Eigenvalues of Laplacian]\label{example of noise: eigenvalues of laplacian}

For $k\in\mathbb{N}$, set $f_k=e_k$, where $\{e_k\}_{k\in\mathbb{N}}$ are the homogeneous Dirichlet eigenfunctions of the Laplacian on $U$ with corresponding eigenvalues $\{\lambda_k\}_{k\in\mathbb{N}}$.
That is, solutions to the equation
\begin{equation}
    \begin{cases}\label{eq:eigenvector and eigenvalue of lap}
        -\Delta e_k=\lambda_k e_k,& x\in U\\
        e_k=0,& x\in\partial U.
    \end{cases}
\end{equation}
The noise is then defined for every $K\in\mathbb
N$ by 
\[\xi^K(x,t):=\sum_{k:\lambda_k\leq\lambda_K}e_kB_t^k.\]
We show in the Appendix, Section \ref{subsec: Appendix example of noise}, that this choice of noise satisfies the conditions in Assumption \ref{assumption on noise}.
\end{example} 
For the subsequent results, we will have to make the below simplifying assumption. 
It is similar to the assumption required on the torus, see Assumption 2.1 of \cite{clini2023central}.
\begin{assumption}[Simplifying Assumption]\label{assumption of constant boundary data and constant initial condition}
Suppose that the $\mathcal{F}_0$-measurable initial data and boundary data are the same random positive constant.
That is,
\[\Phi(\rho^{\epsilon,K})|_{\partial U} = \rho_0=M \quad \text{for some}\quad M>0.\]
\end{assumption}
\begin{remark}\label{rmk: remark on hydrodynamic limit solved by constant}
    The consequence of Assumption \ref{assumption of constant boundary data and constant initial condition} is that the solution $\bar\rho$ of the hydrodynamic limit equation \eqref{heat equation hydrodynamic limit} is solved uniquely\footnote{The uniqueness follows by the same arguments as Theorem 3.6 of \cite{popat2025well} in the special case $\sigma=0$ (no noise), and the equivalence of stochastic kinetic and weak solutions.} by the constant $\bar\rho=M$, and therefore we are essentially looking at fluctuations around a constant function.
   If the boundary data was different to the initial condition, provided that the transport term $\nu(\bar\rho)$ was sufficiently nice so that it is dominated by the diffusion $\Phi$, we would expect the solution to converge to the value of the boundary data.
   However, in this case we would not have that $\partial_t\bar\rho= 0$ for every $t\in[0,T]$, so the energy estimates below (see Proposition \ref{ppn: energy estimate for difference of regularised equation and hydrodynamic limit}) would not be possible. 
\end{remark}

\section{Quantitative law of large numbers and central limit }\label{sec: quantitative law of large numbers}
For every $\epsilon\in(0,1)$ and $K\in\mathbb{N}$, the  well-posedness of stochastic kinetic solutions of equation \eqref{generalised Dean-Kawasaki Equation Stratonovich with epsilon}
\[\partial_t\rho^{\epsilon,K}=\Delta\Phi(\rho^{\epsilon,K})-\sqrt{\epsilon}\nabla\cdot(\sigma(\rho^{\epsilon,K})\circ\dot{\xi}^K)-\nabla\cdot \nu(\rho^{\epsilon,K})\]
with Dirichlet boundary conditions is presented in \cite{popat2025well}. 
For completeness and clarity of presentation, we compile the list of necessary assumptions on the non-linearities in Assumption \ref{asm: assumptions for well posedenss of equation} and assumptions on the boundary data in Assumption \ref{asm: assumptions on boundary data for well-posedness} in Section \ref{sec: assumptions for well-posedness} of the Appendix. 
We emphasise that Assumption \ref{asm: assumptions for well posedenss of equation} captures all of the cases of interest, $\Phi(\xi)=\xi^m$ for the full range $m\in(0,\infty)$, the critical square root $\sigma(\xi)=\sqrt{\xi}$ case, and lipschitz $\nu$.
In order to consider weak solutions, we need to smooth the singularity $\sigma$, and to obtain energy estimates we need to add a regularisation.
\begin{definition}[Regularised equation]\label{defn: regularised equation}
  Let $\{\sigma_n\}_{n\in\mathbb{N}}$ be a sequence such that for each $n\in\mathbb{N}$, let $\sigma_n\in C([0,\infty))\cap C^\infty((0,\infty))$ with $\sigma_n(0)=0$ and $\sigma'_n\in C_c^\infty([0,\infty))$ as defined in Assumption 4.10 of \cite{popat2025well}.
  For $\alpha\in(0,1)$, let $\rho^{n,\epsilon,K}$ denote the solution to the regularised equation 
\begin{equation}\label{equation for rho n epsilon k}
    \partial_t\rho^{n,\epsilon,K}=\Delta\Phi(\rho^{n,\epsilon,K})+\alpha\Delta\rho^{n,\epsilon,K}-\sqrt{\epsilon}\nabla\cdot(\sigma_n(\rho^{n,\epsilon,K})\circ\dot{\xi}^K)-\nabla\cdot \nu(\rho^{n,\epsilon,K}),
\end{equation}
with boundary data $\Phi(\rho^{n,\epsilon,K})=\bar{f}$ and initial condition $\rho_0$.
\end{definition}
Lemma 5.22 of \cite{fehrman2024well} shows that we can pick $\{\sigma_n\}_{n\in\mathbb{N}}$ to approximate $\sigma$, in the sense that the sequence satisfies the regularity properties in Definition \ref{defn: regularised equation} uniformly in $n$ and satisfies $\sigma_n\to\sigma$ in $C^1_{loc}((0,\infty))$ as $n\to\infty$.
Furthermore, regularisation in $\alpha$ is a technical requirement and will not play a key role here, so for ease of notation we will omit it in the superscript throughout the text.\\
The It\^o-to-Stratonovich conversion of \eqref{equation for rho n epsilon k} is 
\begin{multline}\label{eq: ito form of dk equation}
    \partial_t\rho^{n,\epsilon,K}=\Delta\Phi(\rho^{n,\epsilon,K})+\alpha\Delta\rho^{n,\epsilon,K}-\sqrt{\epsilon}\nabla\cdot(\sigma_n(\rho^{n,\epsilon,K})\dot{\xi}^K)-\nabla\cdot \nu(\rho^{n,\epsilon,K})\\
    +\frac{\epsilon}{2}\nabla\cdot(F_1^K(\sigma_n'(\rho^{n,\epsilon,K})^2\nabla\rho^{n,\epsilon,K}+\sigma_n(\rho^{n,\epsilon,K})\sigma'_n(\rho^{n,\epsilon,K}) F^K_2).
\end{multline}
For this equation, we study the CLT about the solution to
 \begin{equation}\label{eq: regularised hydrodynamic limit}
 \partial_t\bar{\rho}=\Delta\Phi(\bar{\rho})+\alpha\Delta \bar{\rho}-\nabla\cdot \nu(\bar{\rho}),
 \end{equation}
 where we abused notation as initially $\bar\rho$ denoted the solution to the above equation in the case $\alpha=0$.
 However, owing to Assumption \ref{assumption of constant boundary data and constant initial condition}, the unique weak solution of \eqref{eq: regularised hydrodynamic limit} is still $\bar\rho=M$ for every $\alpha\in(0,1)$.\\
With the regularisations mentioned in Definition \ref{defn: regularised equation}, we can make sense of weak solutions to \eqref{equation for rho n epsilon k}.
\begin{definition}[Weak solution of regularised equation]\label{weak solution of regularised equation in CLT and LDP}
    A weak solution of equation \eqref{equation for rho n epsilon k} satisfies
    \begin{enumerate}
        \item The regularity and boundary condition: 
        For every $k\in\mathbb{N}$, \[\left((\rho^{n,\epsilon,K}\wedge k)\vee 1/k\right)-\left((\Phi^{-1}(\bar f)\wedge k)\vee 1/k\right)\in L^2(\Omega \times [0,T];H^1_0(U)).\]
        \item The equation:     For all test functions $\psi\in C_c^\infty(U)$ and every $t\in[0,T]$,
\begin{multline*}
    \int_U\rho^{n,\epsilon,K}(x,t)\psi(x)=\int_U\rho_0(x)\psi(x)-\int_0^t\int_U\Phi'(\rho^{n,\epsilon,K})\nabla\psi\cdot\nabla\rho^{n,\epsilon,K}-\alpha\int_0^t\int_U\nabla\psi\cdot\nabla\rho^{n,\epsilon,K}\\
    +\sqrt{\epsilon}\int_0^t\int_U\sigma_n(\rho^{n,\epsilon,K})\nabla\psi\cdot\,d\xi^K+ \int_0^t\int_U\nabla\psi\cdot\nu(\rho^{n,\epsilon,K})\\
    + \frac{\epsilon}{2}\int_0^t\int_UF^K_1[\sigma_n'(\rho^{n,\epsilon,K})]^2\nabla\psi\cdot\nabla\rho^{n,\epsilon,K}+ \frac{\epsilon}{2}\int_0^t\int_U\sigma_n(\rho^{n,\epsilon,K})\sigma'_n(\rho^{n,\epsilon,K}) \nabla\psi\cdot F^K_2.
\end{multline*}
    \end{enumerate}
\end{definition} 

\begin{proposition}[Well-posedness of weak solutions to regularised equation]\label{proposition on existence of weak solution clt and ldp chapter}
    Under Assumption \ref{asm: assumptions for well posedenss of equation}, \ref{assumption on noise} and \ref{assumption of constant boundary data and constant initial condition}, for each $n\in\mathbb{N}, \alpha\in(0,1)$ there exists a unique weak solution $\rho^{n,\epsilon,K}$ of \eqref{eq: ito form of dk equation}. Furthermore, weak solutions and stochastic kinetic solutions coincide.
\end{proposition}
\begin{proof}
    The argument is the same as that of Proposition 4.30 \cite{popat2025well}.
\end{proof} 
For the below results, we need the auxiliary functions, which we can only define for general $p>2$ due to Assumption \ref{assumption of constant boundary data and constant initial condition} which tells us $\bar\rho$ is a fixed constant, independent of space or time.
\begin{definition}[Auxiliary functions for $\Phi$ and $\sigma$]\label{def: auxhilary functions}
Under Assumption \ref{assumption of constant boundary data and constant initial condition}, let $\bar\rho$ denote the constant solution to the hydrodynamic limit equation \eqref{eq: regularised hydrodynamic limit}.
      For every $p> 2$, define $\Theta_{\Phi,p}$ to be the unique function satisfying 
      \begin{equation}\label{eq: equation for theta phi p}
          \Theta_{\Phi,p}(\bar\rho)=0,\quad \text{and}\quad \Theta'_{\Phi,p}(\xi)=(\xi-\bar\rho)^{\frac{p-2}{2}}(\Phi'(\xi))^{1/2},
      \end{equation}
    and define $\Theta_{\sigma,p}$ to be the unique function satisfying 
    \begin{equation}\label{eq: equation for theta sigma p}
        \Theta_{\sigma,p}(\bar\rho)=0,\quad \text{and}\quad 
        \Theta_{\sigma,p}'(\xi)=(\xi-\bar\rho)^{p-2}\sigma(\xi)\sigma'(\xi).
    \end{equation}
\end{definition}
For $p=2$, the functions $\Theta_{\Phi,2}$ and $\Theta_{\sigma,2}$ are defined with the same antiderivatives as \eqref{eq: equation for theta phi p} and \eqref{eq: equation for theta sigma p} respectively, but with initial conditions replaced by $\Theta_{\Phi,2}(0)= \Theta_{\sigma,2}(0)=0$.
\begin{assumption}\label{asm: new assumptions}
    \begin{enumerate}
      \item Assume that $\Phi,\nu\in C^2(0,\infty)$ are such that there exists constants $c,\beta\in(0,\infty)$ with
      \begin{equation}\label{eq: bound on phi'' and nu''}
          |\Phi''(\xi)|+|\nu''(\xi)|\leq c(1+\xi^\beta).
      \end{equation}
       \item  For $\bar\rho$ the solution to the regularised hydrodynamic limit equation \eqref{eq: regularised hydrodynamic limit}, and for the functions $\Theta_{\Phi,p}$ and $\Theta_{\sigma,p}$ as in equations \eqref{eq: equation for theta phi p} and \eqref{eq: equation for theta sigma p} respectively, for every $p\geq2$, suppose that  there exists constants $q\in(0,2)$ and $c\in(0,\infty)$ such that for every $\xi\in(0,\infty)$,
\begin{align}\label{eq: bound for second energy estimate}
    (\xi-\bar\rho)^{p-2}\sigma^2(\xi)+\Theta_{\sigma,p}(\xi)\leq c(1+\Theta^q_{\Phi,p}(\xi)).
\end{align}
 \item For $\Theta_{\sigma,p}$ as in \eqref{eq: equation for theta sigma p}, suppose that for every $p\geq 2$ there exists constants $\gamma,c\in(0,\infty)$ such that for every $\xi\in(0,\infty)$,        
        \begin{equation}\label{eq: bound for sigma p and quotient to p/2}
            \sigma^p(\xi)+\left(\frac{\Theta_{\sigma,p}(\xi)}{(\xi-\bar\rho)^{p-2}}\right)^{p/2}\leq c(1+\xi^\gamma).
        \end{equation}
         \end{enumerate}
\end{assumption}
\begin{remark}\label{rmk: smuggling constant rho bar into the inequality for point 1 new asms}
    The first point is used in the proof of the central limit theorem for the regularised equation, when we make rigorous the formal computation from \eqref{formal computation for convergence of v epsilon to v}, in particular see the the computation between equations  \eqref{eq: first term in clt approximate equation estimate} and \eqref{eq: needing to bound l4 norm of v} below.
    We note that due to the presence of the constant on the right hand side of \eqref{eq: bound on phi'' and nu''} and Assumption \ref{assumption of constant boundary data and constant initial condition}, the bound can be alternatively stated
     \[|\Phi''(\xi)|+|\nu''(\xi)|\leq c(1+(\xi-\bar\rho)^\beta).\]
     This can be seen rigorously via the inequality $a^\beta\leq 2^\beta((a-b)^\beta+b^\beta)$.
\end{remark}
The second and third points in Assumption \ref{asm: new assumptions} are used to prove two different energy estimates in Propositions \ref{ppn: energy estimate for difference of regularised equation and hydrodynamic limit} and \ref{prop: second energy estimate for rho n epsilon k} respectively.
We give a comment below that justifies why they are reasonable.
\begin{remark}[Verifying point two of Assumption \ref{asm: new assumptions}]\label{rmk: on the second new assumption}
Since for small values of the argument $\xi\in(0,\infty)$, the functions on the left hand side of \eqref{eq: bound for second energy estimate} are bounded, the presence of the constant term on the right hand side gives that the inequality is always true.
Hence, we just need to verify the inequality \eqref{eq: bound for second energy estimate} for large $\xi\in(0,\infty)$, where by \say{large}, we mean large relative to $\bar\rho$ so that we have $(\xi-\bar\rho)\approx\xi$.\\
For large $\xi$, in the model case $\Phi(\xi)=\xi^m$ we have that up to a constant $\Theta_{\Phi,p}'(\xi)\approx \xi^{\frac{p+m-3}{2}}$, and so roughly
\begin{equation}\label{eq: bound for theta phi p for large xi}
    \Theta_{\Phi,p}(\xi)\approx \xi^{\frac{p+m-1}{2}}.
\end{equation}
Furthermore, when $\sigma(\xi)=\Phi^{1/2}(\xi)$ as in the particle system example \eqref{eq: equation linearised DK}, for large $\xi$ we have that $(\xi-\bar\rho)^{p-2}\sigma^2(\xi)\approx\xi^{p+m-2}$ and since $\sigma(\xi)\sigma'(\xi)\approx \xi^{m-1}$, we also have
       $\Theta_{\sigma,p}(\xi)\approx\xi^{m+p-2}$.
   Hence, the assumption is satisfied with the choice
 \begin{align}\label{eq: choice of q}
     q=\frac{2(m+p-2)}{m+p-1}<2.
 \end{align}
 We note that this also holds in the special case of the Dean--Kawasaki equation, $\Phi(\xi)=\xi, \sigma(\xi)=\sqrt{\xi}$ with the choice $m=1$ in \eqref{eq: choice of q}.
\end{remark}
\begin{remark}[Verifying point three of Assumption \ref{asm: new assumptions}]
  Recall from Remark \ref{rmk: on the second new assumption}, that in the case of interest from the particle system perspective $\Phi(\xi)=\xi^m$ and $\sigma(\xi)=\Phi^{1/2}(\xi)$, we had $\Theta_{\sigma,p}(\xi)\approx\xi^{m+p-2}$.
  Hence the quotient in \eqref{eq: bound for sigma p and quotient to p/2} is well defined as $\Theta_{\sigma,p}(\xi)$ decays to zero faster than $(\xi-\bar\rho)^{p-2}$ as $\xi\to\bar\rho$.
      Again by Remark \ref{rmk: on the second new assumption}, we just need to check the inequality \eqref{eq: bound for sigma p and quotient to p/2} for large values of the argument $\xi\in(0,\infty)$.
      
   In this case we have
     \[\sigma^p(\xi)+\left(\frac{\Theta_{\sigma,p}(\xi)}{(\xi-\bar\rho)^{p-2}}\right)^{p/2}\approx \xi^{mp/2},\]
     so $\gamma=\frac{mp}{2}$ suffices.
     By the same reasoning as Remark \ref{rmk: smuggling constant rho bar into the inequality for point 1 new asms}, we can replace the bound \eqref{eq: bound for sigma p and quotient to p/2} with the more convenient
     \begin{equation}\label{eq: restated bound for sigma p and quotient p/2}
         \sigma^p(\xi)+\left(\frac{\Theta_{\sigma,p}(\xi)}{(\xi-\bar\rho)^{p-2}}\right)^{p/2}\leq c(1+(\xi-\bar\rho)^\gamma).
     \end{equation} 
\end{remark}
To bound the right hand side of equation \eqref{eq: bound for second energy estimate} which will appear in the energy estimate, we need the below technical lemma.
\begin{lemma}[Interpolation estimate]\label{lemma: interpolation estimate}
Let $\Psi$ be any function, and $z:U\times[0,T]\to\mathbb{R}$ be measurable with constant boundary $z|_{\partial U}=M$  and $\Psi(z)\in L^2([0,T];H^1(U))$.
Then for every $q\in(0,2)$ and $\delta\in(0,1)$, there exists a constant $c\in(0,\infty)$ depending on the boundary data $M$ such that
\[\int_0^t\int_U\Psi^q(z)\leq c\left(\frac{t}{\delta} +\delta\int_0^t\int_U|\nabla\Psi(z)|^2\right).\]
\end{lemma}
\begin{proof}
 Since $q$ is strictly less than 2, multiplying the integrand by $1$ and using Cauchy-Schwarz and Young's inequality with exponent $\frac{2}{q}>1$ then gives for every $\delta\in(0,1)$
\begin{align*}
      \int_0^T\int_U|\Psi(z)|^{q}&\leq c\delta^{-1}T+ c\delta\int_0^T\int_U\Psi(z)^2.
  \end{align*}
Using the trivial inequality $a^2\leq 2(a-b)^2+2b^2$ with the constant $b=\Psi(M)$, and subsequently applying Poincar\'e inequality gives the claim,
  \begin{align*}
     \int_0^T\int_U|\Psi(z)|^{q}&\leq c\delta^{-1}T+c\delta\int_0^t\int_U\left(\Psi(z) -\Psi(M)\right)^2+c\delta\int_0^t\int_U\Psi(M)^2\\
      &\leq c\delta^{-1}T+ c\delta\int_0^t\int_U\left|\nabla\left(\Psi(z) -\Psi(M)\right)\right|^2\nonumber\\
      &\leq c\delta^{-1}T+ c\delta\int_0^t\int_U\left|\nabla\Psi(z)\right|^2.\nonumber \tag*{\qedhere}
  \end{align*}
\end{proof}
The bound in the above lemma can easily be adapted to functions $z$ with non-constant boundary data, but this will not be needed here.
Our first result will be proving two $L^p(U\times[0,T])$-energy estimates for the regularised equation, both of which will be useful in the sequel.

\begin{proposition}[$p$-independent energy estimate]\label{ppn: energy estimate for difference of regularised equation and hydrodynamic limit}
    For fixed $\alpha\in(0,1)$, let $\rho^{n,\epsilon,K}$ be the weak solution of the regularied equation \eqref{equation for rho n epsilon k} that exists due to Proposition \ref{proposition on existence of weak solution clt and ldp chapter}, and let $\bar\rho$ be the solution to the regularised hydrodynamic limit equation \eqref{eq: regularised hydrodynamic limit}.
    Suppose Assumptions \ref{assumption on noise}, \ref{assumption of constant boundary data and constant initial condition}, \ref{asm: new assumptions} and \ref{asm: assumptions for well posedenss of equation} hold, and let $\Theta_{\Phi,p}$ be the unique function defined in Definition \ref{def: auxhilary functions}.
Then we have for every $\epsilon\in(0,1)$ and $K\in\mathbb{N}$, that there exists a constant $c\in(0,\infty)$ independent of $T, \alpha, \epsilon$, $n$ and $K$ such that for every $p\geq 2$, 
\small
\begin{multline}\label{eq: Lp space-time estimate for rho n epsilon k}
    \frac{1}{Tp(p-1)}\mathbb{E}\left[\int_0^T\int_U \left(\rho^{n,\epsilon,K}-\bar\rho\right)^p\right]+\mathbb{E}\left[\int_0^T\int_U|\nabla\Theta_{\Phi,p}(\rho^{n,\epsilon,K})|^2\right]+\frac{4\alpha}{p^2}\mathbb{E}\left[\int_0^T\int_U|\nabla(\rho^{n,\epsilon,K}-\bar\rho)^{p/2}|^2\right]\\
    \leq c \epsilon T\left(\|\nabla\cdot F^K_2\|_{L^\infty(U)}+\|F^K_3\|
    _{L^\infty(U)}\right),
\end{multline}
\normalsize
and 
\small
\begin{multline}\label{eq: l infinity time lp space estimate for rho n epsilon k}
    \frac{1}{p(p-1)}\sup_{t\in[0,T]}\mathbb{E}\left[\int_U \left(\rho^{n,\epsilon,K}-\bar\rho\right)^p\right]+\mathbb{E}\left[\int_0^T\int_U|\nabla\Theta_{\Phi,p}(\rho^{n,\epsilon,K})|^2\right]+\frac{4\alpha}{p^2}\mathbb{E}\left[\int_0^T\int_U|\nabla(\rho^{n,\epsilon,K}-\bar\rho)^{p/2}|^2\right]\\
    \leq c \epsilon T\left(\|\nabla\cdot F^K_2\|_{L^\infty(U)}+\|F^K_3\|_{L^\infty(U)}\right).
\end{multline}
\normalsize
\end{proposition}

\begin{proof}
 By It\^o's formula and Assumption \ref{assumption of constant boundary data and constant initial condition}, for every $p\geq 2$ we have
 \small
\begin{align*}
    &\frac{1}{p(p-1)}\int_U \left(\rho^{n,\epsilon,K}-\bar\rho\right)^p=-\int_0^t\int_U(\rho^{n,\epsilon,K}-\bar\rho)^{p-2}\nabla (\rho^{n,\epsilon,K}-\bar\rho)\cdot\nabla\Phi(\rho^{n,\epsilon,K})\\
    &-\alpha\int_0^t\int_U(\rho^{n,\epsilon,K}-\bar\rho)^{p-2}\nabla (\rho^{n,\epsilon,K}-\bar\rho)\cdot\nabla\rho^{n,\epsilon,K}\\
    &+\int_0^t\int_U(\rho^{n,\epsilon,K}-\bar\rho)^{p-2}\nabla (\rho^{n,\epsilon,K}-\bar\rho)\cdot\nu(\rho^{n,\epsilon,K})-\frac{\sqrt{\epsilon}}{(p-1)}\int_0^t\int_U (\rho^{n,\epsilon,K}-\bar\rho)^{p-1}\nabla\cdot(\sigma_n(\rho^{n,\epsilon,K})\dot{\xi}^K)\\
    &-\frac{\epsilon}{2}\int_0^t\int_U(\rho^{n,\epsilon,K}-\bar\rho)^{p-2}\nabla (\rho^{n,\epsilon,K}-\bar\rho)\cdot\left(F_1^K(\sigma_n'(\rho^{n,\epsilon,K}))^2\nabla\rho^{n,\epsilon,K}+\sigma_n(\rho^{n,\epsilon,K})\sigma'_n(\rho^{n,\epsilon,K}) F^K_2\right)\\
    &+\frac{\epsilon}{2}\int_0^t\int_U(\rho^{n,\epsilon,K}-\bar\rho)^{p-2}\left(F^K_1(\sigma_n'(\rho^{n,\epsilon,K}))^2|\nabla\rho^{n,\epsilon,K}|^2+2\sigma_n(\rho^{n,\epsilon,K})\sigma_n'(\rho^{n,\epsilon,K})F^K_2\cdot\nabla\rho^{n,\epsilon,K}+F^K_3\sigma_n^2(\rho^{n,\epsilon,K})\right).
\end{align*}
\normalsize
    Using that $\nabla (\rho^{n,\epsilon,K}-\bar\rho)=\nabla\rho^{n,\epsilon,K}$ and the cancellation of the terms in the final two lines gives after re-arranging, that for the function $\Theta_{\Phi,p}$ defined in Definition \ref{def: auxhilary functions}
\begin{multline}\label{eq: lp estimate for rho n epsilon k}
    \frac{1}{p(p-1)}\int_U \left(\rho^{n,\epsilon,K}-\bar\rho\right)^p+\int_0^t\int_U|\nabla\Theta_{\Phi,p}(\rho^{n,\epsilon,K})|^2+\frac{4\alpha}{p^2}\int_0^t\int_U|\nabla(\rho^{n,\epsilon,K}-\bar\rho)^{p/2}|^2=\\
    \int_0^t\int_U(\rho^{n,\epsilon,K}-\bar\rho)^{p-2}\nabla \rho^{n,\epsilon,K}\cdot\nu(\rho^{n,\epsilon,K})-\frac{\sqrt{\epsilon}}{(p-1)}\int_0^t\int_U (\rho^{n,\epsilon,K}-\bar\rho)^{p-1}\nabla\cdot(\sigma_n(\rho^{n,\epsilon,K})\dot{\xi}^K)\\
    +\frac{\epsilon}{2}\int_0^t\int_U(\rho^{n,\epsilon,K}-\bar\rho)^{p-2}\left(\sigma_n(\rho^{n,\epsilon,K})\sigma_n'(\rho^{n,\epsilon,K})F^K_2\cdot\nabla\rho^{n,\epsilon,K}+F^K_3\sigma_n^2(\rho^{n,\epsilon,K})\right).
\end{multline}
The first term on the right hand side of \eqref{eq: lp estimate for rho n epsilon k} vanishes due to the fact that the boundary condition is constant,  
\begin{align}\label{eq: nu vanishing on boundary}
    \int_0^t\int_U(\rho^{n,\epsilon,K}-\bar\rho)^{p-2}\nabla \rho^{n,\epsilon,K}\cdot\nu(\rho^{n,\epsilon ,K})&=\sum_{i=1}^d\int_0^t \int_U\partial_i\left(\Theta_{\nu,p,i}(\rho^{n,\epsilon,K})\right)\nonumber\\
    &=t\sum_{i=1}^d\Theta_{\nu,p,i}(\Phi^{-1}(\bar{f}))\int_
    {\partial U}\hat{\eta_i}=0,
\end{align}
where we defined for $i=1,\hdots,d$ the unique functions $\Theta_{\nu,p,i}$ by $\Theta_{\nu,p,i}(0)=0$ and $\Theta'_{\nu,p,i}(\xi)=(\xi-\bar\rho)^{p-2}\nu_i(\xi)$, and  $\hat{\eta}=(\hat{\eta}_1,\hdots,\hat{\eta}_d)$ denotes the outward pointing unit normal at the boundary.
The final equality follows from the divergence theorem, letting $e_i$ denote the standard basis vector in $\mathbb{R}^d$ in the $i$'th direction, we have
\[\int_{\partial U}\hat\eta_i=\int_{\partial U}e_i\cdot\hat\eta=\int_U\nabla\cdot e_i=\int_U\frac{\partial}{\partial_{x_i}}(1)=0.\]
The noise term in \eqref{eq: lp estimate for rho n epsilon k} is a martingale, so vanishes under expectation.
We are left with handling the two terms in the final line of \eqref{eq: lp estimate for rho n epsilon k}.
For the first, we have using integration by parts
    \begin{align}\label{eq: integrating by parts first term in lp estimate}
    \frac{\epsilon}{2}\int_0^t\int_U(\rho^{n,\epsilon,K}-\bar\rho)^{p-2}\sigma_n(\rho^{n,\epsilon,K})\sigma_n'(\rho^{n,\epsilon,K})F^K_2\cdot\nabla\rho^{n,\epsilon,K}&=-\frac{\epsilon}{2}\int_0^t\int_U\Theta_{\sigma_n,p}(\rho^{n,\epsilon,K})\nabla\cdot F^K_2.
\end{align}
for the function $\Theta_{\sigma_n,p}$ defined in Definition \ref{def: auxhilary functions}.
In the case $p=2$ we do not pick up a boundary term when integrating by parts due to the Dirichlet boundary data of the noise coefficients, see point 3 of Definition \ref{assumption on noise}, that gives
\begin{align}\label{eq: no boundary terms when IBP F_2}
    \frac{\epsilon}{2}\int_0^t\int_{\partial U}\sigma_n^2(M)F_2^K\cdot\hat{\eta}=\frac{\epsilon \sigma_n^2(M)}{2}\sum_{k=1}^K\int_0^t\int_{\partial U}f_k \nabla f_k\cdot\hat\eta=0 .
\end{align}
Using an $L^\infty(U)$-estimate for $\nabla\cdot F^K_2$, and using point 2 of Assumption \ref{asm: new assumptions}, we have that there exists constants $q\in(0,2)$ and $c\in(0,\infty)$ such that
        \begin{multline*}
    \frac{\epsilon}{2}\int_0^t\int_U(\rho^{n,\epsilon,K}-\bar\rho)^{p-2}\sigma_n(\rho^{n,\epsilon,K})\sigma_n'(\rho^{n,\epsilon,K})F^K_2\cdot\nabla\rho^{n,\epsilon,K}\\
    \leq c\epsilon \|\nabla\cdot F^K_2\|_{L^\infty(U)}\int_0^t\int_U(1+\Theta^q_{\Phi,p}(\rho^{n,\epsilon,K})).
\end{multline*}
For the final term of \eqref{eq: lp estimate for rho n epsilon k}, using an $L^\infty(U)$-estimate for $F_3^K$ and point 2 of Assumption \ref{asm: new assumptions} gives the existence of constants $q\in(0,2)$ and $c\in(0,\infty)$ such that 
\begin{align*}
     \frac{\epsilon}{2}\int_0^t\int_U(\rho^{n,\epsilon,K}-\bar\rho)^{p-2}F^K_3\sigma_n^2(\rho^{n,\epsilon,K})\leq c\epsilon \|F^K_3\|_{L^\infty(U)}\int_0^t\int_U(1+\Theta^q_{\Phi,p}(\rho^{n,\epsilon,K})).
\end{align*}
The integrals of $\Theta^q_{\Phi,p}$ are controlled using Lemma \ref{lemma: interpolation estimate}, again picking $\delta>0$ small enough to absorb the resulting gradient term to the left hand side of the estimate.
Putting everything together, we are left with
\begin{multline}\label{eq: getting to the energy estimates}
    \frac{1}{p(p-1)}\mathbb{E}\left[\int_U \left(\rho^{n,\epsilon,K}-\bar\rho\right)^p\right]+\mathbb{E}\left[\int_0^t\int_U|\nabla\Theta_{\Phi,p}(\rho^{n,\epsilon,K})|^2\right]+\frac{4\alpha}{p^2}\mathbb{E}\left[\int_0^t\int_U|\nabla(\rho^{n,\epsilon,K}-\bar\rho)^{p/2}|^2\right]\\
    \leq c \epsilon t\left(\|\nabla\cdot F^K_2\|_{L^\infty(U)}+\|F^K_3\|_{L^\infty(U)}\right).
\end{multline}

The first estimate \eqref{eq: Lp space-time estimate for rho n epsilon k} then follows from integrating the inequality \eqref{eq: getting to the energy estimates} over $[0,T]$ with respect to $t$, changing the order of integration, and noticing that when $f$ is non-negative,
\[\int_0^T\int_0^tf(s)\,ds\,dt=\int_0^T(T-s)f(s)\,ds\leq T\int_0^Tf(s)\,ds.\]
The second estimate \eqref{eq: l infinity time lp space estimate for rho n epsilon k} follows by taking supremum in time on both sides of \eqref{eq: getting to the energy estimates}.
\end{proof}
\begin{remark}\label{rmk: choice of scaling in energy estimate}
The above estimate is stable for every fixed $K\in\mathbb{N}$.
However in order to take the $K\to\infty$ limit and ensure the right hand side of the estimates converge to zero, one needs a joint scaling regime $\epsilon\to0, K\to\infty$ such that
    \[\epsilon\left(\|\nabla\cdot F^K_2\|_{L^\infty(U)}+\|F^K_3\|_{L^\infty(U)}\right)\to0.\]
Based on other energy estimates below, we will deduce a final joint scaling regime later in Remark \ref{rmk: choice of joint scaling}.

\end{remark}

Using the third point of Assumption \ref{asm: new assumptions}, we can also obtain a bound depending on $p$ on the right hand side. 

\begin{proposition}[$p$-dependent estimate]\label{prop: second energy estimate for rho n epsilon k}
  For fixed $\alpha\in(0,1)$, let $\rho^{n,\epsilon,K}$ be the weak solution of the regularised equation \eqref{equation for rho n epsilon k} that exists due to Proposition \ref{proposition on existence of weak solution clt and ldp chapter}, and let $\bar\rho$ be the solution to the regularised hydrodynamic limit equation \eqref{eq: regularised hydrodynamic limit}.
    Suppose Assumptions \ref{assumption on noise}, \ref{assumption of constant boundary data and constant initial condition}, \ref{asm: new assumptions} and \ref{asm: assumptions for well posedenss of equation} hold, and let $\Theta_{\Phi,p}$ be the unique function defined in Definition \ref{def: auxhilary functions}.
Then we have for every $\epsilon\in(0,1)$ and $K\in\mathbb{N}$, that there exists a constant $c\in(0,\infty)$ independent of $T, \alpha, \epsilon$, $n$ and $K$ such that for every $p\geq 2$,
\small
\begin{multline*}
    \frac{1}{Tp(p-1)}\mathbb{E}\left[\int_0^T\int_U \left(\rho^{n,\epsilon,K}-\bar\rho\right)^p\right]+\mathbb{E}\left[\int_0^T\int_U|\nabla\Theta_{\Phi,p}(\rho^{n,\epsilon,K})|^2\right]+\frac{4\alpha}{p^2}\mathbb{E}\left[\int_0^T\int_U|\nabla(\rho^{n,\epsilon,K}-\bar\rho)^{p/2}|^2\right]\\
    \leq c \epsilon^{p/2} T\left(\|\nabla\cdot F^K_2\|^{p/2}_{L^\infty(U)}+\|F^K_3\|^{p/2}_{L^\infty(U)}\right)\left(1+\epsilon T\left(\|\nabla\cdot F^K_2\|_{L^\infty(U)}+\|F^K_3\|
    _{L^\infty(U)}\right)\right),
\end{multline*}
\normalsize
and 
\small
\begin{multline*}
    \frac{1}{p(p-1)}\sup_{t\in[0,T]}\mathbb{E}\left[\int_U \left(\rho^{n,\epsilon,K}-\bar\rho\right)^p\right]+\mathbb{E}\left[\int_0^T\int_U|\nabla\Theta_{\Phi,p}(\rho^{n,\epsilon,K})|^2\right]+\frac{4\alpha}{p^2}\mathbb{E}\left[\int_0^T\int_U|\nabla(\rho^{n,\epsilon,K}-\bar\rho)^{p/2}|^2\right]\\
    \leq c \epsilon^{p/2} T\left(\|\nabla\cdot F^K_2\|^{p/2}_{L^\infty(U)}+\|F^K_3\|^{p/2}_{L^\infty(U)}\right)\left(1+\epsilon T\left(\|\nabla\cdot F^K_2\|_{L^\infty(U)}+\|F^K_3\|
    _{L^\infty(U)}\right)\right).
\end{multline*}
\normalsize
\end{proposition}
\begin{proof}
    The proof follows via the same estimates as the proof of Proposition \ref{ppn: energy estimate for difference of regularised equation and hydrodynamic limit}.
    It is just a matter of how the terms in the final line of equation \eqref{eq: lp estimate for rho n epsilon k} are handled.
For the first term, after integrating by parts as in \eqref{eq: integrating by parts first term in lp estimate} and using an $L^\infty(U)$-estimate for $\nabla\cdot F_2^K$, H\"older's inequality and Young's inequality with exponent $\frac{p}{p-2}>1$, and point 3 of Assumption \ref{asm: new assumptions} (specifically \eqref{eq: restated bound for sigma p and quotient p/2}) gives that there exists constants $\beta, c\in(0,\infty)$ such that for every $\delta\in(0,1)$, 
\begin{multline}\label{bound for theta sigma n p}
    \frac{\epsilon}{2}\int_0^t\int_U\Theta_{\sigma_n,p}(\rho^{n,\epsilon,K})\nabla\cdot F^K_2\leq c\epsilon \|\nabla\cdot F_2^K\|_{L^\infty(U)}\int_0^t\int_U(\rho^{n,\epsilon,K}-\bar\rho)^{p-2}\frac{\Theta_{\sigma_n,p}(\rho^{n,\epsilon,K})}{(\rho^{n,\epsilon,K}-\bar\rho)^{p-2}}\\
    \leq \delta\int_0^t\int_U (\rho^{n,\epsilon,K}-\bar\rho)^{p}+c\delta^{-1}\epsilon^{p/2}\|\nabla\cdot F_2^K\|^{p/2}_{L^\infty(U)}\int_0^t\int_U \left(\frac{\Theta_{\sigma_n,p}(\rho^{n,\epsilon,K})}{(\rho^{n,\epsilon,K}-\bar\rho)^{p-2}}\right)^{p/2}\\
    \leq \delta\int_0^t\int_U (\rho^{n,\epsilon,K}-\bar\rho)^{p}+c\delta^{-1}\epsilon^{p/2}\|\nabla\cdot F_2^K\|^{p/2}_{L^\infty(U)}\int_0^t\int_U (1+(\rho^{n,\epsilon,K}-\bar\rho)^\gamma).
\end{multline}
In precisely the same way, by taking the $L^\infty(U)$-norm of $F_3^K$, H\"older's inequality and Young's inequality with exponent $\frac{p}{p-2}>1$, and point 3 of Assumption \ref{asm: new assumptions} (specifically \eqref{eq: restated bound for sigma p and quotient p/2}), we get for the final term of \eqref{eq: lp estimate for rho n epsilon k} that there exists constants $\beta, c\in(0,\infty)$ such that for every $\delta\in(0,1)$,
\begin{multline*}
    \frac{\epsilon}{2}\int_0^t\int_U(\rho^{n,\epsilon,K}-\bar\rho)^{p-2}F^K_3\sigma_n^2(\rho^{n,\epsilon,K})\leq \frac{\epsilon}{2}\|F^K_3\|_{L^\infty(U)}\int_0^t\int_U(\rho^{n,\epsilon,K}-\bar\rho)^{p-2}\sigma_n^2(\rho^{n,\epsilon,K})\\
    \leq \delta\int_0^t\int_U (\rho^{n,\epsilon,K}-\bar\rho)^{p}+c\delta^{-1}\epsilon^{p/2}\|F_3^K\|^{p/2}_{L^\infty(U)}\int_0^t\int_U \sigma_n^p(\rho^{n,\epsilon,K})\\
        \leq \delta\int_0^t\int_U (\rho^{n,\epsilon,K}-\bar\rho)^{p}+c\delta^{-1}\epsilon^{p/2}\|F_3^K\|^{p/2}_{L^\infty(U)}\int_0^t\int_U (1+(\rho^{n,\epsilon,K}-\bar\rho)^\gamma).
\end{multline*}
For both the above term and the final line of equation \eqref{bound for theta sigma n p}, we pick $\delta>0$ small enough so that the first term can be absorbed onto the left hand side.
For final term above, without loss of generality we can pick $\gamma\geq 2$, for which we can use the first energy estimate Proposition \ref{ppn: energy estimate for difference of regularised equation and hydrodynamic limit}
to obtain the bound
\begin{multline*}
    \frac{1}{p(p-1)}\mathbb{E}\left[\int_U \left(\rho^{n,\epsilon,K}-\bar\rho\right)^p\right]+\mathbb{E}\left[\int_0^t\int_U|\nabla\Theta_{\Phi,p}(\rho^{n,\epsilon,K})|^2\right]+\frac{4\alpha}{p^2}\mathbb{E}\left[\int_0^t\int_U|\nabla(\rho^{n,\epsilon,K}-\bar\rho)^{p/2}|^2\right]\\
    \leq c \epsilon^{p/2} t\left(\|\nabla\cdot F^K_2\|^{p/2}_{L^\infty(U)}+\|F^K_3\|^{p/2}_{L^\infty(U)}\right)\left(1+\epsilon t\left(\|\nabla\cdot F^K_2\|_{L^\infty(U)}+\|F^K_3\|
    _{L^\infty(U)}\right)\right).
\end{multline*}
The two estimates are then obtained in the same way as Proposition \ref{ppn: energy estimate for difference of regularised equation and hydrodynamic limit}.
\end{proof}
\begin{remark}
Based on the joint scaling from Remark \ref{rmk: choice of scaling in energy estimate}, whenever we have that\\  $\epsilon\left(\|\nabla\cdot F^K_2\|_{L^\infty(U)}+\|F^K_3\|_{L^\infty(U)}\right)\leq 1$, the rate of convergence to zero of the right hand of Proposition \ref{prop: second energy estimate for rho n epsilon k} is governed by
\[ \epsilon^{p/2} \left(\|\nabla\cdot F^K_2\|^{p/2}_{L^\infty(U)}+\|F^K_3\|^{p/2}_{L^\infty(U)}\right).\]
Whenever $p>2$, this is a faster rate of convergence than the $p$-independent estimate of Proposition \ref{ppn: energy estimate for difference of regularised equation and hydrodynamic limit}.
\end{remark}


\subsection{Central limit theorem for approximate SPDE}

In this section we will prove a central limit theorem for the regularised equation \eqref{equation for rho n epsilon k}. 
We prove fluctuations converge strongly in the space $L^2([0,T];H^{-s}(U))$, for every $s>d/2$, to the linearised SPDE \eqref{eq: linearised spde v}. 
First we define a strong solution to the linearised SPDE, which we again emphasise is distribution valued.
\begin{definition}[Strong solution to linearised SPDE]\label{definition strong solution to ou process}
Let $\dot\xi$ denote an $\mathbb{R}^d$-valued space-time white noise and let $\bar\rho$ denote the weak solution to the hydrodynamic limit equation \eqref{heat equation hydrodynamic limit}.
A strong solution to the linearised SPDE \eqref{eq: linearised spde v} is an $\mathcal{F}_t$-adapted and almost surely continuous $H^{-s}(U)$-valued process $v$ satisfying that for every $\psi\in H^s_0(U)$,
\[\langle v(t),\psi\rangle_s=\int_0^t\langle v(r),\Phi'(\bar\rho)\Delta\psi(r)+\nu'(\bar\rho)\cdot\nabla\psi\rangle_s\,dr+\int_0^t\int_U\sigma(\bar\rho)\nabla\psi\cdot\dot{\xi}\,dr.\]
We denote by $\langle\cdot,\cdot\rangle_s$ the dual pairing between $H^{-s}(U)$ and $H^s(U)$.
\end{definition} 
\begin{remark}\label{remark on boundary data of OU process} The boundary condition is implicit in the above expression. In particular, for the first term when we applied integration by parts and did not pick up any boundary terms.
\end{remark}
\begin{remark}\label{rmk: no need to smooth v}
In the central limit theorem results we always assume that Assumption \ref{assumption of constant boundary data and constant initial condition} holds.
This implies that $\sigma(\bar\rho)$ is either identically zero (in which case there is no noise in the system), or is a constant, uniformly bounded away from zero.
In both cases we do not need to smooth $\sigma$.
Furthermore, owing to point 2 of Assumption \ref{asm: assumptions for well posedenss of equation}, $\Phi'(\bar\rho)$ is a strictly positive constant and so the corresponding term acts as a (smoothing) Laplacian term.
Hence, regularisation in way of equations \eqref{equation for rho n epsilon k} and \eqref{eq: regularised hydrodynamic limit} is not required for energy estimates of the equation.
\end{remark}
\begin{proposition}[Existence and uniqueness of strong solutions to the linearised SPDE]\label{proposition existence of strong solutions to OU process}
    Under Assumptions \ref{assumption on noise}, \ref{assumption of constant boundary data and constant initial condition} and \ref{asm: assumptions for well posedenss of equation}, there is a unique strong solution of the linearised stochastic PDE \eqref{eq: linearised spde v} as defined in Definition \ref{definition strong solution to ou process}.
\end{proposition}  
The existence in the below proof follows Proposition 3.7 of \cite{dirr2020conservative}.

\begin{proof}
     Let $\xi^K$ be the finite dimensional noise from Definition \ref{definition truncated noise}, satisfying Assumption \ref{assumption on noise} and $\bar\rho$ the solution of hydrodynamic limit equation \eqref{heat equation hydrodynamic limit}.
Let $v_{K}$ denote the solution of the linearised SPDE with the finite dimensional noise, that is, the solution to 
\[\partial_t v_{K}=\Delta(\Phi'(\bar\rho)v_{K}))-\nabla\cdot(\sigma(\bar\rho)\dot\xi^K+\nu'(\bar\rho)v_{K}),\]
with zero initial data and boundary condition.
Simplified versions of Theorem 3.6 and Theorem 4.38 of \cite{popat2025well} proves uniqueness and existence of stochastic kinetic solutions and weak solutions for every $k\in\mathbb{N}$.\\
For $s>\frac{d+2}{2}$ consider $z_{K}:=(-\Delta)^{-s/2}v_{K}$, which is function valued.\\
Applying It\^o's formula to $z^2_K$ and using Assumption \ref{assumption of constant boundary data and constant initial condition} gives for the first order term
\begin{multline*}
    \int_0^t\int_U z_{K} dz_{K}=\int_0^t\int_U z_{K}\left((-\Delta)^{-s/2}\Delta(\Phi'(\bar\rho)v_{K}))-(-\Delta)^{-s/2}\nabla\cdot(\sigma(\bar\rho)\dot\xi^K+\nu'(\bar\rho)v_{K})\right)\\
    =-\int_0^t\int_U\nabla z_{K}\cdot \left((-\Delta)^{-s/2}\nabla(\Phi'(\bar\rho)v_{K}))-(-\Delta)^{-s/2}(\sigma(\bar\rho)\dot\xi^K+\nu'(\bar\rho)v_{K})\right)\\
    =-\int_0^t\int_U \Phi'(\bar\rho) |\nabla z_{K}|^2 + \int_0^t\int_U\sigma(\bar\rho)(-\Delta)^{-s/2}\nabla z_{K}\cdot d\xi^K+\int_0^t\int_U\nabla z_{K}\cdot\nu'(\bar\rho)z_{K}.
\end{multline*}

The final term vanishes by analogous reasoning to the term involving $\nu$ in Proposition \ref{ppn: energy estimate for difference of regularised equation and hydrodynamic limit}, see equation \eqref{eq: nu vanishing on boundary} there.

The It\^o correction can be bounded by Assumption \ref{assumption on noise} since $s>\frac{d+2}{2}$ to give for a constant independent of $K$,
\begin{align*}
    \frac{1}{2}\int_0^t\int_U \,d\langle z_{K}\rangle_t&=\frac{1}{2}\sum_{k=1}^K\int_0^t\int_U \sigma^2(\bar\rho)|(-\Delta)^{-s/2}\nabla f_k|^2=\frac{t|\sigma^2(\bar\rho)|}{2}\sum_{k=1}^K\|f_k\|^2_{H^{-s+1}(U)}\leq ct.
\end{align*}
The final inequality follows from the fact that we can upper bound the partial sum by the infinite sum, which is bounded and independent of $K$.
Putting everything together and taking supremum over time and an expectation, we have
\small
\begin{align*}
    &\mathbb{E}\left[\sup_{t\in[0,T]} \|z_{K}\|^2_{L^2(U)}+\Phi'(\bar\rho)\int_0^T\int_U|\nabla z_{K}|^2\right]\leq \mathbb{E}\left[\sup_{t\in[0,T]}\left|\int_0^t\int_U\sigma(\bar\rho)\left((-\Delta)^{-s/2}\nabla z_{K}\right)\cdot d\xi^K\right|\right]+cT.
\end{align*}
\normalsize
For the noise term in the above equation, Assumption \ref{assumption of constant boundary data and constant initial condition} alongside Burkholder-Davis-Gundy inequality, H\"older and Young inequalities, the fact that $s>\frac{3}{2}$, meaning it compensates the gradient, shows that for a constant $c\in(0,\infty)$ depending on $\bar\rho$, but independent of $T$ and $K$,
\begin{align*}
    \mathbb{E}\left[\max_{t\in[0,T]}\left|\int_0^t\int_U\sigma(\bar\rho)\left((-\Delta)^{-s/2}\nabla z_{K}\right)\cdot d\xi^K\right|\right]&\leq c|\sigma(\bar\rho)|\mathbb{E}\|(-\Delta)^{-s/2}\nabla z_{K}\|_{L^2(U\times[0,T])}\\
    &\leq c\mathbb{E}\|z_K\|_{L^2(U\times[0,T])}.
\end{align*}
Putting everything together, we get that there is a constant $c\in(0,\infty)$ independent of $T$ and $K$ such that
\[\mathbb{E}\left[\sup_{t\in[0,T]} \|z_{K}\|^2_{L^2(U)}+\Phi'(\bar\rho)\int_0^T\int_U|\nabla z_{K}|^2\right]\leq c(T+ \mathbb{E}\|z_K\|_{L^2(U\times[0,T])}).\]
Transferring the estimate of $z_K$ to $v_K$ and using Gr\"onwall's inequality, it follows that there is a constant $c\in(0,\infty)$ independent of $K$ such that
\begin{align}\label{l2 estimate for smoothed ou process}
    \mathbb{E}\left[ \|v_K\|^2_{L^\infty([0,T];H^{-s}(U))}+\|v_K\|^2_{L^2([0,T];H^{-s+1}(U))}\right]\leq cT.
\end{align}
Furthermore, to estimate higher order time regularity, observe that distributionally we have due to Assumption \ref{assumption of constant boundary data and constant initial condition}
\[ v_{K}(\cdot,t)=\int_0^t\Phi'(\bar\rho)\Delta v_{K}-\nabla\cdot(\nu'(\bar\rho)v_{K})\,dt-\int_0^t\nabla\cdot(\sigma(\bar\rho)d\xi^K)=I^{f.v.}_t(\cdot)+I^{mart}_t(\cdot).\]
By equation \eqref{l2 estimate for smoothed ou process} we get that the finite variation term satisfies for constant $c\in(0,\infty)$ independent of $K$,
\[\|I^{f.v.}_\cdot\|_{W^{1,2}([0,T];H^{-s+1}(U))}\leq c\|v_K\|_{L^2([0,T];H^{-s}(U))}.\]
From Assumption \ref{assumption of constant boundary data and constant initial condition}, the Burkholder-Davis-Gundy inequality and Lemma \ref{lemma sum of norm of e_k}, we have for a constant $c$ depending on $T,\beta,\bar\rho$, but independent of $K$, such that
\begin{align}\label{eq: higher order time rgularity}
    \mathbb{E}\|I_\cdot^{mart}\|^2_{W^{\beta,2}([0,T];H^{-(s+1)}(U))}&=\mathbb{E}\int_0^T\int_0^T|s-t|^{-(1+2\beta)}\left\|\sum_{k=1}^K\int_s^t\sigma(\bar\rho) f_k\,dB_t^k\right\|^2_{H^{-s}(U)}\nonumber\\
     &\leq c\sigma^2(\bar\rho)\mathbb{E}\int_0^T\int_0^T|s-t|^{-(1+2\beta)}\sum_{k=1}^K\int_s^t\|f_k\|^2_{H^{-s}(U)}\nonumber\\
     &\leq c\int_0^T\int_0^T|s-t|^{-2\beta}\leq c,
\end{align}
where again in the final inequality we upper bounded the partial sum up to $K$ by the infinite sum and used point 2 of Assumption \ref{assumption on noise}.
By the above estimates and the compact embedding of $H^{-s}$ into $H^{-s'}$ whenever $s<s'\in(0,\infty)$, the Aubin-Lions-Simon lemma (\cite{aubin1963theoreme,lions1969quelques,simon1986compact}) tells us that due to estimates \eqref{l2 estimate for smoothed ou process} and \eqref{eq: higher order time rgularity}, the laws of $\{v_K\}_{K\in\mathbb{N}}$ are tight on $L^2([0,T];H^{-s}(U))$ for every $s>d/2$.
Following arguments similar to Theorem 3.6 of \cite{popat2025well}, we get that in the $K\to\infty$ limit, $\{v_K\}_{k\in\mathbb{N}}$ converge in law in $L^2([0,T];H^{-s}(U))$ to an element $v\in L^2([0,T]\times\Omega;H^{-s}(U))$ for every $s>d/2$, satisfying for every $\psi\in C_c^\infty(U)$,
\[\langle v(t),\psi\rangle_s=\int_0^t\langle v(r),\Phi'(\bar\rho)\Delta\psi(r)+\nu'(\bar\rho)\cdot\nabla\psi\rangle_s\,dr+\int_0^t\int_U\sigma(\bar\rho)\nabla\psi\cdot\,d\xi\,dr.\]
Both terms on the right hand side are continuous in time, so it follows that for every $k\in\mathbb{N}$, and $\{f_k\}_{k\in\mathbb{N}}$ as in Definition \ref{definition truncated noise}, $t\mapsto \langle v(t),f_k\rangle_s$ has a  continuous modification in $L^2([0,T])$.
By the $\mathbb{P}$-a.s. boundedness of $v$ in $L^2([0,T];H^{-s}(U))$, we know that there exists a $H^{-s}(U)$-continuous modification, denoted again by $v$, satisfying the above equation for every $\psi\in C_c^\infty(U)$ and $t\in[0,T]$.
This completes the proof of existence.\\
For uniqueness, let $v,\tilde{v}$ be two strong solutions to the linearised SPDE \eqref{eq: linearised spde v} as in Definition \ref{definition strong solution to ou process}. Then $w=v-\tilde{v}$ is a distributional solution to 
\begin{equation*}
   \begin{cases}
    \partial_tw=\Phi'(\bar\rho)\Delta w-\nabla\cdot(\nu'(\bar\rho)w) & \text{on} \hspace{5pt} U\times(0,T],\\
w=0,&\text{on}\hspace{5pt} \partial U\times[0,T],\\
    w(\cdot,t=0)=0 ,&\text{on}\hspace{5pt} U\times\{t=0\}.
   \end{cases} 
\end{equation*}
Let $\{e_k\}_{k\in\mathbb{N}}$ denote the eigenfunctions of the Laplacian with zero Dirichlet boundary conditions as in Example \ref{example of noise: eigenvalues of laplacian}, which we recall forms an orthonormal basis in $L^2(U)$.
Since $v$ and $\tilde{v}$ are elements of $H^{-s}(U)$, so is $w=v-\tilde{v}$. 
Hence $w$ can be represented weakly in terms of the eigenfunctions $\{e_k\}$ by
\begin{equation}\label{expansion of w in terms of e_k}
    w(x,t)=\sum_{k=1}^\infty w_k(t)e_k(x),\,\text{where}\,\, w_k(t):=\langle w(t),e_k\rangle_s,
\end{equation}
and the convergence of the above series is in $H^{-s}(U)$.
The expansion is justified because $\{e_k\}_{k\in\mathbb{N}}$ still form a basis of the space $H^{-s}(U)$ when equipped with appropriate weights, even though they are not orthonormal in $H^{-s}(U)$.
For every fixed $t\in[0,T]$ and $k\in\mathbb{N}$, testing the equation for $w$ against $e_k$ gives
\begin{align}\label{equation for well posedness of OU process}
    \partial_t\langle w(t),e_k\rangle_s&= \langle\partial_t w(t),e_k\rangle_s \nonumber\\
    &=\Phi'(\bar\rho)\langle\Delta w(t),e_k \rangle_s-\langle\nabla\cdot(\nu'(\bar\rho)w(t)),e_k\rangle_s\nonumber\\
    &=-\lambda_k\Phi'(\bar\rho) \langle w(t), e_k\rangle_s +\langle w(t),\nu'(\bar\rho)\cdot\nabla e_k\rangle_s,
\end{align}
where $\lambda_k$ is the eigenvalue corresponding to eigenfunction $e_k$ as in Example \ref{example of noise: eigenvalues of laplacian}.
When performing integration by parts, we used that $e_k=w(\cdot,t)=0$ on $\partial U$.
Furthermore, for the final term of \eqref{equation for well posedness of OU process}, using the expansion in \eqref{expansion of w in terms of e_k} gives
\begin{align*}
    \langle w(t),\nu'(\bar\rho)\cdot\nabla e_k\rangle_s=\sum_{j=1}^\infty w_j(t)\langle e_j,\nu'(\bar\rho)\cdot\nabla e_k\rangle_s.
\end{align*}
Since it is not the case that the off-diagonal terms in the above sum vanish, we see that this is a coupled system, where the time evolution of $w_k(t)$ depends on $w_j(t)$ for every $j\in\mathbb{N}$.
We deal with the evolution of the whole system $\{w_j(t)\}_{k\in\mathbb{N}}$ at once. 
Define for $t\in[0,T]$,
\[E(t):=\sum_{j=1}^\infty|w_j(t)|^2.\]
The definition of $E$ and equation \eqref{equation for well posedness of OU process} give that
\begin{align}\label{eq: equation for evolution of E}
    \partial_t E(t)=2\sum_{j=1}^\infty w_j(t)\partial_tw_j(t)&=-2\Phi'(\bar\rho)\sum_{j=1}^\infty\lambda_j|w_j(t)|^2+2\sum_{j=1}^\infty w_j(t)\sum_{k=1}^\infty w_k(t)\langle e_j(t),\nu'(\bar\rho)\cdot\nabla e_k\rangle_s\nonumber\\
    &\leq-2\Phi'(\bar\rho)\lambda_{min}E(t)+2\sum_{j=1}^\infty\sum_{k=1}^\infty w_j(t) w_k(t)\langle e_j(t),\nu'(\bar\rho)\cdot\nabla e_k\rangle_s,
\end{align}
where $\lambda_{min}$ is the smallest eigenvalue.
We need to deal with the second term.
Recall that the weighted functions $\psi_j:=\lambda_j^{s/2}e_j$ form an orthonormal basis in $H^s(U)$.
It then follows by applying Cauchy-Schwarz inequality followed by Bessel's inequality and Assumption \ref{assumption of constant boundary data and constant initial condition}, that for $s>\frac{d+2}{2}$,
\begin{align*}
    \left|\sum_{j=1}^\infty\sum_{j,k=1}^\infty w_j(t) w_k(t)\langle e_j(t),\nu'(\bar\rho)\cdot\nabla e_k\rangle_s\right|&=\left|\sum_{j=1}^\infty\sum_{j,k=1}^\infty\lambda_j^{-s/2} w_j(t) w_k(t)\langle \psi_j(t),\nu'(\bar\rho)\cdot\nabla e_k\rangle_s\right|\\
    &\hspace{-50pt}\leq \left(\sum_{j=1}^\infty\sum_{k=1}^\infty \lambda_j^{-s}w^2_j(t) w^2_k(t)\right)^{1/2}\left(\sum_{j=1}^\infty\sum_{k=1}^\infty |\langle \psi_j(t),\nu'(\bar\rho)\cdot\nabla e_k\rangle_s|^2\right)^{1/2}\\
    &\leq \lambda_{min}^{-s} E(t) \left(\sum_{k=1}^\infty \|\nu'(\bar\rho)\cdot\nabla e_k\|^2_{H^{-s}(U)}\right)^{1/2}\\
    &\leq \lambda_{min}^{-s}|\nu'(\bar\rho)| E(t) \left(\sum_{k=1}^\infty \|e_k\|^2_{H^{-s+1}(U)}\right)^{1/2},
\end{align*}
 where the sum on the right hand side is bounded by Lemma \ref{lemma sum of norm of e_k} in the Appendix, noting that here $s>\frac{d+2}{2}$.
Putting equation \eqref{eq: equation for evolution of E} and subsequent computation together, we have for the evolution of $E$ that there exists a running constant $c\in(0,\infty)$ such that
\[\partial_tE(t)\leq\left(-2\Phi'(\bar\rho)\lambda_{min}+2c\lambda_{min}^{-s}|\nu'(\bar\rho)|\right) E(t)\leq cE(t).\]

Gr\"onwall's inequality alongside the fact that $E(0)=0$ due to the initial condition of $w$, implies that
\[E(t)=0 \,\,\text{for every}\,\, t\in[0,T],\]
which in turn says that for every $j\in\mathbb{N}$, we have $\langle w(t),e_j\rangle_s=0$.
Using the fact that $\{e_k\}_{k\in\mathbb{N}}$ is an orthonormal basis of $L^2(U)$, we deduce that $w=0$ $\mathbb{P}-a.s.$ in $L^2([0,T];H^{-s}(U))$ for every $s>d/2$, which completes the proof.
\end{proof}
Now that we have proved well-posedness of the limiting equation $v$, let us turn our focus on formulating the central limit theorem result for the regularised equation.

Analogous to equation \eqref{eq: for v epsilon k}, define 
\begin{equation}\label{eq: definition of v n epsilon k with regularised equations}
    v^{n,\epsilon,K}:=\epsilon^{-1/2}(\rho^{n,\epsilon,K}-\bar\rho)
\end{equation}
for $\rho^{n,\epsilon,K}$ the weak solution to the regularised equation \eqref{equation for rho n epsilon k} in the sense of Proposition \ref{proposition on existence of weak solution clt and ldp chapter}, and $\bar\rho$ the solution of the regularised hydrodynamic limit equation \eqref{eq: regularised hydrodynamic limit}.\\
Just as the formal computation in equation \eqref{formal computation for convergence of v epsilon to v} shows, we expect that $v^{n,\epsilon,K}$ converges to $v$
as $\epsilon\to0$, and $v\in H^{-s}(U)$ is only distribution valued. 
So we would expect all $L^p(U\times[0,T])$-norms of $v^{n,\epsilon,K}$ to blow up as $\epsilon \to0, K\to\infty$.
This statement is quantified by the below.
\begin{proposition}\label{prop: blow up of lp norm of v n epsilon K}
Fix $K,n\in\mathbb{N}$ and $\alpha,\epsilon\in(0,1)$.
Let $v^{n,\epsilon,K}$ denote the weak solution to equation \eqref{eq: definition of v n epsilon k with regularised equations} above.
Under Assumptions \ref{assumption on noise}, \ref{assumption of constant boundary data and constant initial condition}, \ref{asm: new assumptions} and \ref{asm: assumptions for well posedenss of equation}, for every $p\geq2$ we have
\small
    \begin{align}\label{eq: estimate for lp space-time norm of v n epsilon k}
        \mathbb{E}\|v^{n,\epsilon,K}\|^p_{L^p(U\times[0,T])}\leq c T^2\left(\|\nabla\cdot F^K_2\|^{p/2}_{L^\infty(U)}+\|F^K_3\|^{p/2}_{L^\infty(U)}\right)\left(1+\epsilon T\left(\|\nabla\cdot F^K_2\|_{L^\infty(U)}+\|F^K_3\|
    _{L^\infty(U)}\right)\right).
    \end{align}
    \normalsize
\end{proposition}
\begin{proof}
    The proof is a direct consequence of the $p$-dependent energy estimate, Proposition \ref{prop: second energy estimate for rho n epsilon k}, by which it follows that
     \begin{multline*}
         \int_0^T\int_U (v^{n,\epsilon,K})^p=\epsilon^{-p/2}\int_0^T\int_U(\rho^{n,\epsilon,K}-\bar\rho)^p\\
         \leq c T^2\left(\|\nabla\cdot F^K_2\|^{p/2}_{L^\infty(U)}+\|F^K_3\|^{p/2}_{L^\infty(U)}\right)\left(1+\epsilon T\left(\|\nabla\cdot F^K_2\|_{L^\infty(U)}+\|F^K_3\|
    _{L^\infty(U)}\right)\right).
     \end{multline*}
\end{proof}

\begin{remark}[Divergence of {$L^p(U\times[0,T])$}-norms of $v^{n,\epsilon,K}$]
    In contrast to Remark \ref{rmk: choice of scaling in energy estimate} where we can find a joint scaling $\epsilon\to0, K\to\infty$ such that the $L^p(U\times[0,T])$-norms of $(\rho^{n,\epsilon,K}-\bar\rho)$ converge to zero, clearly the first term in the product  of \eqref{eq: estimate for lp space-time norm of v n epsilon k} is independent of $\epsilon$, so we have divergence of the $L^p(U\times[0,T])$-norms of $v^{n,\epsilon,K}$ under any scaling regime when the $L^\infty(U)$-norms of $\nabla\cdot F_2^K$ and $F_3^K$ diverge as $K\to\infty$.
    Furthermore, based on the joint scaling from Remark \ref{rmk: choice of scaling in energy estimate}, whenever $\epsilon,K$ satisfy $\epsilon\left(\|\nabla\cdot F^K_2\|_{L^\infty(U)}+\|F^K_3\|_{L^\infty(U)}\right)\leq 1$, the left hand side of \eqref{eq: estimate for lp space-time norm of v n epsilon k} diverges at the rate
   \[\left(\|\nabla\cdot F^K_2\|^{p/2}_{L^\infty(U)}+\|F^K_3\|^{p/2}_{L^\infty(U)}\right).\]
    For the choice of noise as in Example \ref{example of noise: eigenvalues of laplacian}, equations \eqref{eq: bound for divergence of F2 in laplacian eigenfunction example} and \eqref{eq: bounds for F_k in eigenvalue of laplacian case} in the Appendix show that the rate of divergence is solely governed by 
    \[\|F^K_3\|^{p/2}_{L^\infty(U)}\sim K^{p(d+2)}.\]
\end{remark}

We finish this section with the first main result of the work, proving the central limit theorem for the approximating SPDE.

Just as in equations \eqref{equation for rho n epsilon k} and \eqref{eq: regularised hydrodynamic limit}, for $\alpha\in(0,1),$ define the regularised equation $v$ as the strong solution to
\begin{equation}\label{eq: regularised OU process}
    \partial_t v=\Delta(\Phi'(\bar\rho)v))+\alpha\Delta v-\nabla\cdot(\sigma(\bar\rho)\dot\xi+\nu'(\bar\rho)v).
\end{equation} 
As usual we abused notation as initially $v$ denoted the solution to the above equation in the case $\alpha=0$.
We mentioned in Remark \ref{rmk: no need to smooth v} that the regularisation in \eqref{eq: regularised OU process} is not necessary to prove well-posedness in light of Assumption \ref{assumption of constant boundary data and constant initial condition}. 
However, it will be necessary when grouping terms below, see equation \eqref{eq: for z n epsilon K}.
Existence of strong solutions to the regularised equation \eqref{eq: regularised OU process} follows in the same way as Proposition \ref{proposition existence of strong solutions to OU process}.\\
Fundamentally, the below estimate is different from that of  Propositions \ref{ppn: energy estimate for difference of regularised equation and hydrodynamic limit}, \ref{prop: second energy estimate for rho n epsilon k} and \ref{prop: blow up of lp norm of v n epsilon K}, since neither quantity in the difference $v^{n,\epsilon,K}-v$ is solved by a constant, so we will see in the proof below that one really needs to handle the difference of the non-linear terms.
\begin{theorem}[CLT for approximating equation]\label{theorem CLT for approximating equation}
    Fix $\alpha,\epsilon\in(0,1)$ and $N,K\in\mathbb{N}$.
    Let $v^{n,\epsilon,K}$ denote the weak solution to the regularised equation \eqref{eq: definition of v n epsilon k with regularised equations}, and let $v$ be the solution of the regularised, linearised stochastic PDE \eqref{eq: regularised OU process} in the sense of Definition \ref{definition strong solution to ou process}.
    Furthermore, for fixed $s>\frac{d+2}{2}$, for every $K\in\mathbb{N}$, denote by $\mathcal{T}_s(K)$ the tail sum $\mathcal{T}_s(K):=\sum_{k=K}^\infty \|f_k\|_{H^{-s+1}(U)}$.
    Then under Assumptions \ref{assumption on noise}, \ref{assumption of constant boundary data and constant initial condition}, \ref{asm: new assumptions} and \ref{asm: assumptions for well posedenss of equation}, for $\beta\in(0,\infty)$ as in point 1 of Assumption \ref{asm: new assumptions}, there exists a constant $c\in(0,\infty)$ independent of $n$, $\epsilon$, $K$ and $\alpha$ such that 
    \begin{align}\label{eq: clt for approximate equation}
        \mathbb{E}&\|v^{n,\epsilon,K}-v\|^2_{L^2([0,T];H^{-s}(U))}\leq c\left(1+\epsilon T \left(\|\nabla\cdot F^K_2\|_{L^\infty(U)}+\|F^K_3\|
    _{L^\infty(U)}\right)\right)\nonumber\\
        &\times\left[\epsilon T^2\left(\|\nabla\cdot F^K_2\|^{2}_{L^\infty(U)}+\|F^K_3\|^{2}_{L^\infty(U)}\right)+\epsilon^{\beta+1} T^2\left(\|\nabla\cdot F^K_2\|^{\beta+2}_{L^\infty(U)}+\|F^K_3\|^{\beta+2}_{L^\infty(U)}\right)\right. \nonumber\\
        &\left.+\epsilon T\left(\|F_1^K\|^2_{L^\infty(U)}\|\sigma_n'\|^2_{L^\infty(U)}+\|F_2^K\|^2_{L^\infty(U;\mathbb{R}^d)}\right)\right] +c\mathbb{E}\|\sigma_n(\rho^{n,\epsilon,K})-\sigma(\bar\rho)\|_{L^2(U\times[0,T])}+cT\,\mathcal{T}_s(K).
    \end{align}
\end{theorem}

\begin{proof}
    For $s>\frac{d+2}{2}$ let $z^{n,\epsilon,K}:=(-\Delta)^{-s/2}\left(v^{n,\epsilon,K}-v\right)$.
    Since $s>\frac{d+2}{2}$, we have that $\mathbb{P}-$a.s. $z^{n,\epsilon,K}\in L^2([0,T];H^1(U))$ and satisfies
\begin{align}\label{eq: for z n epsilon K}
    \partial_t z^{n,\epsilon,K}&=\epsilon^{-1/2}(-\Delta)^{-s/2}(\Delta\Phi(\rho^{n,\epsilon,K})-\Delta\Phi(\bar\rho))-\Phi'(\bar\rho)(-\Delta)^{-s/2}\Delta v\nonumber\\
    &+\alpha\epsilon^{-1/2}(-\Delta)^{-s/2}(\Delta\rho^{n,\epsilon,K}-\Delta\bar\rho)-\alpha(-\Delta)^{-s/2}\Delta v\nonumber\\
    &+\epsilon^{-1/2}(-\Delta)^{-s/2}\left(\nabla\cdot\nu(\rho^{n,\epsilon,K})-\nabla\cdot\nu(\bar\rho)\right)-(-\Delta)^{-s/2}\nabla\cdot(\nu'(\bar\rho)v) \nonumber\\
       &-(-\Delta)^{-s/2}\nabla\cdot(\sigma_n(\rho^{n,\epsilon,K})\dot\xi^K)+(-\Delta)^{-s/2}\nabla\cdot(\sigma(\bar\rho)\dot\xi)\nonumber\\
       &+\frac{\epsilon^{1/2}}{2}(-\Delta)^{-s/2}\nabla\cdot(F^K_1(\sigma_n'(\rho^{n,\epsilon,K}))^2\nabla\rho^{n,\epsilon,K}+\sigma_n(\rho^{n,\epsilon,K})\sigma'_n(\rho^{n,\epsilon,K})F^K_2).
\end{align}
By It\^o formula, we have
\begin{equation} \label{ito formula for z n epsilon k}
    \int_0^t\int_Ud(z^{n,\epsilon,K})^2=2\int_0^t\int_U z^{n,\epsilon,K}\,dz^{n,\epsilon,K}+\int_0^t\int_U\,d\langle z^{n,\epsilon,K}\rangle.
\end{equation}
First we compute the quadratic variation term in \eqref{ito formula for z n epsilon k} and bound it right away. 
By Assumption \ref{assumption on noise} we have that $\dot\xi=\sum_{k\in\mathbb
{N}}f_kdB_t^k$, so we can re-write the noise terms in \eqref{eq: for z n epsilon K} as 
\begin{align*}
    &-(-\Delta)^{-s/2}\nabla\cdot(\sigma_n(\rho^{n,\epsilon,K})\dot\xi^K)+(-\Delta)^{-s/2}\nabla\cdot(\sigma(\bar\rho)\dot\xi)\\
    &=-\sum_{k=1}^K(-\Delta)^{-s/2}\nabla((\sigma_n(\rho^{n,\epsilon,K})-\sigma(\bar\rho))f_k)\cdot dB_t^k+\sum_{k=K+1}^\infty(-\Delta)^{-s/2}\nabla(\sigma(\bar\rho))f_k)\cdot dB_t^k.
\end{align*}
It follows by the $L^2(U)$-orthonormality of $\{f_k\}_{k\in\mathbb{N}}$ and the fact that $s>3/2$ that the first of the two terms can be bounded by
\begin{align*}
    \sum_{k=1}^K\int_0^t\int_U\left((-\Delta)^{-s/2}\nabla((\sigma_n(\rho^{n,\epsilon,K})-\sigma(\bar\rho))f_k)\right)^2\leq c\|\sigma_n(\rho^{n,\epsilon,K})-\sigma(\bar\rho)\|_{L^2(U\times[0,T])},
\end{align*}
and for the second, using Assumption \ref{assumption of constant boundary data and constant initial condition} alongside the definition of the $H^{-s}(U)$-norm gives
\begin{align*}
    \sum_{k=K+1}^\infty\int_0^t\int_U\left((-\Delta)^{-s/2}\nabla(\sigma(\bar\rho))f_k)\right)^2 &\leq \sigma^2(\bar\rho)t\sum_{k=K+1}^\infty\int_U\left((-\Delta)^{-s/2}\nabla f_k\right)^2\\
    &    \leq ct\mathcal{T}_s(K),
\end{align*}
where $\mathcal{T}_s(K)$ is defined in the statement of the theorem and in the final inequality we used $\mathcal{T}_s(K+1)\leq \mathcal{T}_s(K)$.
By point 2 of Assumption \ref{assumption on noise}, $\mathcal{T}_s(K)$ decays to zero as $K\to\infty$ with an explicit rate that depends on the specific choice of $\{f_k\}_{k\in\mathbb{N}}$, see Example \ref{example value of G(K) for eigenvalue of laplacian} below.\\
\begin{multline}\label{outcome of ito formula z n epsilon k}
       \frac{1}{2}\int_0^t\int_U(z^{n,\epsilon,K})^2\leq\int_0^t\int_Uz^{n,\epsilon,K}\left[\epsilon^{-1/2}(-\Delta)^{-s/2}(\Delta\Phi(\rho^{n,\epsilon,K})-\Delta\Phi(\bar\rho))-\Phi'(\bar\rho)(-\Delta)^{-s/2}\Delta v\right.\\
       +\alpha\epsilon^{-1/2}(-\Delta)^{-s/2}(\Delta\rho^{n,\epsilon,K}-\Delta\bar\rho)-\alpha(-\Delta)^{-s/2}\Delta v\\
    +\epsilon^{-1/2}(-\Delta)^{-s/2}\left(\nabla\cdot\nu(\rho^{n,\epsilon,K})-\nabla\cdot\nu(\bar\rho)\right)-(-\Delta)^{-s/2}\nabla\cdot(\nu'(\bar\rho)v)\\
       -(-\Delta)^{-s/2}\nabla\cdot(\sigma_n(\rho^{n,\epsilon,K})\dot\xi^K)+(-\Delta)^{-s/2}\nabla\cdot(\sigma(\bar\rho)\dot\xi)\\
       \left.+\frac{\epsilon^{1/2}}{2}(-\Delta)^{-s/2}\nabla\cdot(F^K_1(\sigma_n'(\rho^{n,\epsilon,K}))^2\nabla\rho^{n,\epsilon,K}+\sigma_n(\rho^{n,\epsilon,K})\sigma'_n(\rho^{n,\epsilon,K})F^K_2)\right]\\
       +c\|\sigma_n(\rho^{n,\epsilon,K})-\sigma(\bar\rho)\|_{L^2(U\times[0,T])}+\sigma^2(\bar\rho)t\,\mathcal{T}_s(K).
\end{multline}
We deal with each term of \eqref{outcome of ito formula z n epsilon k} in turn.\\
First of all, for the terms in the second line involving the $\alpha$-regularisation, they can be re-written as
\[2\alpha\int_0^t\int_Uz^{n,\epsilon,K}\Delta z^{n,\epsilon,K}=-2\alpha\int_0^t\int_U|\nabla z^{n,\epsilon,K}|^2,\]
which can be moved over to the left hand side of the estimate.
Next, the terms in the first line on the right hand side can be written as
\[\int_0^t\int_Uz^{n,\epsilon,K}\left((-\Delta)^{-s/2}\Delta\left(v^{n,\epsilon,K}\frac{\Phi(\rho^{n,\epsilon,K})-\Phi(\bar\rho)}{\rho^{n,\epsilon,K}-\bar\rho}-\Phi'(\bar\rho) v\right)\right).\]
We write the above term in the inner bracket as \say{what we want} plus a correction,
\begin{align}\label{eq: first term in clt approximate equation estimate}
    &\left(v^{n,\epsilon,K}\frac{\Phi(\rho^{n,\epsilon,K})-\Phi(\bar\rho)}{\rho^{n,\epsilon,K}-\bar\rho}-\Phi'(\bar\rho) v\right)=\Phi'(\bar\rho)(v^{n,\epsilon,K}-v)+v^{n,\epsilon,K}\left(\frac{\Phi(\rho^{n,\epsilon,K})-\Phi(\bar\rho)}{\rho^{n,\epsilon,K}-\bar\rho}-\Phi'(\bar\rho) \right).
\end{align}
Using integration by parts and Assumption \ref{assumption of constant boundary data and constant initial condition}, we get that the first term in \eqref{eq: first term in clt approximate equation estimate} gives
\begin{align}\label{eq: moving L2 norm of gradient of z to lhs}
    \int_0^t\int_U z^{n,\epsilon,K} (-\Delta)^{-s/2}\Delta \Phi'(\bar\rho)(v^{n,\epsilon,K}-v)=-\Phi'(\bar\rho)\int_0^t\int_U|\nabla z^{n,\epsilon,K}|^2,
\end{align}
which can be moved onto the left hand side of the estimate.
By Assumption \ref{assumption of constant boundary data and constant initial condition} and point 2 of Assumption \ref{asm: assumptions for well posedenss of equation} we have $\Phi'(\bar\rho)$ is strictly positive.
This allows us to absorb space-time integrals of $|\nabla z^{n,\epsilon,K}|^2$ to the left hand side in a way that does not depend on the regularisation $\alpha$.\\
We now aim to bound the final term corresponding to the correction in \eqref{eq: first term in clt approximate equation estimate}.
Fix $(x,t)\in U\times[0,T]$.
By the mean value theorem, we have that there exists a $\xi=\xi(x,t)\in(\bar\rho,\rho^{n,\epsilon,K})$ (if $\rho^{n,\epsilon,K}(x,t)>\bar\rho$), otherwise $\xi\in(\rho^{n,\epsilon,K},\bar\rho)$ such that 
\[\Phi(\rho^{n,\epsilon,K})-\Phi(\bar\rho)=\Phi'(\xi)(\rho^{n,\epsilon,K}-\bar\rho).\]
Consequently, substituting this into the final term of \eqref{eq: first term in clt approximate equation estimate} and applying the mean value theorem again gives for $\xi'=\xi'(x,t)\in(\bar\rho,\xi)$ (if $\xi(x,t)>\bar\rho)$ otherwise for $\xi'\in(\xi,\bar\rho)$,
\begin{align*}
   \left(\frac{\Phi(\rho^{n,\epsilon,K})-\Phi(\bar\rho)}{\rho^{n,\epsilon,K}-\bar\rho}-\Phi'(\bar\rho) \right)&=\Phi'(\xi)-\Phi'(\bar\rho) =\Phi''(\xi')(\xi-\bar\rho).
\end{align*}
Note that the final term exists due to point 1 of Assumption \ref{asm: new assumptions}.
Then, using  Cauchy-Schwarz and Young's inequalities, the fact that $s>\frac{d+2}{2}$, the bound $|\xi-\bar\rho|\leq |\rho^{n,\epsilon,K}-\bar\rho|$ and the bound on $\Phi''$ given in point 1 of Assumption \ref{asm: new assumptions} then proves that there exists constants $c,\beta\in(0,\infty)$ such that for every $\delta\in(0,1)$,
\begin{align}\label{eq: needing to bound l4 norm of v}
    \int_0^t\int_Uz^{n,\epsilon,K}&(-\Delta)^{-s/2}\Delta\left(v^{n,\epsilon,K}\Phi''(\xi')(\xi-\bar\rho)\right)\nonumber\\
    &\leq \frac{\delta}{2}\int_0^t\int_U (z^{n,\epsilon,K})^2+\frac{1}{2\delta}\int_0^T\int_U\left((-\Delta)^{-s/2}\Delta\left(v^{n,\epsilon,K}\Phi''(\xi')(\xi-\bar\rho)\right)\right)^2\nonumber\\
    & \leq \frac{\delta}{2}\int_0^t\int_U (z^{n,\epsilon,K})^2+\frac{c}{\delta}\int_0^T\int_U\left(v^{n,\epsilon,K}(\rho^{n,\epsilon,K}-\bar\rho)(1+(\rho^{n,\epsilon,K}-\bar{\rho})^\beta)\right)^2\nonumber\\
    &\leq \frac{\delta}{2}\int_0^t\int_U (z^{n,\epsilon,K})^2+\frac{c}{\delta}\int_0^T\int_U\left(\epsilon^{1/2}(v^{n,\epsilon,K})^2+\epsilon^{\frac{\beta+1}{2}}(v^{n,\epsilon,K})^{\beta+2}\right)^2\nonumber\\
    &\leq \frac{\delta}{2}\int_0^t\int_U (z^{n,\epsilon,K})^2+\frac{c\epsilon}{\delta}\int_0^T\int_U(v^{n,\epsilon,K})^4+\frac{c\epsilon^{\beta+1}}{\delta}\int_0^T\int_U(v^{n,\epsilon,K})^{2\beta+4}.
\end{align}
We can pick $\delta>0$ sufficiently small so that the first term can be absorbed onto the left hand side of \eqref{outcome of ito formula z n epsilon k}.
The terms involving $\nu$ in the third line of \eqref{outcome of ito formula z n epsilon k} can be handled in the same way as the above terms, again using point 1 of Assumption \ref{asm: new assumptions}. 
We have that there exists constants $c,\beta\in(0,\infty)$ such that for every $\delta>0$, $\xi,\xi'$ in the same intervals as above,
\small
\begin{align*}
    &\int_0^t\int_Uz^{n,\epsilon,K}\left[\epsilon^{-1/2}(-\Delta)^{-s/2}\left(\nabla\cdot\nu(\rho^{n,\epsilon,K})-\nabla\cdot\nu(\bar\rho)\right)-(-\Delta)^{-s/2}\nabla\cdot(\nu'(\bar\rho)v)\right]\\
    &=\int_0^t\int_U\nabla z^{n,\epsilon,K}\cdot\nu'(\bar\rho) z^{n,\epsilon,K}
    +\int_0^t\int_U z^{n,\epsilon,K}(-\Delta)^{-s/2}\nabla\cdot\left(v^{n,\epsilon,K}\left(\frac{\nu(\rho^{n,\epsilon,K})-\nu(\bar\rho)}{\rho^{n,\epsilon,K}-\bar\rho}-\nu'(\bar\rho)\right)\right)\\
    &=\int_0^t\int_U\nabla z^{n,\epsilon,K}\cdot\nu'(\bar\rho) z^{n,\epsilon,K}
    +\int_0^t\int_U z^{n,\epsilon,K}(-\Delta)^{-s/2}\nabla\cdot\left(v^{n,\epsilon,K}\nu''(\xi')(\xi-\bar\rho)\right)\\
    &\leq \int_0^t\int_U\nabla z^{n,\epsilon,K}\cdot\nu'(\bar\rho) z^{n,\epsilon,K}+c\delta\int_0^t\int_U|\nabla z^{n,\epsilon,K}|^2+\frac{c\epsilon}{\delta}\int_0^T\int_U(v^{n,\epsilon,K})^4+\frac{c\epsilon^{\beta+1}}{\delta}\int_0^T\int_U(v^{n,\epsilon,K})^{2\beta+4}.
\end{align*}
\normalsize
The first term on the right hand side in the final line vanishes due to Assumption \ref{assumption of constant boundary data and constant initial condition} in analogy with terms involving $\nu$ in previous estimates, see for instance equation \eqref{eq: nu vanishing on boundary}.
The second term can be absorbed onto the left hand side of the estimate by picking $\delta>0$ sufficiently small, independently of $\alpha$. We will deal with the remaining terms below.\\
The noise terms in the fourth line on the right hand side of \eqref{outcome of ito formula z n epsilon k} vanish under expectation.\\
For the first term in the penultimate line of \eqref{outcome of ito formula z n epsilon k}, we have using integration by parts, Cauchy-Schwarz inequality, the fact that $s>d/2$ and Young's inequality, that there exists a constant $c\in(0,\infty)$ such that for every $\delta\in(0,1)$,
\begin{align}\label{eq: problematic term in approximate CLT}
    &\frac{\epsilon^{1/2}}{2}\int_0^t\int_U z^{n,\epsilon,K} (-\Delta)^{-s/2}\nabla\cdot(F_1^K(\sigma_n'(\rho^{n,\epsilon,K}))^2\nabla\rho^{n,\epsilon,K})\nonumber\\
    &=-\frac{\epsilon^{1/2}}{2}\int_0^t\int_U \nabla z^{n,\epsilon,K}\cdot (-\Delta)^{-s/2}(F_1^K(\sigma_n'(\rho^{n,\epsilon,K}))^2\nabla\rho^{n,\epsilon,K})\nonumber\\
    &\leq\frac{\epsilon^{1/2}\|F_1^K\|_{L^\infty(U)}\|\sigma_n'\|_{L^\infty(U)}}{2}\int_0^t\int_U |\nabla z^{n,\epsilon,K}|\,|(-\Delta)^{-s/2}\nabla\sigma_n(\rho^{n,\epsilon,K})|\nonumber\\
    &\leq\frac{\delta}{4}\int_0^t\int_U|\nabla z^{n,\epsilon,K}|^2+\frac{\epsilon\|F_1^K\|^2_{L^\infty(U)}\|\sigma_n'\|^2_{L^\infty(U)}}{\delta}\int_0^t\int_U|\sigma_n(\rho^{n,\epsilon,K})|^2
\end{align}
Similarly, for the second term in the penultimate line, integration by parts, Cauchy-Schwarz and Young's inequalities imply that there exists a constant $c\in(0,\infty)$ such that for every $\delta\in(0,1)$,
\begin{align*}
    &\frac{\epsilon^{1/2}}{2}\int_0^t\int_U z^{n,\epsilon,K} (-\Delta)^{-s/2}\nabla\cdot(\sigma_n(\rho^{n,\epsilon,K})\sigma'_n(\rho^{n,\epsilon,K})F^K_2)\\
    &=-\frac{\epsilon^{1/2}}{2}\int_0^t\int_U \nabla z^{n,\epsilon,K} (-\Delta)^{-s/2}(\sigma_n(\rho^{n,\epsilon,K})\sigma'_n(\rho^{n,\epsilon,K})F^K_2)\\
    &\leq \frac{\delta}{4}\int_0^t\int_U|\nabla z^{n,\epsilon,K}|^2+\frac{\epsilon \|F_2^K\|^2_{L^\infty(U;\mathbb{R}^d)}}{\delta}\int_0^t\int_U(\sigma_n(\rho^{n,\epsilon,K})\sigma'_n(\rho^{n,\epsilon,K}))^2
\end{align*}

Again pick $\delta$ sufficiently small, independent of $\alpha$, so that the first term of the above equation and equation \eqref{eq: problematic term in approximate CLT} can be moved to the left hand side.
Putting everything together, it follows from \eqref{outcome of ito formula z n epsilon k} and subsequent computations that we have the bound
\small
\begin{multline}\label{eq: final equation in clt estimate}
       \mathbb{E}\left[\int_0^t\int_U\frac{1}{2}(z^{n,\epsilon,K})^2+(\Phi'(\bar\rho)+\alpha)|\nabla z^{n,\epsilon,K}|^2\right]\leq c\epsilon\mathbb{E}\left[\int_0^T\int_U(v^{n,\epsilon,K})^4\right]+c\epsilon^{\beta+1}\mathbb{E}\left[\int_0^T\int_U(v^{n,\epsilon,K})^{2\beta+4}\right]\\
       +\epsilon\|F_1^K\|^2_{L^\infty(U)}\|\sigma_n'\|^2_{L^\infty(U)}\mathbb{E}\left[\int_0^t\int_U|\sigma(\rho^{n,\epsilon,K})|^2\right]+\epsilon \|F_2^K\|^2_{L^\infty(U;\mathbb{R}^d)}\mathbb{E}\left[\int_0^t\int_U(\sigma_n(\rho^{n,\epsilon,K})\sigma'_n(\rho^{n,\epsilon,K}))^2\right]\\
    +c\mathbb{E}\|\sigma_n(\rho^{n,\epsilon,K})-\sigma(\bar\rho)\|_{L^2(U\times[0,T])}+\sigma^2(\bar\rho)t\,\mathcal{T}_s(K).
\end{multline}
\normalsize
 For the first two terms on the right hand side, we use Proposition \ref{prop: blow up of lp norm of v n epsilon K} and the resulting terms appear on the right hand side of the estimate.
 Importantly, these terms have factors of $\epsilon$ in front, implying the existence of a joint scaling in $\epsilon,K$ under which the terms vanish.

To estimate the terms in the second line on the right hand side, we use the bounds on 
$\sigma^2$ and $(\sigma\sigma')^2$ given in equations \eqref{eq: bound on sigma and so sigma squared} and \eqref{eq: bound on nu squared and sigma sigma' squared} in Assumption \ref{asm: assumptions for well posedenss of equation} in the Appendix alongside the assumed growth on $\Phi$ in equation \eqref{eq: growth bound on phi} which we alter in the same way as in Remark \ref{rmk: smuggling constant rho bar into the inequality for point 1 new asms} to give
\begin{multline*}
    \epsilon\|F_1^K\|^2_{L^\infty(U)}\|\sigma_n'\|^2_{L^\infty(U)}\mathbb{E}\left[\int_0^t\int_U|\sigma(\rho^{n,\epsilon,K})|^2\right]+\epsilon \|F_2^K\|^2_{L^\infty(U;\mathbb{R}^d)}\mathbb{E}\left[\int_0^t\int_U(\sigma_n(\rho^{n,\epsilon,K})\sigma'_n(\rho^{n,\epsilon,K}))^2\right]\\
    \leq c\epsilon\left(\|F_1^K\|^2_{L^\infty(U)}\|\sigma_n'\|^2_{L^\infty(U)}+\|F_2^K\|^2_{L^\infty(U;\mathbb{R}^d)}\right)\mathbb{E}\left[\int_0^t\int_U\left(1+(\rho^{n,\epsilon,K}-\bar\rho)^m\right)\right].
\end{multline*}
Using the $p$-independent energy estimate of Proposition \ref{ppn: energy estimate for difference of regularised equation and hydrodynamic limit} and Assumption \ref{assumption of constant boundary data and constant initial condition} gives that the above term is bounded above by
\[c\epsilon t\left(\|F_1^K\|^2_{L^\infty(U)}\|\sigma_n'\|^2_{L^\infty(U)}+\|F_2^K\|^2_{L^\infty(U;\mathbb{R}^d)}\right)\left(1+\epsilon t \left(\|\nabla\cdot F^K_2\|_{L^\infty(U)}+\|F^K_3\|
    _{L^\infty(U)}\right)\right).\]
Putting everything together then gives the estimate.
\end{proof}
  \begin{remark}[Simplifying the estimate]
   The above estimate (equation \eqref{eq: clt for approximate equation}) has a delicate balance of coefficients involving powers of $\epsilon$ and powers of the $L^\infty(U)$-norms of sums of noise coefficients, making it impossible to simplify further in a straightforward way.  
   It comes down to how we bound the two norms of $v^{n,\epsilon,K}$ in equation \eqref{eq: final equation in clt estimate}.
       Using Proposition \ref{prop: blow up of lp norm of v n epsilon K} followed by the facts that that we can upper bound the $L^\infty(U)$-norms of the noise coefficients by their squares and $\epsilon^{\beta+1}\leq \epsilon^{(\beta+2)/2}$ for every $\epsilon\in(0,1),\beta>1$, and the trivial inequality $(1+x)(x+x^{\frac{\beta+2}{2}})\leq 8(x+x^{\frac{\beta+4}{2}})$ we could have alternative stated the bound for these two terms as
\begin{align*}
    &\epsilon\mathbb{E}\left[\int_0^T\int_U(v^{n,\epsilon,K})^4\right]+\epsilon^{\beta+1}\mathbb{E}\left[\int_0^T\int_U(v^{n,\epsilon,K})^{2\beta+4}\right]\\
    &\leq c T^2\left(1+\epsilon T\left(\|\nabla\cdot F^K_2\|_{L^\infty(U)}+\|F^K_3\|
    _{L^\infty(U)}\right)\right)\\
    &\hspace{50pt}\times\left(\epsilon\left(\|\nabla\cdot F^K_2\|^{2}_{L^\infty(U)}+\|F^K_3\|^{2}_{L^\infty(U)}\right)+\epsilon^{\beta+1}\left(\|\nabla\cdot F^K_2\|^{\beta+2}_{L^\infty(U)}+\|F^K_3\|^{\beta+2}_{L^\infty(U)}\right)\right)\\
    &\leq cT^3\left(\epsilon\left(\|\nabla\cdot F^K_2\|^{2}_{L^\infty(U)}+\|F^K_3\|^{2}_{L^\infty(U)}\right)+\epsilon^{\frac{\beta+4}{2}}\left(\|\nabla\cdot F^K_2\|^{\beta+4}_{L^\infty(U)}+\|F^K_3\|^{\beta+4}_{L^\infty(U)}\right)\right).
\end{align*}
This makes the scaling between $\epsilon$ and $K$ more explicit.
   \end{remark}
\begin{remark}
     It is not possible to directly take the limit of the regularisation $n\to\infty$ in equation \eqref{eq: clt for approximate equation} and get a statement for the singular equation. 
     This is due to the presence of the term $\|\sigma_n'\|_{L^\infty(0,\infty)}$ which would diverge in the critical case when $\sigma_n$ are approximating the square root.
\end{remark}
 \begin{remark}[Necessity of Assumption \ref{assumption of constant boundary data and constant initial condition}]
       It was illustrated in Proposition 4.14 of \cite{popat2025well} that if the boundary data is non-constant, then one can only estimate $L^2(U\times[0,T])$-norms of the solutions $\rho^{n,\epsilon,K}$ of \eqref{generalised Dean-Kawasaki Equation Stratonovich with epsilon} rather than $L^p(U\times[0,T])-$norms for general $p\geq2$, see Proposition \ref{prop: l2 estimate from 4.14 of popat} below.
       Consequently we could only hope to obtain $L^2(U\times[0,T])$-estimates for $v^{n,\epsilon,K}$ rather than the result of Proposition \ref{prop: blow up of lp norm of v n epsilon K}.
       However, in the proof of the central limit theorem above, we saw in equation \eqref{eq: needing to bound l4 norm of v} that we needed to bound the $L^4(U\times[0,T])$-norm and $L^{2\beta+4}(U\times[0,T])$-norm of $v^{n,\epsilon,K}$ for $\beta\in(0,\infty)$. This further highlights the necessity of Assumption \ref{assumption of constant boundary data and constant initial condition}.
   \end{remark}
\begin{example}[Explicit bound for $\mathcal{T}_s(K)$ for eigenvalues of Laplacian]\label{example value of G(K) for eigenvalue of laplacian}
    For $k\in\mathbb{N}$, if $f_k=e_k$ are the eigenfunctions of the Laplacian as in Example \ref{example of noise: eigenvalues of laplacian}, the computation in Lemma \ref{lemma sum of norm of e_k} in the Appendix shows that by Weyl's law \cite{weyl1911asymptotische}, we have for $s>\frac{d+2}{2}$ the approximation
\[\mathcal{T}_s(K)=\sum_{k=K}^\infty\|f_k\|^2_{H^{-s+1}(U)}\sim\int_K^\infty k^{-2(s-1)/d}\,dk\sim K^{1-2(s-1)/d},\]
which is an explicit rate of decay for the tail sum in terms of $K$ (note that $1-2(s-1)d<0$ so the right hand side does indeed converge to zero as $K\to\infty$).\\
In general if we have the growth $\|f_k\|_{H^{-s}(U)}^2\sim k^{-\alpha}$ for $\alpha>1$, then the tail sum grows like $K^{1-\alpha}$.
\end{example}

\subsection{Central limit theorem for the stochastic PDE with singular coefficients}\label{sec: clt for singular equation}
In this section we first establish an $L^\infty(U\times[0,T])$-estimate for the solution $\rho^{\epsilon,K}$ of the singular equation \eqref{generalised Dean-Kawasaki Equation Stratonovich with epsilon}.
This estimate is used to establish the quantitative central limit theorem in probability.\\
In the Theorem below, the computations following equation \eqref{eq: bounding lp norm by l p minus 2 norm} follow closely to Theorem 3.9 of \cite{dirr2020conservative}, so we just present the main idea.
\begin{theorem}[{$L^\infty(U\times[0,T])$}-estimate for equation with singular coefficients]\label{thm: l infinity estimate for singular equation}
Let $T\in(0,\infty)$, $\epsilon\in(0,1)$, $K\in\mathbb{N}$, and under Assumption \ref{assumption of constant boundary data and constant initial condition}, let $\rho_0\equiv M>0$.
Assume also that Assumptions \ref{assumption on noise}, \ref{asm: new assumptions} and \ref{asm: assumptions for well posedenss of equation} hold.
Then, for stochastic kinetic solutions $\rho^{\epsilon,K}$ of \eqref{generalised Dean-Kawasaki Equation Stratonovich with epsilon}, we have that there exists constants $c,\gamma\in(0,\infty)$ independent of $\epsilon$ such that
\[\mathbb{E}\|(\rho^{\epsilon,K}-M)_-\|_{L^\infty(U\times[0,T])}\leq c\left(\epsilon\left(\|F_3^K\|_{L^\infty(U)}+\|\nabla\cdot F_2^K\|_{L^\infty(U)}\right)\right)^\gamma. \]
\end{theorem}
\begin{proof}
    By Assumption \ref{assumption of constant boundary data and constant initial condition}, the boundary data is also $\Phi(\rho^{\epsilon,K})|_{\partial U}=M$.\\
    In order to apply It\^o's formula, we will need to work with the regularised equation \eqref{equation for rho n epsilon k}.
    However, all of the estimates are completely stable with respect to $\alpha\in(0,1)$ and the smoothing $n\in\mathbb{N}$, so in the end we can just apply dominated convergence and take the limits $\alpha\to 0, n\to\infty$.
    For fixed $p\geq 2,n\in\mathbb{N},K\in\mathbb{N},\alpha\in(0,1),\epsilon\in(0,1)$ and $\rho^{n,\epsilon,K}$ as in \eqref{equation for rho n epsilon k}, It\^o's formula tells us
    \[d(\rho^{n,\epsilon,K}-M)_-^p=\alpha(\rho^{n,\epsilon,K}-M)_-^{p-1}d(\rho^{n,\epsilon,K}-M)+\frac{1}{2}p(p-1)(\rho^{n,\epsilon,K}-M)_-^{p-2} d\langle\rho^{n,\epsilon,K}-M\rangle.\]

This gives after integrating by parts and rearranging, that
\small
\begin{multline}\label{eq: ito formula for l infinity estimate}
    \frac{1}{p(p-1)}\int_U (\rho^{n,\epsilon,K}_t-M)_-^p+\int_0^t\int_U(\rho^{n,\epsilon,K}-M)_-^{p-2}\Phi'(\rho^{n,\epsilon,K})|\nabla\rho^{n,\epsilon,K}|^2+\alpha\int_0^t\int_U(\rho^{n,\epsilon,K}-M)_-^{p-2}|\nabla\rho^{n,\epsilon,K}|^2\\
=\int_0^t\int_U(\rho^{n,\epsilon,K}-M)_-^{p-2}\left(\sqrt{\epsilon}\sigma_n(\rho^{n,\epsilon,K})\nabla\rho^{n,\epsilon,K} \cdot d\xi^K +\nabla\rho^{n,\epsilon,K} \cdot\nu(\rho^{n,\epsilon,K})\right)\\
+\int_0^t\int_U(\rho^{n,\epsilon,K}-M)_-^{p-2}\left(\frac{\epsilon}{2} F^K_3\sigma_n^2(\rho^{n,\epsilon,K})+\frac{\epsilon}{2}\sigma_n(\rho^{n,\epsilon,K})\sigma_n'(\rho^{n,\epsilon,K})\nabla\rho^{n,\epsilon,K}\cdot F^K_2\right).
\end{multline}
\normalsize
The key observation to make is that, since $\rho^{n,\epsilon,K}$ is non-negative, the term $(\rho^{n,\epsilon,K}-M)_-$ is non-zero only when $\rho^{n,\epsilon,K}\in[0,M]$, which is a compact set.
On this set, since $\sigma\in C_{loc}([0,\infty))$, we have that uniformly in the $n$ regularisation, 
\begin{align}\label{eq: sigma bounded}
    \left|\sigma_n(\rho^{n,\epsilon,K})\mathbbm{1}_{(\rho^{n,\epsilon,K}-M)_-\neq 0}\right|\leq c.
\end{align}
We want to use this to handle the first term on the right hand side of \eqref{eq: ito formula for l infinity estimate}. 
For every $n\in\mathbb{N}$, there exists a constant $c\in(0,\infty)$ independent of $n$ such that for the unique function 
$\Psi_{\sigma_n,p}(\xi)$, 
\begin{align}\label{eq: bd for Theta sigma}
\Psi_{\sigma_n,p}(\xi):=\int_0^\xi (\xi'-M)_-^{p-2}\sigma_n(\xi')\,d\xi'\leq c\int_0^\xi (\xi'-M)_-^{p-2}=\frac{c}{p-1}(\xi-M)_-^{p-1},    
\end{align}
Hence for the first term on the right hand side, using integration by parts,
\[\int_0^t\int_U(\rho^{n,\epsilon,K}-M)_-^{p-2}\sigma_n(\rho^{n,\epsilon,K})\nabla\rho^{n,\epsilon,K} \cdot d\xi^K=-\int_0^t\int_U\Psi_{\sigma_n,p}(\rho^{n,\epsilon,K})\nabla\cdot d\xi^K.\]
Analogous reasoning to \eqref{eq: no boundary terms when IBP F_2} illustrates that we do not pick up additional boundary terms in the case $p=2$.
The Burkholder–Davis–Gundy inequality, H\"older’s inequality and Young’s inequality prove that there exists a constant $c\in(0,\infty)$ such that for every $\delta\in(0, 1)$, 
\begin{align}\label{eq: handling noise term in l infinity estimate}
    \mathbb{E}&\left[\sqrt{\epsilon}\max_{t\in[0,T]}\left|\int_0^t\int_U(\rho^{n,\epsilon,K}-M)_-^{p-2}\sigma_n(\rho^{n,\epsilon,K})\nabla\rho^{n,\epsilon,K} \cdot d\xi^K\right|\right]\nonumber\\
    &\leq \frac{c\sqrt{\epsilon}}{p-1}\mathbb{E}\left[\left(\int_0^T\sum_{i=1}^d\sum_{k=1}^K\left(\int_U(\rho^{n,\epsilon,K}-M)_-^{p-1}\partial_i f_k\right)^2\right)^{1/2}\right]\nonumber\\
    &\leq c\sqrt{\epsilon}\mathbb{E}\left[\left(\sup_{t\in[0,T]}\int_U(\rho^{n,\epsilon,K}-M)_-^{p}\right)^{1/2}\left(\int_0^T\int_U(\rho^{n,\epsilon,K}-M)_-^{p-2}F^K_3\right)^{1/2}\right]\nonumber\\
    &\leq \frac{\delta}{2}\mathbb{E}\left[\sup_{t\in[0,T]}\int_U(\rho^{n,\epsilon,K}-M)_-^{p}\right]+\frac{c\epsilon\|F^K_3\|_{L^\infty(U)}}{2\delta}\mathbb{E}\left[\int_0^T\int_U(\rho^{n,\epsilon,K}-M)_-^{p-2}\right].
\end{align}
We take $\delta$ sufficiently small so that the first term can be absorbed onto the left hand side (after taking supremum in time).
Let us show how to deal with the other terms on the right hand side of equation \eqref{eq: ito formula for l infinity estimate}. 
The term involving $\nu$ vanishes due to Assumption \ref{assumption of constant boundary data and constant initial condition} by the same reasoning as usual, see for example equation \eqref{eq: nu vanishing on boundary}.
For the first term in the final line of \eqref{eq: ito formula for l infinity estimate}, we again use equation \eqref{eq: sigma bounded} which allows us to obtain for constant $c\in(0,\infty)$ independent of $n$,
\begin{align*}
   \frac{\epsilon}{2} \int_0^t\int_U(\rho^{n,\epsilon,K}-M)_-^{p-2} F^K_3\sigma_n^2(\rho^{n,\epsilon,K})\leq c\epsilon\|F^K_3\|_{L^\infty(U)}\int_0^t\int_U(\rho^{n,\epsilon,K}-M)_-^{p-2},
\end{align*}
Due to the assumption $\left(\sigma(\xi)\sigma'(\xi)\right)^2\leq c(1+\xi^m)$ which can be deduced from equations \eqref{eq: growth bound on phi} and \eqref{eq: bound on nu squared and sigma sigma' squared} from Assumption \ref{asm: assumptions for well posedenss of equation} in the Appendix, and the trivial observation that if $\xi\in[0,M]$ then $\xi^m\leq M^m$ is bounded, we have using integration by parts that for the final term on the right hand side of \eqref{eq: ito formula for l infinity estimate}, there exists a constant $c\in(0,\infty)$ such that
\begin{multline*}
\frac{\epsilon}{2}\int_0^t\int_U(\rho^{n,\epsilon,K}-M)_-^{p-2}\sigma(\rho^{n,\epsilon,K})\sigma'(\rho^{n,\epsilon,K})\nabla\rho^{n,\epsilon,K}\cdot F_2^K\leq c\epsilon\int_0^t\int_U(\rho^{n,\epsilon,K}-M)_-^{p-2}\nabla\rho^{n,\epsilon,K}\cdot F_2^K\\
\leq c\epsilon\|\nabla\cdot F^K_2\|_{L^\infty(U)}\int_0^t\int_U (\rho^{n,\epsilon,K}-M)_-^{p-1}.
\end{multline*}
As above, by an analogous computation to \eqref{eq: no boundary terms when IBP F_2}, in the case $p=2$ we don't pick up boundary terms when integrating by parts.
This term can then handled in the same way as the noise term \eqref{eq: handling noise term in l infinity estimate}, which can be easily seen by multiplying the integrand by $1$ and applying H\"older's inequality.

Putting equation \eqref{eq: ito formula for l infinity estimate} together with the subsequent computations, and taking supremum in time and expectation gives
\small
\begin{align}\label{eq: bounding lp norm by l p minus 2 norm}
       &\frac{1}{p(p-1)}\mathbb{E}\left[\sup_{t\in[0,T]}\int_U (\rho^{n,\epsilon,K}_t-M)_-^p\right]
       +\mathbb{E}\left[\int_0^T\int_U(\rho^{n,\epsilon,K}-M)_-^{p-2}\Phi'(\rho^{n,\epsilon,K})|\nabla\rho^{\epsilon,K}|^2\right]\nonumber\\
&+\alpha\mathbb{E}\left[\int_0^T\int_U\left|\nabla(\rho^{n,\epsilon,K}-M)_-^{p/2}\right|^2\right]\leq c\epsilon\left(\|F_3^K\|_{L^\infty(U)}+\|\nabla\cdot F_2^K\|_{L^\infty(U)}\right)\mathbb{E}\left[\int_0^T\int_U(\rho^{n,\epsilon,K}-M)_-^{p-2}\right].
\end{align}
\normalsize
The right hand side of the above estimate is independent of $\alpha$ and $n$, as it does not depend on the regularised $\sigma_n$. Hence, by dominated convergence, we take the limits as $\alpha\to0, n\to\infty$, which is represented in the below by replacing $\rho^{n,\epsilon,K}$ with $\rho^{\epsilon,K}$.
Furthermore, it follows by \cite{revuz2013continuous} Chapter 4, Proposition 4.7 and Exercise 4.30, and H\"older's inequality, that for every $p\in(1,\infty)$ and $n_p:=p^{-1}$ that we can raise both sides to the power $n_p$ at the cost of picking up an additional constant,
\begin{align*}
       &\mathbb{E}\left[\left(\sup_{t\in[0,T]}\int_U (\rho^{\epsilon,K}_t-M)_-^p+p(p-1)\int_0^T\int_U(\rho^{\epsilon,K}-M)_-^{p-2}\Phi'(\rho^{\epsilon,K})|\nabla\rho^{\epsilon,K}|^2\right.\right.\\
       &\hspace{180pt}\left.\left.+\alpha p(p-1)\int_0^T\int_U\left|\nabla(\rho^{\epsilon,K}-M)_-^{p/2}\right|^2\right)^{1/p}\right]
       \nonumber\\
&\leq \frac{n_p^{-n_p}}{1-n_p}\left(cp^2\epsilon\left(\|F_3^K\|_{L^\infty(U)}+\|\nabla\cdot F_2^K\|_{L^\infty(U)}\right)\right)^{1/p}\mathbb{E}\left[\|(\rho^{\epsilon,K}-M)_-\|_{L^{p-2}(U\times[0,T])}\right]^{\frac{p-2}{p}},
\end{align*} 
where on the right hand side we used the bound $p(p-1)<p^2$.
We now want to use a Moser iteration to conclude. 
The above bound, bounding the $L^p(U)$-norm by the $L^{p-2}(U)$-norm will not be sufficient for us, because if we used that to iterate, the constants that we accumulate would diverge. 
Instead, we we will need the below exponential bound.
By the same computation as in Theorem 3.9 of \cite{dirr2020conservative}, using the log-convex inequality, the Sobolev inequality and H\"older's inequality gives the improved estimate
\begin{multline}\label{eq: bound for negaitve part of rho - M}
    \mathbb{E}\|(\rho^{\epsilon,K}-M)_-\|_{L^\frac{(2+d)p}{d}(U\times [0,T])}\\
    \leq \frac{n_p^{-n_p}}{1-n_p}\left(c p^2\epsilon\left(\|F_3^K\|_{L^\infty(U)}+\|\nabla\cdot F_2^K\|_{L^\infty(U)}\right)\right)^{1/p}\mathbb{E}\left[\|(\rho^{\epsilon,K}-M)_-\|_{L^{p-2}(U\times[0,T])}\right]^{\frac{p-2}{p}}.
\end{multline}

Now, in the Mozer iteration below, at each step we pick up a multiplicative factor of $\frac{2+d}{d}$ that ensures the constants remain bounded in the limit.
Proceed inductively, define $p_0=0$ and for $k\in\mathbb{N}$ define iteratively
\[p_k=\frac{2+d}{d}p_{k-1}+2.\]
Trivially the bound $p_k>\left(\frac{d+2}{d}\right)^k$ is satisfied.
Inequality \eqref{eq: bound for negaitve part of rho - M} gives that for $q_k:=p_{k-1}+2$,
\small
\begin{multline}\label{eq: mozer iteration}
   \mathbb{E}\|(\rho^{\epsilon,K}-M)_-\|_{L^{p_k}(U\times [0,T])}\\
   \leq \frac{n_{q_k}^{-n_{q_k}}}{1-n_{q_k}}\left(c\epsilon q_k^2\left(\|F_3^K\|_{L^\infty(U)}+\|\nabla\cdot F_2^K\|_{L^\infty(U)}\right)\right)^{n_{q_k}}\mathbb{E}\left[\|(\rho^{\epsilon,K}-M)_-\|_{L^{p_{k-1}}(U\times[0,T])}\right]^{\frac{q_k}{q_k+2}}\\
   \leq \prod_{r=1}^{k}\left(\frac{n_{q_r}^{-n_{q_{r}}}}{1-n_{q_{r}}}\left(c q_{r}\right)^{2n_{q_r}}\right)^{\prod_{s=r+1}^{k}\frac{q_s}{q_{s}+2}}\left(\epsilon \left(\|F_3^K\|_{L^\infty(U)}+\|\nabla\cdot F_2^K\|_{L^\infty(U)}\right)\right)^{\sum_{r=1}^{k}n_{q_{r}}\prod_{s=r+1}^{k}\frac{q_{s}}{q_{s}+2}}.
\end{multline}
\normalsize
We conclude by showing that the constant and exponent in the above equation do not diverge in the $k\to\infty$ limit.
We begin analysing what happens to the exponent on the final term of \eqref{eq: mozer iteration}.
First of all, since we have $q_k>\left(\frac{d+2}{d}\right)^{k-1}+2$, and the fact that $\log(1-x)\approx -x$ for small $x>0$, we have
\begin{align*}
    \liminf_{N\to\infty}\log\left(\prod_{s=2}^{N}\frac{q_s}{q_s+2}\right)&=\liminf_{N\to\infty}\sum_{s=2}^{N}\log\left(1-\frac{2}{q_s+2}\right)\\
    &\approx \liminf_{N\to\infty}\sum_{s=2}^{N}-\frac{2}{q_s+2}> \liminf_{N\to\infty}\sum_{s=2}^{N}-2\left(\frac{d}{d+2}\right)^{s-1}.
\end{align*}
We note that the right hand side is a geometric series with ratio $r=\frac{d}{d+2}<1$ and so the series converges, and hence the left hand side is bigger than $-\infty$.\\
Furthermore, since $1/q_r$ decays exponentially due to the bound on $q_r$, we have that there is a constant $\gamma\in(0,\infty)$ that bounds the exponent on the final term of \eqref{eq: mozer iteration},
\[\lim_{k\to\infty}\left(\sum_{r=1}^{k}n_{q_{r}}\prod_{s=r+1}^{k}\frac{q_{s}}{q_{s}+2}\right)=\gamma.\]
Now we need to analyse the product of constants on the right hand side of \eqref{eq: mozer iteration}.
The definition of $n_{q_r}$ tells us that
\begin{align*}
    \limsup_{k\to\infty}\prod_{r=1}^{k}\left(\left(c q_{r}\right)^{2n_{q_r}}\frac{n_{q_r}^{-n_{q_{r}}}}{1-n_{q_{r}}}\right)^{\prod_{s=r}^{k-1}\frac{q_{s}}{q_{s}+2}}\leq\limsup_{k\to\infty}\prod_{r=1}^{k}\left(\left(c q_{r}\right)^{2n_{q_r}}\frac{q_r^{1+1/q_r}}{q_r-1}\right),
\end{align*}
and so for every $k\in\mathbb{N}$,
\begin{align*}
    \log\left(\prod_{r=1}^{k}\left(c q_{r}\right)^{2n_{q_r}}\frac{q_r^{1+1/q_r}}{q_r-1}\right)&=\sum_{r=1}^{k}2/q_r\log(cq_r)+1/q_r\log(q_r)+\log\left(\frac{q_r}{q_r-1}\right)\\
    &=\sum_{r=1}^{k}2/q_r\log(cq_r)+1/q_r\log(q_r)-\log(1+(q_r-1)^{-1}).
\end{align*}
Since  $q_r>\left(\frac{d+2}{d}\right)^{r-1}$, we get that
\[\limsup_{k\to\infty}\sum_{r=1}^{k}2n_{q_r}\log(cq_r)+1/q_r\log(q_r)-\log(1+(q_r-1)^{-1})<\infty.\]
This implies that there is a constant $\tilde{c}\in(0,\infty)$ such that
\[ \limsup_{k\to\infty}\prod_{r=1}^{k}\left(c\frac{n_{q_r}^{-n_{q_{r}}}}{1-n_{q_{r}}}\right)^{\prod_{s=r+1}^{k}\frac{q_{s}}{q_{s}+2}}\leq \tilde{c}.\]
Therefore, passing to the limit as $k\to\infty$ in \eqref{eq: mozer iteration}, we have the existence of constants $c,\gamma\in(0,\infty)$ such that
\[ \mathbb{E}\|(\rho^{\epsilon,K}_t-M)_-\|_{L^{\infty}(U\times [0,T])} \leq c\left(\epsilon\left(\|F_3^K\|_{L^\infty(U)}+\|\nabla\cdot F_2^K\|_{L^\infty(U)}\right)\right)^\gamma,\]
as required.
\end{proof}
We are now in a position to prove the central limit theorem for the equation with singular coefficients.
The proof combines the central limit theorem result for the regularised equation, theorem \ref{theorem CLT for approximating equation}, with the uniform estimate above.
We note that the methods are identical to the proof of Theorem 3.10 of \cite{dirr2020conservative}.
We provide the proof for completeness.
    \begin{theorem}[Central limit theorem for singular equation]\label{thm: CLT for singular equation}
    Let $K\in\mathbb{N}$, $\epsilon\in(0,1)$ and suppose Assumptions \ref{assumption on noise}, \ref{assumption of constant boundary data and constant initial condition}, \ref{asm: new assumptions} and \ref{asm: assumptions for well posedenss of equation} hold with $\beta\in(0,\infty)$ as in point 1 of Assumption \ref{asm: new assumptions}.
        Let $\rho^{\epsilon,K}$ be a stochastic kinetic solution of equation \eqref{generalised Dean-Kawasaki Equation Stratonovich with epsilon} and $\bar\rho$ be the solution of the hydrodynamic limit equation \eqref{heat equation hydrodynamic limit}.
        Let $v^{\epsilon,K}$ be a weak solution to \eqref{eq: for v epsilon k}, and $v$ a strong solution to in the linearised SPDE \eqref{eq: linearised spde v} sense of Definition \ref{definition strong solution to ou process}. 
        Then for every $s>\frac{d+2}{2}$
        there exists constants $c,\gamma\in(0,\infty)$ such that for every $\eta\in(0,1)$,
        \begin{multline}\label{eq: CLT for singular equation}
            \mathbb{P}\left(\|v^{\epsilon,K}-v\|_{L^2([0,T];H^{-s}(U))}\geq \eta\right)\leq c\eta^{-2}\left(1+\epsilon T\left(\|\nabla\cdot F^K_2\|_{L^\infty(U)}+\|F^K_3\|
    _{L^\infty(U)}\right)\right)\\
        \times\left[\epsilon T^2\left(\|\nabla\cdot F^K_2\|^{2}_{L^\infty(U)}+\|F^K_3\|^{2}_{L^\infty(U)}\right)+\epsilon^{\beta+1} T^2\left(\|\nabla\cdot F^K_2\|^{\beta+2}_{L^\infty(U)}+\|F^K_3\|^{\beta+2}_{L^\infty(U)}\right)\right. \\
        \left.+\epsilon T\left(\|F_1^K\|^2_{L^\infty(U)}+\|F_2^K\|^2_{L^\infty(U;\mathbb{R}^d)}\right)\right]+c\eta^{-2}T\,\mathcal{T}_s(K)+ c\left(\epsilon\left(\|F_3^K\|_{L^\infty(U)}+\|\nabla\cdot F_2^K\|_{L^\infty(U)}\right)\right)^\gamma.
        \end{multline}
\end{theorem}
\begin{proof}
By Assumption \ref{assumption of constant boundary data and constant initial condition}, the initial condition of $\rho^{\epsilon,K}$ is the constant $\rho_0=M>0$.
Furthermore, let $\tilde\sigma\in C([0,\infty))\cap C^\infty((0,\infty))$ be an arbitrary function satisfying for a constant $c\in(0,\infty)$ potentially depending on $M$\footnote{  In the square root diffusion case we would have $|\tilde{\sigma}'(\xi)|\leq \sqrt{2}M^{-1/2}$.},
    \begin{align}\label{def: tilde sigma}
        \tilde\sigma(\xi)=\sigma(\xi) \,\text{for every}\, \xi\in[M/2,\infty),\, \text{and}\, |\tilde\sigma'(\xi)|\leq c\,\text{for every $\xi\in(0,\infty)$}.
    \end{align}
   
Let $\tilde\rho^{\epsilon,K}$ be the unique stochastic kinetic solution to \eqref{generalised Dean-Kawasaki Equation Stratonovich with epsilon} with $\sigma$ replaced by $\tilde\sigma$.
Due to pathwise uniqueness of solutions, Theorem 3.6 of \cite{popat2025well}, the two solutions $\rho^{\epsilon,K}$, and $\tilde\rho^{\epsilon,K}$ coincide on the event when $\sigma$ and $\tilde\sigma$ coincide.
That is,
\[\rho^{\epsilon,K}=\tilde\rho^{\epsilon,K}\quad \text{in}\quad L^1(U\times[0,T]),\]
on the event $(\mathcal{S}\cap \tilde{S})\subseteq \Omega$ defined by
\[\mathcal{S}=\{\|\rho^{\epsilon,K}\|_{L^\infty(U\times[0,T])}\geq M/2\}\]
and
\[\tilde{\mathcal{S}}=\{\|\tilde\rho^{\epsilon,K}\|_{L^\infty(U\times[0,T])}\geq M/2\}.\]
Furthermore, define $\tilde v^{\epsilon,K}:=\epsilon^{-1/2}(\tilde\rho^{\epsilon,K}-\bar\rho)$, and let $v$ denote the solution to the linearised equation \eqref{eq: linearised spde v}.
We have the upper bound
\begin{equation}\label{eq: upper bound for clt}
    \mathbb{P}\left(\|v^\epsilon-v\|_{L^2([0,T];H^{-s}(U))}\geq \eta\right)\leq\mathbb{P}\left(\|\tilde{v}^\epsilon-v\|_{L^2([0,T];H^{-s}(U))}\geq \eta\right)+\mathbb{P}(\mathcal{S}^c)+\mathbb{P}(\tilde{\mathcal{S}}^c).
\end{equation}
For the first term, by Chebyshev's inequality and the central limit theorem for the approximate equation Theorem
\ref{theorem CLT for approximating equation}, replacing the computation in \eqref{eq: problematic term in approximate CLT} by incorporating the bound uniform bound \eqref{def: tilde sigma} on $\sigma'$, we get
\begin{multline}\label{eq: equation deriving CLT in prob for singular eq}
    \mathbb{P}\left(\|\tilde{v}^{\epsilon,K}-v\|_{L^2([0,T];H^{-s}(U))}>\eta\right)\leq c\eta^{-2}\left(1+\epsilon T\left(\|\nabla\cdot F^K_2\|_{L^\infty(U)}+\|F^K_3\|
    _{L^\infty(U)}\right)\right)\\
        \times\left[\epsilon T^2\left(\|\nabla\cdot F^K_2\|^{2}_{L^\infty(U)}+\|F^K_3\|^{2}_{L^\infty(U)}\right)+\epsilon^{\beta+1} T^2\left(\|\nabla\cdot F^K_2\|^{\beta+2}_{L^\infty(U)}+\|F^K_3\|^{\beta+2}_{L^\infty(U)}\right)\right. \\
        \left.+\epsilon T\left(\|F_1^K\|^2_{L^\infty(U)}+\|F_2^K\|^2_{L^\infty(U;\mathbb{R}^d)}\right)\right]+c\eta^{-2}T\,\mathcal{T}_s(K).
\end{multline}
We simplified the above estimate using the uniform bound on $\sigma$ from equation \eqref{def: tilde sigma}, which allows us to obtain using Proposition \ref{ppn: energy estimate for difference of regularised equation and hydrodynamic limit} with $p=2$ that
\begin{align*}
   \mathbb{E} \|\sigma(\tilde\rho^{\epsilon,K})-\sigma(\bar\rho)\|_{L^2(U\times[0,T])}&\leq c\mathbb{E}\|\tilde\rho^{\epsilon,K}-\bar\rho\|_{L^2(U\times[0,T])}\leq c\epsilon T^2 \left(\|\nabla\cdot F^K_2\|_{L^\infty(U)}+\|F^K_3\|_{L^\infty(U)}\right),
\end{align*}
where the right hand side is already present on the right hand side of the estimate \eqref{eq: equation deriving CLT in prob for singular eq}.
For the final two terms of \eqref{eq: upper bound for clt}, it follows from the $L^\infty(U\times[0,T])$-estimate of Theorem \ref{thm: l infinity estimate for singular equation}, which equally applies to $\tilde\rho^{\epsilon,K}$ due to \eqref{def: tilde sigma}, that
\[\mathbb{P}(\mathcal{S}^c)+\mathbb{P}(\tilde{\mathcal{S}}^c)\leq c\left(\epsilon\left(\|F_3^K\|_{L^\infty(U)}+\|\nabla\cdot F_2^K\|_{L^\infty(U)}\right)\right)^\gamma.\]
Putting both parts of estimate \eqref{eq: upper bound for clt} together finishes the proof.
\end{proof}
\begin{remark}
   We point out that the convergence in probability in Theorem \ref{thm: CLT for singular equation} above is weaker than the $L^1(\Omega;L^2([0,T];H^{-s}(U)))$-convergence for the regularised equation stated in Theorem \ref{theorem CLT for approximating equation}.
   We can not hope to achieve the convergence in expectation for the singular equation due to a lack of control of the expectation
    \[\mathbb{E}\left[\|v^{\epsilon,K}-v\|^2_{L^2([0,T];H^{-s}(U))}\mathbbm{1}_{\mathcal{S}^c}\right].\]
    \end{remark}
\begin{remark}[Choice of joint scaling]\label{rmk: choice of joint scaling}
    Based on the previous estimates (Propositions \ref{ppn: energy estimate for difference of regularised equation and hydrodynamic limit} and \ref{prop: second energy estimate for rho n epsilon k}, Theorems \ref{theorem CLT for approximating equation}, \ref{thm: l infinity estimate for singular equation} and \ref{thm: CLT for singular equation}), we choose a joint scaling such that the right hand side of the estimates converge to zero as $\epsilon\to0$ and $K\to\infty$. 
    A sufficient joint scaling is one which ensures that as $\epsilon\to0,K\to\infty$ we have
    \begin{equation}\label{eq: choice of joint scaling}
        \epsilon\left(\|F_1^K\|^2_{L^\infty(U)}+\|F_2^K\|^2_{L^\infty(U;\mathbb{R}^d)}+\|F_3^K\|^2_{L^\infty(U)}+\|\nabla\cdot F_2^K\|^2_{L^\infty(U)}\right)\to0.
    \end{equation}
    In the case of eigenvalues of Laplacian from Example \ref{example of noise: eigenvalues of laplacian}, equation \eqref{eq: bounds for F_k in eigenvalue of laplacian case} in the Appendix gives explicitly the scaling relation needed between $\epsilon$ and $K$ as
    \[\epsilon K^{d+2}\to0.\]
\end{remark}
\section{The large deviations principle}\label{sec: ldp}
The goal of this section is to prove that under the joint scaling regime $\epsilon\to0, K(\epsilon)\to\infty$ of Remark \ref{rmk: choice of joint scaling}, solutions $\rho^{\epsilon, K}$ of the generalised Dean--Kawasaki equation with truncated noise \eqref{generalised Dean-Kawasaki Equation Stratonovich with epsilon},
\[\partial_t\rho^{\epsilon,K}=\Delta\Phi(\rho^{\epsilon,K})-\sqrt{\epsilon}\nabla\cdot(\sigma(\rho^{\epsilon,K})\circ\dot{\xi}^K)-\nabla\cdot \nu(\rho^{\epsilon,K})\]
satisfy a large deviations principle, and also to identify the corresponding rate function.
As mentioned in the final part of Section \ref{sec: application to particle system}, the motivation is that one can study the large deviations principle for the zero range particle system by looking at the large deviations of the SPDE.\\
From Remark \ref{rmk: choice of joint scaling}, we know that $K=K(\epsilon)$ can be chosen entirely as a function of $\epsilon$ (and integer valued via the relabelling $K(\epsilon):=\lfloor K(\epsilon)\rfloor$), so for convenience of notation, in this section we denote solutions of \eqref{generalised Dean-Kawasaki Equation Stratonovich with epsilon} by
\[\rho^{\epsilon,K}\equiv \rho^\epsilon.\]
We begin by stating what it means for the solutions of the SPDE $\rho^{\epsilon}$ to satisfy the large deviations principle.
From the well-posedness results of \cite{popat2025well}, we know that solutions $\rho^{\epsilon}$ live in the space $L^1(U\times[0,T])$.
The theory of large deviations is concerned with events $A$ in the Borel set $\mathcal{B}(L^1(U\times[0,T]))$ such that the probability $\mathbb{P}(\rho^{\epsilon}\in A)$ decays exponentially under the joint scaling of Remark \ref{rmk: choice of joint scaling}.
The exponential rate is quantified in terms of a \say{rate function} $I:L^1(U\times[0,T])\to\mathbb{R}$.
The below definitions are standard, for example see Definitions 1 and 2 of \cite{budhiraja2008large}.
\begin{definition}[Rate function]\label{def: rate function}
    A function $I:L^1(U\times[0,T])\to\mathbb{R}$ is called a rate function on $L^1(U\times[0,T])$ if for every constant $M<\infty$, the level set 
    \[\{\rho\in L^1(U\times[0,T]):I(\rho)\leq M\}\]
    is a compact subset of $L^1(U\times[0,T])$.
    For $A\in \mathcal{B}(L^1(U\times[0,T]))$, define $I(A):=\inf_{\rho\in A}I(f)$.
    \end{definition}
\begin{definition}[Large deviations principle]\label{def: ldp}
    Let $I$ be a rate function on $L^1(U\times[0,T])$. 
    The sequence $\{\rho^{\epsilon}\}_{\epsilon\in(0,1)}$ satisfies the large deviation principle on $L^1(U\times[0,T])$ if the following two conditions hold
    \begin{enumerate}
        \item Large deviation upper bound: For each closed subset $F$ of $L^1(U\times[0,T])$,
        \[\limsup_{\epsilon\to0}\epsilon\log\mathbb{P}(\rho^\epsilon\in F)\leq -I(F).\]
        \item Large deviation lower bound:  For each open subset $G$ of $L^1(U\times[0,T])$,
        \[\liminf_{\epsilon\to0}\epsilon\log\mathbb{P}(\rho^\epsilon\in G)\geq -I(G).\]
    \end{enumerate}
\end{definition}
It is well known that if a sequence of random variables satisfies the large deviation principle with some rate function, then the rate function is unique, see for example Theorem 1.3.1 of \cite{dupuis2011weak}.\\
The result we will rely on is presented in \cite{budhiraja2008large} and is based on the variational representation of infinite dimensional Brownian.
To set up the problem and to state the relevant result from \cite{budhiraja2008large}, we first need to recall the spaces where various objects live.
First of all, based on the results in \cite{popat2025well}, in particular the assumptions on the boundary data stated in Assumption \ref{asm: assumptions on boundary data for well-posedness} in the Appendix, we require that the boundary data $\bar{f}$ lives in the space $H^{1}(\partial U)$. 
Throughout we denote functions defined on the boundary $\partial U$ with overbars.
\begin{definition}[The space $H^{1}(\partial U)$]\label{def: boundary sobolev space}
    The space of real valued functions defined on $\partial U$ with finite $H^1(\partial U)$-Sobolev norm is
    \[H^1(\partial U):=\{\bar{h}:\partial U\to\mathbb{R}: \|\bar{h}\|_{L^2(\partial U)}+\|\nabla_{\partial U}\bar{h}\|_{L^2(\partial U)} <\infty \},\]
    where $\nabla_{\partial U}\bar{h}$ denotes the tangential derivative of $\bar{h}$. That is, the gradient of $\bar{h}$ in directions tangent to $\partial U$.
    \end{definition}
On page 176 of \cite{fabes1978potential} this norm is rigorously defined on a $C^1$-regular domain in terms of a local coordinate system of $\partial U$.
Briefly, the definition involves considering a covering of $\partial U$ by a local co-ordinate chart and measuring the $L^2(\partial U)$-norms of the function and its derivative in these charts. 
The norm is well-defined and independent of the choice of chart, see \cite{seeley1959singular}.\\
Due to the definition of stochastic kinetic solutions, see Definition 2.8 of \cite{popat2025well}, a priori one considers that the non-negative initial data lives in the space $L^1(U)$.
However, formal arguments in Section \ref{sec: Correct space for initial data for the large deviations principle} in the Appendix illustrate that we need to take the smaller space, the space of functions with finite entropy (with respect to $\Phi$), $Ent_{\Phi}(U)$.
\begin{definition}[Entropy space]\label{def: entropy space}
Let $\Phi\in C([0,\infty))\cap C^1_{loc}((0,\infty))$ be non-negative and $v_0$ be the harmonic PDE on $U$ with boundary data $log(\bar{f})$ as defined in point 4 of Assumption \ref{asm: assumptions on boundary data for well-posedness}.
Define
\begin{align}\label{eq: defn of Psi}
    \Psi_{\Phi}(\xi):=\int_0^\xi (\log(\Phi(x))-v_0(x))\,dx.
\end{align}
The space of functions with finite entropy (with respect to $\Phi$) is the space
\begin{equation}\label{eq: entropy space}
    Ent_\Phi(U):=\left\{\rho_0\in L^1(U): \rho_0\geq 0\quad\text{and}\quad \int_U\Psi_\Phi(\rho_0(x))\,dx<\infty\right\}.
\end{equation}
\end{definition}
The integrability of the $v_0$ term in \eqref{eq: defn of Psi} is a consequence of point 4 of Assumption \ref{asm: assumptions on boundary data for well-posedness}, so the boundedness condition in \eqref{eq: entropy space} is just a condition on the integrability of the first term on the right hand side of \eqref{eq: defn of Psi}.\\
As in Definition \ref{definition truncated noise}, we view the infinite sequence of Brownian motions
\[B:=(B^k)_{k\in\mathbb{N}}\] 
used to define the noise $\xi^K$ in Definition \ref{definition truncated noise} as living in the space $C([0,T];(\mathbb{R}^d)^\infty)$ equipped with the metric topology of co-ordinate wise convergence.
In this section we again assume that the noise satisfies Assumption \ref{assumption on noise}.\\
In the notation of \cite{budhiraja2008large}, for $\epsilon\in(0,1)$ and fixed boundary data $\bar{f}\in H^1(\partial U)$, if it exists, we will denote by
\[\mathcal{G}^{\bar{f},\epsilon}: Ent_\Phi(U)\times C([0,T];(\mathbb{R}^d)^\infty)\to L^1(U\times[0,T])\]
to be the solution map for equation \eqref{generalised Dean-Kawasaki Equation Stratonovich with epsilon}. 
That is, if we let $\rho^{\epsilon}(\bar{f},\rho_0)$ denote the stochastic kinetic solution of \eqref{generalised Dean-Kawasaki Equation Stratonovich with epsilon} with boundary data $\Phi(\rho^\epsilon)|_{\partial U}=\bar{f}\in H^1(\partial U)$ and initial condition $\rho_0\in Ent_\Phi(U)$, then distributionally we have
\begin{equation}\label{eq: solution map for dean kawasaki}
    \rho^{\epsilon}(\bar{f},\rho_0)=\mathcal{G}^{\bar{f},\epsilon}(\rho_0,\sqrt{\epsilon}B).
\end{equation}
The existence of such a map $\mathcal{G}^{\bar{f},\epsilon}$ is shown as part of the proof of Proposition \ref{prop: existence of solution map for controlled SPDE} below.
We state an assumption under which we have the large deviations principle for \eqref{generalised Dean-Kawasaki Equation Stratonovich with epsilon}, see  Assumption 2 of \cite{budhiraja2008large}.
\begin{assumption}[Assumption for large deviations principle]\label{asm: Assumption for large deviations principle}
     Fix $\bar{f}\in H^1(\partial U)$. Suppose that there exists a measurable map $\mathcal{G}^{\bar{f},0}: Ent_\Phi(U)\times C([0,T];(\mathbb{R}^d)^\infty)\to L^1(U\times[0,T])$ such that
     \begin{enumerate}
         \item For every $R<\infty$ and compact subset $K\subset Ent_{\Phi}(U)$, the set
         \[\Gamma_{R,K}:=\left\{\mathcal{G}^{\bar{f},0}\left(\rho_0,\int_0^\cdot g(s)\,ds\right): g\in L^2_R(U\times[0,T]), \rho_0\in K\right\}\]
    is a compact subset of $L^1(U\times[0,T])$,
    where we defined the bounded $L^2(U\times[0,T])$ space by
    \[L^2_R(U\times[0,T]):=\{g\in L^2(U\times[0,T];\mathbb{R}^d):\|g\|_{L^2(U\times[0,T])}\leq R\}.\]
    \item For an arbitrary family of initial conditions $\{\rho_0^\epsilon\}_{\epsilon\in(0,1)}\subset Ent_\Phi(U)$ and controls $\{g^\epsilon\}_{\epsilon\in(0,1)}\subset L^\infty(\Omega; L^2(U\times[0,T];\mathbb{R}^d))$,  whenever we have the weak convergences\footnote{We ask for weak convergence of the initial data in the space $L^1(U)$ since it is more intuitive to understand what it means to converge weakly in $L^1(U)$. That is to say, the dual space of $L^1(U)$, the space of bounded linear functionals on $L^1(U)$, is better understood than the dual space of $Ent_\Phi(U)$.} $\rho_0^\epsilon\to \rho_0$ in $L^1(U)$ and $g^\epsilon \to g$ in $L^2(U\times[0,T])$ as $\epsilon\to 0$, then we have that distributionally in $L^1(U\times[0,T])$, as $\epsilon\to0$,
    \[\mathcal{G}^{\bar{f},\epsilon}\left(\rho_0^\epsilon,\sqrt{\epsilon}B+\left(\int_0^\cdot g_k^\epsilon(s)\,ds\right)_{k\in\mathbb{N}}\right)\to \mathcal{G}^{\bar{f},0}\left(\rho_0,\left(\int_0^\cdot g_k(s)\,ds\right)_{k\in\mathbb{N}}\right),\]
    where for $k\in\mathbb{N}$, $g^\epsilon_k(s):=\langle g^\epsilon(\cdot,s),f_k\rangle_{L^2(U)}$ and $g_k(s):=\langle g(\cdot,s),f_k\rangle_{L^2(U)}$ denote projections onto the spatial components of the noise.
    \end{enumerate}
\end{assumption}
\begin{theorem}[Sufficient condition for large deviations principle to hold]\label{thm: LDP theorem from BDM}
    For $\bar{f}\in H^1(\partial U)$, $\rho_0\in Ent_\Phi(U)$ and $\rho\in L^1(U\times[0,T])$, define 
    \begin{align}\label{eq: rate function}
        I_{\bar{f},\rho_0}(\rho)&:=\inf_{g\in L^2(U\times[0,T];\mathbb{R}^d): \rho=\mathcal{G}^{\bar{f},0}(\rho_0,\int_0^\cdot g(s)\,ds)}\left\{\frac{1}{2}\int_0^T\|g(\cdot, s)\|_{L^2(U;\mathbb{R}^d)}^2\,ds\right\}\nonumber\\
        &=\frac{1}{2}\inf_{g\in L^2(U\times[0,T];\mathbb{R}^d)}\left\{\|g\|^2_{L^2(U\times[0,T];\mathbb{R}^d)}:\partial_t\rho=\Delta\Phi(\rho)-\nabla\cdot(\sigma(\rho)g+ \nu(\rho))\right.\nonumber\\
        &\hspace{200pt}\left.: \Phi(\rho)|_{\partial U}=\bar{f}\,\, \text{and} \,\, \rho(\cdot,0)=\rho_0\right\}.
    \end{align}
    Suppose that Assumption \ref{asm: Assumption for large deviations principle} holds, and also that for all $\bar{f}\in H^1(\partial U)$ and $\rho\in L^1(U\times[0,T])$, the map $\rho_0\mapsto I_{\bar{f},\rho_0}(\rho)$ is lower semi-continuous from $Ent_\Phi(U)$ to $[0,\infty]$.\\
    Then, for arbitrary fixed fixed $\bar{f}\in H^1(\partial U)$ and every $\rho_0\in Ent_\Phi(U)$, the map $\rho\mapsto I_{\bar{f},\rho_0}(\rho)$ is a rate function on $L^1(U\times[0,T])$, the family $\{I_{\bar{f},\rho_0}(\cdot): \bar{f}\in H^1(\partial U),\, \rho_0\in Ent_\Phi(U)\}$ of rate functions has compact level sets on compacts, and $\{\rho^{\epsilon}(\bar{f},\rho_0)\}_{\epsilon\in(0,1)}$ satisfies the large deviations principle on $L^1(U\times[0,T])$ with rate function $I_{\bar{f},\rho_0}$.
 \end{theorem}
\begin{remark}[Uniform large deviations principle]\label{rmk: uniform ldp}
The conditions in Theorem \ref{thm: LDP theorem from BDM} do not just imply a large deviations principle, but in fact imply the stronger \say{uniform} large deviations principle.
That is, the large deviations principle holds uniformly with respect to compact subsets of the initial condition space $Ent_{\Phi,R}(U)$ defined by
\[Ent_{\Phi,R}(U):=\left\{\rho_0\in Ent_{\Phi}(U): \int_U\Psi_\Phi(\rho_0)\,dx<R\right\}.\]
See Definition 5 of \cite{budhiraja2008large} for further details.   
\end{remark}
    \subsection{Existence of solution maps}\label{sec: Weak solutions of the skeleton equation, existence of solution map for controlled SPDE}
Our goal is to unpack and prove the various conditions stated in Assumption \ref{asm: Assumption for large deviations principle} and Theorem \ref{thm: LDP theorem from BDM}.
In this section we begin by proving well-posedness of weak solutions of the limiting parabolic-hyperbolic PDE that appears in the rate function, which we will refer to as the skeleton equation.
An $L^1(U)$-contraction then implies the existence of a solution map $\mathcal{G}^{\bar{f},0}$ for the skeleton equation.
We end the section by proving the existence of solution maps $\{\mathcal{G}^{\bar{f},\epsilon},\epsilon\in(0,1)\}$ for the controlled SPDE for every $\epsilon\in(0,1)$.\\
The skeleton equation is the equation 
\begin{equation}\label{eq: skeleton equation}
   \begin{cases}
    \partial_t\rho=\Delta\Phi(\rho)-\nabla\cdot(\sigma(\rho)g+ \nu(\rho)), & \text{on} \hspace{5pt} U\times(0,T],\\
\Phi(\rho)=\bar{f},&\text{on}\hspace{5pt} \partial U\times[0,T],\\
    \rho(\cdot,t=0)=\rho_0 ,&\text{on}\hspace{5pt} U\times\{t=0\}.
   \end{cases} 
\end{equation}
We first define what it means to be a weak solution of the skeleton equation. 
\begin{definition}[Weak solution of the skeleton equation]\label{def: weak solution of skeleton equation}
    Let the non-linear functions $\Phi,\sigma,\nu$ satisfy Assumptions \ref{asm: assumptions for well posedenss of equation} and \ref{asm: assumptions on boundary data for well-posedness}, and let $g\in L^2(U\times[0,T];\mathbb{R}^d)$ and $\rho_0\in Ent_\Phi(U)$. 
    For boundary data $\bar{f}\in H^1(\partial U)$, define the  PDE $h$ as the unique harmonic function on $U$ with boundary data $\Phi^{-1}(\bar{f})$ as in equation \eqref{eq: pde h} of the Appendix.
    A weak solution of \eqref{eq: skeleton equation} with initial data $\rho_0$ and boundary data $\bar{f}$ is a continuous $L^1(U)$-valued function $\rho\in L^\infty([0,T];L^1(U))$ satisfying
    \begin{enumerate}
        \item Boundary condition and regularity: $(\Phi^{1/2}(\rho)-\Phi^{1/2}(h))\in L^2([0,T];H^1_0(U))$.
        \item The equation: For every $t\in[0,T]$ and $\psi\in C_c^\infty(U)$ we have
        \begin{multline*}
            \int_U\rho(x,t)\psi(x)\,dx=\int_U\rho_0(x)\psi(x)\,dx-\int_0^t\int_U\Phi'(\rho)\nabla\rho\cdot\nabla\psi\,dx\,dt\\
            +\int_0^t\int_U\sigma(\rho) g\cdot\nabla\psi\,dx\,dt+\int_0^t\int_U\nu(\rho)\cdot\nabla\psi\,dx\,dt.
        \end{multline*}
    \end{enumerate}
\end{definition}
\begin{theorem}[Well-posedness of the skeleton equation]\label{thm: well-posedness of skeleton equation assumed}
    For every $g\in L^2(U\times[0,T];\mathbb{R}^d)$, $\rho_0\in Ent_\Phi(U)$ and $\bar{f}\in H^1(\partial U)$, under some assumptions\footnote{See Assumptions 6 and 10 of \cite{fehrman2023non}.} on the non-linear coefficients $\Phi,\sigma$ and $\nu$, there exists a unique weak solution to the skeleton equation \eqref{eq: skeleton equation}.
    Furthermore, we have the below equivalence of stochastic kinetic solutions and weak solutions,
    \begin{align*}
        &\rho \, \text{is a stochastic kinetic solution of \eqref{eq: skeleton equation} with initial data $\rho_0$ and boundary data $\bar{f}$}\\
        &\iff \rho \,\text{is a weak solution of \eqref{eq: skeleton equation} in the sense of Definition \ref{def: weak solution of skeleton equation}.}
    \end{align*}
    Finally, if $\rho^1,\rho^2$ are two weak solutions to \eqref{eq: skeleton equation} with initial data $\rho^1_0,\rho^2_0\in Ent_\Phi(U)$ and the same control $g\in L^2(U\times[0,T];\mathbb{R}^d)$, then we have the contraction
    \begin{equation}\label{eq: l1 contraction for skeleton equation}
        \sup_{t\in[0,T]}\|\rho^1(\cdot,t)-\rho_2(\cdot,t)\|_{L^1(U)}\leq \|\rho^1_0-\rho^2_0\|_{L^1(U)}.
    \end{equation}
\end{theorem}
\begin{proof}
    The proof on the torus for the choice of coefficients $\nu=0, \sigma=\Phi^{1/2}$ is precisely the content of Theorems 8, 14 and Proposition 20 of \cite{fehrman2023non}. 
    Theorem 8 proves the existence of stochastic kinetic solutions and the $L^1(U)$-contraction \eqref{eq: l1 contraction for skeleton equation}, and is similar to Theorem 3.6 of \cite{popat2025well}. 
    Proposition 20 subsequently proves existence of weak solutions. 
    These two are tied together by the equivalence of weak and stochastic kinetic solutions, which is the content of Theorem 14 of \cite{fehrman2023non}.\\
    For the bounded domain with general coefficients considered in the present work, the result follows by adapting the arguments contained in \cite{fehrman2023non} with arguments from from \cite{popat2025well}.
    \end{proof}
\begin{remark}[Point 1 of Assumption \ref{asm: Assumption for large deviations principle}]\label{rmk: verifying first point in assumption of ldp}
    After we have the well-posedness of the skeleton equation, we turn to Section \ref{sec: Correct space for initial data for the large deviations principle} in the Appendix for an entropy estimate for the equation.
This estimate alongside the tightness result of Proposition 17 in \cite{fehrman2023non} proves that the solution set $\Gamma_{R,K}$ of the skeleton equation is compact in $L^1(U\times[0,T])$.
Hence the first point of Assumption \ref{asm: Assumption for large deviations principle} is satisfied.
\end{remark}
We conclude the discussion on the skeleton equation with a result on the existence of a solution map, which verifies another condition in Assumption \ref{asm: Assumption for large deviations principle}.
\begin{proposition}[Existence of solution map for skeleton equation]\label{prop: existence of solution map for skeleton equation}
    Suppose Assumptions \ref{asm: assumptions for well posedenss of equation} and \ref{asm: assumptions on boundary data for well-posedness} are satisfied, and fix the boundary data $\bar{f}\in H^1(\partial U)$.
    For initial condition $\rho_0\in Ent_\Phi(U)$, $\mathcal{F}_t$-predictable control $g\in L^2(U\times[0,T];\mathbb{R}^d)$ and $g_k(s):=\langle g(\cdot,s),f_k\rangle_{L^2(U)}$, there exists a measurable map
    \[\mathcal{G}^{\bar{f},0}:Ent_\Phi(U)\times C([0,T];(\mathbb{R}^d)^\infty)\to L^1(U\times[0,T])\]
    such that
    \[\rho(\bar{f},\rho_0)=\mathcal{G}^{\bar{f},0}\left(\rho_0,\left(\int_0^\cdot g_k(s)\,ds\right)_{k\in\mathbb{N}}\right),\]
    where $\rho(\bar{f},\rho_0)$ denotes the unique weak solution of \eqref{eq: skeleton equation} with initial data $\rho_0$ and boundary data $\Phi(\rho)|_{\partial U}=\bar{f}$ in the sense of Definition \ref{def: weak solution of skeleton equation} above.
\end{proposition}
\begin{proof}
The proof is a simplified version of the proof for Proposition \ref{prop: existence of solution map for controlled SPDE} below, relying on two main ingredients.
The first is the well-posedness of the equation and the second is the $L^1(U)-$contraction, both conditions were verified above in Theorem \ref{thm: well-posedness of skeleton equation assumed}.
\end{proof}

Let us turn our attention to proving the existence of a solution map 
\begin{equation}\label{eq: for G epsilon}
    \mathcal{G}^{\bar{f},\epsilon}\left(\rho_0^\epsilon,\sqrt{\epsilon}B+\left(\int_0^\cdot g_k^\epsilon(s)\,ds\right)_{k\in\mathbb{N}}\right)
\end{equation}
appearing in point 2 of Assumption \ref{asm: Assumption for large deviations principle}. 
That is, for fixed $\epsilon\in(0,1)$ and boundary data $\bar{f}\in H^1(U)$, we want to prove the existence of a solution map for the controlled SPDE
\begin{equation}\label{eq: controlled SPDE}
    \partial_t\rho^{\epsilon,g^\epsilon}=\Delta\Phi(\rho^{\epsilon,g^\epsilon})-\sqrt{\epsilon}\nabla\cdot(\sigma(\rho^{\epsilon,g^\epsilon})\circ\dot{\xi}^{K(\epsilon)})-\nabla\cdot(\sigma(\rho^{\epsilon,g^\epsilon})P_{K(\epsilon)} g^\epsilon)-\nabla\cdot \nu(\rho^{\epsilon,g^\epsilon}),
\end{equation} 
with initial data $\rho_0^\epsilon\in Ent_\Phi(U)$, boundary data $\Phi(\rho^{\epsilon,g^\epsilon})|_{\partial U}=\bar{f}$, and recall that $P_{K(\epsilon)} g^\epsilon$ denotes the projection of the control $g^\epsilon\in L^2(U\times[0,T];\mathbb{R}^d)$ onto the span of $\{f_k\}_{k=1}^{K(\epsilon)}$.
\begin{remark}[Notation in equation \eqref{eq: controlled SPDE}]
    Initially it may appear as if we are abusing notation and using $\epsilon$ to denote both the scaling in front of the noise term in \eqref{eq: controlled SPDE} and also to denote the approximating sequences $\{g^\epsilon\}_{\epsilon\in(0,1)},\{\rho_0^\epsilon\}_{\epsilon\in(0,1)}$.
    Put differently, $\epsilon$ appears in the solution map \eqref{eq: for G epsilon} both when writing the approximating solution $\mathcal{G}^\epsilon$ which corresponds to the scaling in front of the noise term and the arguments inside the function corresponding to the approximating sequence of coefficients. 
    This is the convention in \cite{budhiraja2008large} and there is no issue with choosing the same scaling in $\epsilon$ for both objects.
\end{remark}

We will need the following result, see Proposition 27 of \cite{fehrman2023non} or Proposition 7.12 of \cite{fehrman2019large} for proof.
\begin{proposition}[Separability of $Ent_\Phi(U)$]\label{prop: entropy space is separable}
    Let $\Phi\in C([0,\infty))\cap C^1_{loc}((0,\infty))$ be non-negative and satisfy Assumption \ref{asm: assumptions for well posedenss of equation}. 
    The space  $Ent_{\Phi}(U)$ equipped with the $L^1(U)$-topology is a complete, separable metric space.
\end{proposition}
The existence of a solution map for the controlled SPDE \eqref{eq: controlled SPDE} below is similar to Theorem 23 and Proposition 24 of \cite{fehrman2023non}.
In the below proof we fix the boundary data $\bar{f}\in H^1(\partial U)$ and  prove that the solution map exists for fixed (but arbitrary) boundary data.
\begin{proposition}[Existence of solution map for controlled SPDE]\label{prop: existence of solution map for controlled SPDE}
    Suppose Assumptions \ref{asm: assumptions for well posedenss of equation} and \ref{asm: assumptions on boundary data for well-posedness} are satisfied and fix the boundary data $\bar{f}\in H^1(\partial U)$.
    Then for every $\epsilon\in(0,1)$, initial condition $\rho_0\in Ent_\Phi(U)$, $\mathcal{F}_t$-predictable control $g\in L^2(U\times[0,T];\mathbb{R}^d)$ and $g^\epsilon_k(s):=\langle g^\epsilon(\cdot,s),f_k\rangle_{L^2(U)}$, there exists a measurable map
    \[\mathcal{G}^{\bar{f},\epsilon}:Ent_\Phi(U)\times C([0,T];(\mathbb{R}^d)^\infty)\to L^1(U\times[0,T])\]
    such that
    \[\rho^{\epsilon,g}(\bar{f},\rho_0)=\mathcal{G}^{\bar{f},\epsilon}\left(\rho_0,\sqrt{\epsilon}B+\left(\int_0^\cdot g_k^\epsilon(s)\,ds\right)_{k\in\mathbb{N}}\right),\]
    where $\rho^{\epsilon,g}(\bar{f},\rho_0)$ is the unique stochastic solution of \eqref{eq: controlled SPDE} with initial data $\rho_0$ and boundary data $\Phi(\rho^{\epsilon,g})|_{\partial U}=\bar{f}$ in the sense of Definition 2.8 of \cite{popat2025well}.
\end{proposition}
\begin{proof}
    The proof consists of two steps.
    Firstly we will we will show the existence of a map $\mathcal{G}^{\bar{f},\epsilon}$ such that for the uncontrolled system \eqref{generalised Dean-Kawasaki Equation Stratonovich with epsilon}, we have
    \[\rho^{\epsilon}(\bar{f},\rho_0)=\mathcal{G}^{\bar{f},\epsilon}\left(\rho_0,\sqrt{\epsilon}B\right).\]
    The second step will be to extend the result to the controlled SPDE $\rho^{\epsilon,g}$ \eqref{eq: controlled SPDE} via Girsanov theorem.\\
    The first step follows the same lines as Theorem 23 of \cite{fehrman2023non}, see also Proposition 4.5 of \cite{dirr2020conservative} for a similar argument.
    We will therefore be brief.\\
    For every fixed $\epsilon\in(0,1)\in\mathbb{N}$, the well-posedness of stochastic kinetic solutions of equation \eqref{generalised Dean-Kawasaki Equation Stratonovich with epsilon} is a consequence of \cite{popat2025well}.
    Now, due to the well-posedness, we know that for fixed boundary data $\bar{f}\in H^1(\partial U)$ and initial condition $\rho_0\in Ent_\Phi(U)$, there exists a measurable function $\mathcal{G}^{\epsilon}_{\bar{f},\rho_0}:C([0,T];(\mathbb{R}^d)^\infty)\to L^1(U\times[0,T])$ that takes the noise to the solution of \eqref{generalised Dean-Kawasaki Equation Stratonovich with epsilon}, that is
    \[\mathcal{G}^{\epsilon}_{\bar{f},\rho_0}(\sqrt{\epsilon}B)=\rho^{\epsilon}(\bar{f},\rho_0).\]
    Since it follows from Proposition \ref{prop: entropy space is separable} above that $Ent_\Phi(U)$ is separable with the $L^1(U)$-topology,
    by looking at a countable dense subset $\{\rho_n\}_{n\in\mathbb{N}}$ of $Ent_\Phi(U)$,
    Theorem 3.6 of \cite{popat2025well} implies that for every $n,m\in\mathbb{N}$,
      \begin{align}\label{eq: l1 contraction for solution map}
          \sup_{t\in[0,T]}\|\mathcal{G}^{\epsilon}_{\bar{f},\rho_n}\left(\sqrt{\epsilon}B\right)-\mathcal{G}^{\epsilon}_{\bar{f},\rho_m}\left(\sqrt{\epsilon}B\right)\|_{L^1(U)}\leq \|\rho_n-\rho_m\|_{L^1(U)}.
      \end{align}
      As a result of the density of $\{\rho_n\}_{n\in\mathbb{N}}$ with respect to the $L^1(U)$-norm,
      it follows that there exists a measurable function $\mathcal{G}^{\epsilon}_{\bar{f},\sqrt{\epsilon}B}: Ent_\Phi(U)\to L^1(U\times[0,T])$ that maps for every $n\in\mathbb{N}$ 
      the initial condition $\rho_n$ to the solution of \eqref{generalised Dean-Kawasaki Equation Stratonovich with epsilon},
    \[\mathcal{G}^{\epsilon}_{\bar{f},\sqrt{\epsilon}B}(\rho_n)=\rho^{\epsilon}(\bar{f},\rho_n)=\mathcal{G}^{\epsilon}_{\bar{f},\rho_n}(\sqrt{\epsilon}B).\]
   The fact that $\{\rho_n\}_{n\in\mathbb{N}}$ 
   is dense then implies that for arbitrary $\rho_0\in Ent_\Phi(U)$,
   on a subset of full probability depending on $\rho_0$, 
    \[ \mathcal{G}^{\epsilon}_{\bar{f},\sqrt{\epsilon}B}(\rho_0)=\rho^{\epsilon}(\bar{f},\rho_0)=\mathcal{G}^{\epsilon}_{\bar{f},\rho_0}(\sqrt{\epsilon}B).\]
    Finally, the desired solution map $\mathcal{G}^{\bar{f},\epsilon}$ is defined by
    \[\mathcal{G}^{\bar{f},\epsilon}\left(\rho_0,\sqrt{\epsilon}B\right):=\mathcal{G}^{\epsilon}_{\bar{f},\sqrt{\epsilon}B}\left(\rho_0\right).\]
    The measurability of $\mathcal{G}^{\bar{f},\epsilon}$ follows from the measurability of the map
    $\mathcal{G}^{\bar{f},\epsilon}\left(\rho_0,\cdot\right)$
    viewed as a function of the noise, and the strong continuity of the map
     $\mathcal{G}^{\bar{f},\epsilon}\left(\cdot,\sqrt{\epsilon}B\right)$ 
     viewed as a function of the initial datum.\\
     Let us now move onto the second step, which we note follows from Theorem 10 of \cite{budhiraja2008large}.
     Fix $\epsilon\in(0,1)$ and measurable control $g\in L^2(U\times[0,T];\mathbb{R}^d)$. 
     The measure
     \[d\gamma^{\epsilon,g}:=\exp\left\{-\frac{1}{\sqrt{\epsilon}}\sum_{k=1}^K\int_0^T\int_U g(x,s)\,dx\,dB_s^k-\frac{1}{2\epsilon}\sum_{k=1}^K\int_0^T\int_U g^2(x,s)\,dx\,ds\,\right\}\,d\mathbb{P}\]
     is a probability measure on $(\Omega,\mathcal{F},\mathbb{P})$, where the probability space is as in Definition \ref{definition truncated noise}.
    The measure $\gamma^{\epsilon,g}$ is absolutely continuous with respect to $\mathbb{P}$, and by Girsanov's theorem (see Theorem 10.14 of \cite{da2014stochastic}) the process
    \[\{\tilde{B}^k_\cdot\}_{k\in\mathbb{N}}:=\left\{B^k_\cdot+\epsilon^{-1/2}\int_0^\cdot g_k \right\}_{k\in\mathbb{N}}\]
     on the probability space $(\Omega,\mathcal{F},\{\mathcal{F}_t\}_{t\geq0},\gamma^{\epsilon,g})$ is an infinite sequence of independent $d-$dimensional Brownian motions.
    Hence, by the first step and uniqueness of solutions of \eqref{generalised Dean-Kawasaki Equation Stratonovich with epsilon} as shown in \cite{popat2025well}, we have that the quantity of interest
    \[\mathcal{G}^{\bar{f},\epsilon}\left(\rho_0,\sqrt{\epsilon}B+\left(\int_0^\cdot g_k^\epsilon(s)\,ds\right)_{k\in\mathbb{N}}\right)\]
    is the unique stochastic kinetic solution to \eqref{generalised Dean-Kawasaki Equation Stratonovich with epsilon} with $B$ replaced by $\tilde{B}$ on $(\Omega,\mathcal{F},\{\mathcal{F}_t\}_{t\geq0},\gamma^{\epsilon,g})$.
    But \eqref{generalised Dean-Kawasaki Equation Stratonovich with epsilon} with $\tilde{B}$ is precisely the same as $\rho^{\epsilon,g}$ as in \eqref{eq: controlled SPDE}.
    The result then follows on the original probability space $(\Omega,\mathcal{F},\{\mathcal{F}_t\}_{t\geq0},\mathbb{P})$ using the fact that the measures $\mathbb{P}$ and $\gamma^{\epsilon,g}$ are absolutely continuous.
\end{proof}

\subsection{Lack of uniformity of solution map with respect to initial condition} \label{sec: lack of uniformity of solution map with respect to initial condition}
It is natural to wonder whether it is possible to extend the uniformity of the large deviations principle to include uniformly bounded subsets of boundary data as well.
Rather than looking at the solution map $\mathcal{G}^{\bar{f},\epsilon}$ as in \eqref{eq: solution map for dean kawasaki}, we would need the existence of a solution map
\[\mathcal{G}^{\epsilon}: H^1(\partial U)\times Ent_\Phi(U)\times C([0,T];(\mathbb{R}^d)^\infty)\to L^1(U\times[0,T])\]
that also takes the boundary data as input.
Recall that a key step in the proof of Proposition \ref{prop: existence of solution map for controlled SPDE} was the $L^1(U)$-contraction
\[\sup_{t\in[0,T]}\|\rho^{\epsilon,1}(\cdot,t)-\rho^{\epsilon,2}(\cdot,t)\|_{L^1(U)}\leq \|\rho^{1}_0-\rho^{2}_0\|_{L^1(U)}\]
for two stochastic kinetic solutions $\rho^{\epsilon,1},\rho^{\epsilon,2}$ of \eqref{generalised Dean-Kawasaki Equation Stratonovich with epsilon} with the same boundary data $\bar{f}$ and different initial conditions $\rho_0^1,\rho_0^2$.
To have uniformity with respect to the boundary data we would need to prove that if the two stochastic kinetic solutions had different initial data, say $\bar{f}^1,\bar{f}^2\in H^1(\partial U)$ respectively, then for some constant $c\in(0,\infty)$ we have an inequality of the form
   \begin{equation}\label{eq: l1 contraction with different boundary data}
       \sup_{t\in[0,T]}\|\rho^{\epsilon,1}(\cdot,t)-\rho^{\epsilon,2}(\cdot,t)\|_{L^1(U)}\leq \|\rho^{1}_0-\rho^{2}_0\|_{L^1(U)}+cT\|\bar{f}^1-\bar{f}^2\|_{H^1(\partial U)}.
   \end{equation}    
In this section we briefly discuss how the uniqueness proof, Theorem 3.6 of \cite{popat2025well}, would be amended and subsequently conclude that it is not possible to obtain an inequality of the form \eqref{eq: l1 contraction with different boundary data}.
\\
    A key observation in the uniqueness proof is that the $L^1(U)-$difference of solutions can be written in terms of the corresponding kinetic functions $\chi^1,\chi^2$ of the two solutions, see Definition \ref{def: kinetic function} in the Appendix below. 
    That is, for every $t\in[0,T]$ we have the identity
    \begin{multline*}
        \left.\int_U|\rho^{\epsilon,1}(x,s)-\rho^{\epsilon,2}(x,s)|\,dx\right|_{s=0}^t\\
        =\left.\int_{\mathbb{R}}\int_U\chi^1(\xi,\rho^{\epsilon,1}(x,s))+\chi^2(\xi,\rho^{\epsilon,1}(x,s))-2\chi^1(\xi,\rho^{\epsilon,1}(x,s))\chi^2(\xi,\rho^{\epsilon,1}(x,s))\right|_{s=0}^t.
    \end{multline*}
    The terms on the right hand side are handled using the kinetic equation, see for example equation \eqref{eq: kinetic equation} in the Appendix below.
    Since test functions in the kinetic equation need to be compactly supported in the spatial and velocity variables, for the first two terms above we need to test against (regularised versions) of relevant cutoff functions. 
    For the final term we need to include the kinetic function as part of the test function, which requires us to smooth the kinetic function via convolution.\\
    The presence of the cutoff function in space $\iota_\gamma$ with parameter $\gamma\in(0,1)$, see Definition 3.2 of \cite{popat2025well}, enables us to integrate by parts and perform analysis without the presence of additional boundary dependent terms. 
    It is precisely when we take the limit in the cutoff parameter $\gamma\to 0$ that boundary terms arise. 
    Specifically, for the cutoff term $I_t^{cut}$ in the proof of Theorem 3.6 of \cite{popat2025well}, we notice from equation (15) there that before taking limits in the cutoff functions, we have the term
    \begin{equation}\label{eq: equation 15 from popat}
        \gamma^{-1}\left(\int_0^t\int_{\partial U}|\Phi_{M,\beta}(\rho^{\epsilon,2})-\Phi_{M,\beta}(\rho^{\epsilon,1})|-\int_0^t\int_{\partial U_\gamma}|\Phi_{M,\beta}(\rho^{\epsilon,2})-\Phi_{M,\beta}(\rho^{\epsilon,1})|\right),
    \end{equation}
    where $\Phi_{M,\beta}$ is the unique function given by $\Phi_{M,\beta}(0)=0$, $\Phi_{M,\beta}'(\xi)=\zeta_M(\xi)\phi_\beta(\xi)\Phi'(\xi)$, and $\zeta_M$ and $\phi_\beta$ are smooth cutoff functions for large and small values respectively, converging pointwise to the constant function $1$ as $M\to\infty$ and $\beta\to0$, see Definition 3.2 of \cite{popat2025well}.
    For $\gamma\geq 0$ sufficiently small, in \eqref{eq: equation 15 from popat} we defined the level set $\partial U_\gamma:=\{x\in U: d(x,\partial U)=\gamma\}$, where $d(x,\partial U)$ denotes the standard minimum Euclidean distance from a point to a set.\\
The below technical lemma illustrates why once can not expect a bound such as \eqref{eq: l1 contraction with different boundary data} to hold for equation \eqref{eq: equation 15 from popat}. 
The main point is that the $H^1(\partial U)$-norm as in Definition \ref{def: boundary sobolev space} looks at tangential derivatives of $f$ along the boundary, whereas the lemma below illustrates that a bound can only be obtained in terms of normal derivatives at the boundary.\\
To easily differentiate between points in the interior of $U$ and points on the boundary $\partial U$, for the below lemma we denote points on the boundary with an asterisk, for example $x^*$. 
\begin{lemma}\label{lemma: bounding difference of surface integrals}
Let $f:U\to\mathbb{R}$ satisfy $f\in H^s(U)$, for some $s>3/2$, so that by an extension of the trace theorem, see Theorem 1.3.1 of \cite{quarteroni2008numerical}, we have that its trace satisfies $f|_{\partial U}\in H^1(\partial U)$.\\
    Then we have the estimate
    \begin{align}\label{eq: equality for differece in boundary integrals}
        \lim_{\gamma\to0}\left(\gamma^{-1}\left(\int_{\partial U}f(x^*)-\int_{\partial U_\gamma}f(x)\right)\right)=\int_{\partial U}H(x^*)f(x^*)+\int_{\partial U}\nabla f(x^*)\cdot\hat\eta(x^*),
    \end{align}
    where $H(x^*)$ denotes the mean curvature of the boundary $\partial U$ at the point $x^*$ and $\hat{\eta}(x^*)$ denotes the outward pointing normal vector field at point $x^*$.
  \end{lemma}
\begin{proof}
   We first observe that for $\gamma$ sufficiently small, one can write points on $\partial U_\gamma$ as a smooth perturbation of points on $\partial U$ by moving along the inward pointing normal vector.
    For every $x\in\partial U_\gamma$, let $x^*=x^*(x):=\Pi_{\partial U}(x)$ denote the closest point on the boundary to $x$. 
    Since $U$ is a $C^2$-domain, we know that for $\gamma$ sufficiently small, such a point exists and is unique, see page 153 of \cite{foote1984regularity}.
    If $\hat{\eta}(x^*)$ denotes the outward pointing unit normal vector field at point $x^*\in \partial U$, then
    \[x=x^*-\gamma\hat\eta(x^*)+O(\gamma^2),\]
    where the final term denotes dependence on terms depending on $\gamma$ of order $\gamma^2$ or higher.
    Similarly, we have for the surface measure $d S_\gamma$ on $\partial U_\gamma$ the first order expansion\footnote{the formula is derived form the first order expansion of the Jacobian determinant associated with the normal flow at the boundary.}
    \[dS_\gamma(x)=\left(1-\gamma H(x^*)+O(\gamma^2)\right)dS(x^*),\]
    where $dS$ denotes the surface measure on $\partial U$, and $H(x^*)$ denotes the mean curvature of the boundary $\partial U$ at point $x^*$.
    Putting these two facts together gives
    \begin{align}\label{eq: converting surface integrals}
        \int_{\partial U_\gamma}f(x)dS_\gamma(x)=\int_{\partial U}f\left(x^*-\gamma\hat\eta(x^*)+O(\gamma^2)\right)\left(1-\gamma H(x^*)+O(\gamma^2)\right)dS(x^*).
    \end{align}
    The first order Taylor expansion
    \[f\left(x^*-\gamma\hat\eta(x^*)+O(\gamma^2)\right)=f(x^*)-\gamma\nabla f(x^*)\cdot\hat{\eta}(x^*)+O(\gamma^2)\]
    implies that we can write \eqref{eq: converting surface integrals} as
    \begin{align*}
         \int_{\partial U_\gamma}f(x)dS_\gamma(x)&=\int_{\partial U}\left(f(x^*)-\gamma\nabla f(x^*)\cdot\hat{\eta}(x^*)+O(\gamma^2)\right)\left(1-\gamma H(x^*)+O(\gamma^2)\right)dS(x^*)\\
         &=\int_{\partial U}f(x^*)dS(x^*)-\gamma\int_{\partial U}\nabla f(x^*)\cdot\hat{\eta}(x^*)-\gamma\int_{\partial U}f(x^*)H(x^*)dS(x^*)+O(\gamma^2),
    \end{align*}
    where we observe in the final line that several of the terms in the product can be absorbed into $O(\gamma^2)$. 
    Finally, the desired difference on the left hand side of \eqref{eq: equality for differece in boundary integrals} is
    \begin{align*}
       \gamma^{-1}\left(\int_{\partial U}f(x^*)-\int_{\partial U_\gamma}f(x)\right)=\int_{\partial U}\nabla f(x^*)\cdot\hat{\eta}(x^*)+\int_{\partial U}f(x^*)H(x^*)dS(x^*)+O(\gamma),
    \end{align*}
    which proves the statement after taking the limit $\gamma\to0$ on both sides.
\end{proof}
The essential difficulty is that stochastic kinetic solutions $\rho^{\epsilon,1},\rho^{\epsilon,1}$ of the Dirichlet problem \eqref{generalised Dean-Kawasaki Equation Stratonovich with epsilon} do not have sufficient regularity to evaluate the normal derivative
\[\int_{\partial U}\frac{\partial}{\partial \hat{\eta}}|\Phi_{M,\beta}(\rho^{\epsilon,2})-\Phi_{M,\beta}(\rho^{\epsilon,1})|,\]
let alone control it by the $H^1(\partial U)$-norm of the difference of the corresponding boundary datum.
We recall that in general, normal derivatives can only be controlled by the tangential derivative whenever a \say{Dirichlet-to-Neumann map} exists.
Such a map only is only known for simple evolution equations such as harmonic PDEs, see Theorem 2.4(iii) of \cite{fabes1978potential}, so that even if the normal derivative could be evaluated, there is no hope to obtain an estimate of the form \eqref{eq: l1 contraction with different boundary data}.
\subsection{Energy estimates for the regularised controlled SPDE}\label{sec: Energy estimates for the regularised controlled SPDE}
In order to prove the weak convergence of solutions required in point 2 of Assumption \ref{asm: Assumption for large deviations principle}, we have to first establish tightness of the laws of the controlled SPDE $\rho^{\epsilon,g}$ \eqref{eq: controlled SPDE}. 
In order to prove relevant energy estimates that allow us to prove tightness, for $\alpha\in(0,1)$ and $n\in\mathbb{N}$ we need to introduce the regularised controlled SPDE,
\begin{equation}\label{eq: regularised controlled SPDE}
\partial_t\rho^{n,\epsilon,g}=\Delta\Phi(\rho^{n,\epsilon,g})+\alpha\Delta\rho^{n,\epsilon,g}-\sqrt{\epsilon}\nabla\cdot(\sigma_n(\rho^{n,\epsilon,g})\circ\dot{\xi}^K)-\nabla\cdot(\sigma_n(\rho^{n,\epsilon,g})P_{K(\epsilon)} g)-\nabla\cdot \nu(\rho^{n,\epsilon,g}),
\end{equation} 
with initial data $\rho^{n,\epsilon,g}(\cdot,0)=\rho_0$ and boundary data $\Phi(\rho^{n,\epsilon,g})|_{\partial U}=\bar{f}$.
As we have done throughout, we do not make explicit the dependence on $\alpha\in(0,1)$ in the notation of the equation.
We recall that the regularisation of the square root enables us to consider weak solutions, and the regularisation in $\alpha$ provides sufficient regularity of solutions enabling the use of It\^o's formula.
The definition of weak solutions for \eqref{eq: regularised controlled SPDE} is similar to the definition for the uncontrolled equation, see Definition \ref{weak solution of regularised equation in CLT and LDP} above.
With the regularisations as in \eqref{eq: regularised controlled SPDE} the existence of weak solutions follows from standard arguments, by smoothing the coefficients, applying a Galerkin projection and taking limits of the projection and coefficients in the correct order. 
See proof of Proposition 4.30 of \cite{popat2025well} for a similar argument.\\
We state four relevant energy estimates for weak solutions of the regularised equation, which are all uniform in the regularisation coefficients $n,\alpha$ and $\epsilon$, and consequently tightness will be shown with respect to subsequences of all three parameters. 
Firstly we will prove an $L^2(U\times[0,T])-$estimate in Proposition \ref{prop: l2 estimate from 4.14 of popat}, which is similar to Proposition 4.14 of \cite{popat2025well}.
This estimate is used in Propositions \ref{prop: higher order spatial regularity for regularised controlled SPDE} and \ref{prop: higher order time regularity for regularised controlled SPDE} below where we prove higher order spatial regularity and higher order time regularity respectively, which are based on Lemma 4.16, Corollary 4.17 and Proposition 4.19 of \cite{popat2025well}.
These higher order regularity estimates allow us to prove tightness of the laws of \eqref{eq: regularised controlled SPDE}, which is done in Corollary \ref{corr: tightness of laws of controlled SPDE}.
Finally we will prove an entropy estimate, which is similar to Proposition 4.24 of \cite{popat2025well} and will be necessary in the sequel.\\
In \cite{popat2025well}, the noise is correlated in space, ensuring that for the noise coefficients $\{F^K_i\}_{i=1,2,3}$ as in Definition \ref{definition truncated noise}, there exists a constant $c\in(0,\infty)$ such that
\[\lim_{K\to\infty} F^K_i<c.\]
This is not the case for the noise considered here, see Remark \ref{rmk: limit of noise coefficients is unbounded}, but the estimates still hold under the scaling regime of Remark \ref{rmk: choice of joint scaling} due to the fact that in equation \eqref{eq: ito form of dk equation}, the noise coefficients are multiplied by a factor of $\epsilon$.
Owing to Remark \ref{rmk: choice of joint scaling}, picking $\epsilon\in(0,1)$ sufficiently small so that
    \begin{equation}\label{eq: scaling so that noise times epsilon is less than 1}
        \sqrt{\epsilon}\left(\|F_1^{K(\epsilon)}\|_{L^\infty(U)}+\|F_2^{K(\epsilon)}\|_{L^\infty(U)}+\|F_3^{K(\epsilon)}\|_{L^\infty(U)}+\|\nabla\cdot F_2^{K(\epsilon)}\|_{L^\infty(U)}\right)<1
    \end{equation} 
ensures that the bounds from \cite{popat2025well} can be directly used.
For simplicity of presentation, we make the below assumption for the below estimates.

\begin{assumption}[Constant boundary data]\label{asm: constant boundary data}
Assume that the non-negative boundary data of the regularised controlled SPDE $\rho^{n,\epsilon,g}$ as in equation \eqref{eq: regularised controlled SPDE} is a constant, $\Phi(\rho^{n,\epsilon,K})|_{\partial U}=M\geq0$.  
\end{assumption}
We emphasise that the above assumption is not necessary, but ensures that the bounds in this section are more digestible.
We can consider all boundary data considered in \cite{popat2025well}, namely non-negative constant data
including zero and all smooth functions bounded away from zero.
The estimates would involve boundary dependence and dependence on norms of certain harmonic PDEs, the precise terms can be found by looking at the corresponding results in \cite{popat2025well}.
All of the additional boundary terms that arise are controlled due to the various points in Assumption \ref{asm: assumptions on boundary data for well-posedness} in the Appendix.\\
The first estimate is a generalisation of the $L^p(U\times[0,T])-$estimates of Propositions \ref{ppn: energy estimate for difference of regularised equation and hydrodynamic limit} and \ref{prop: second energy estimate for rho n epsilon k} to the case of general boundary data, which can only be proven for the choice $p=2$.
\begin{proposition}[{$L^2(U\times[0,T])-$}estimate for the regularised equation]\label{prop: l2 estimate from 4.14 of popat}
      Let the coefficients of \eqref{eq: regularised controlled SPDE} satisfy Assumptions \ref{assumption on noise} and \ref{asm: assumptions for well posedenss of equation}.
    Let $g\in L^2(U\times[0,T];\mathbb{R}^d)$ be an $\mathcal{F}_t$-predictable control, $\rho_0\in L^2(\Omega;L^2(U))$ be non-negative and $\mathcal{F}_0$-measurable initial data, and to simplify the below estimate we assume that the boundary data $\bar{f}\in H^1(\partial U)$ satisfies Assumption \ref{asm: constant boundary data}.
     Suppose that $\epsilon\in(0,1)$ is sufficiently small so that the inequality \eqref{eq: scaling so that noise times epsilon is less than 1} is satisfied.
    For fixed regularisation constants $\alpha\in(0,1), n\in\mathbb{N}$, let $\rho^{n,\epsilon,g}$ denote weak solutions of the regularised controlled SPDE \eqref{eq: regularised controlled SPDE}.
    Then there exists a constant $c\in(0,\infty)$ independent of $\alpha,n,\epsilon$ such that
    \begin{multline}\label{eq: l2 space-time estimate for regularised controlled SPDE}
        \frac{1}{2}\sup_{t\in[0,T]}\mathbb{E}\left[\int_U(\rho^{n,\epsilon,g}-M)^2\right]+\mathbb{E}\left[\int_0^T\int_U|\nabla\Theta_{\Phi,2}(\rho^{n,\epsilon,g})|^2\right]+\alpha\mathbb{E}\left[\int_0^T\int_U|\nabla\rho^{n,\epsilon,g}|^2\right]\\
        \leq \frac{1}{2}\|\rho_0-M\|^2_{L^2(U)} +cT+c\int_0^T\int_U |P_{K(\epsilon)}g|^2,
    \end{multline}
    where $\Theta_{\Phi,2}$ is as in Definition \ref{def: auxhilary functions} with the choice $p=2$.
\end{proposition}
\begin{proof}
    The proof is exactly Proposition 4.14 of \cite{popat2025well}.
   The estimate \eqref{eq: l2 space-time estimate for regularised controlled SPDE} is far more simple than the corresponding estimate in \cite{popat2025well} due to Assumption \ref{asm: constant boundary data}.\\
     We note that the results in \cite{popat2025well} correspond to the uncontrolled equation \eqref{generalised Dean-Kawasaki Equation Stratonovich with epsilon}.
     The presence of the control term in \eqref{eq: regularised controlled SPDE} leads to the additional term
    \begin{align*}
        \int_0^t\int_U\sigma(\rho^{n,\epsilon,g})\nabla\rho^{n,\epsilon,g}\cdot P_{K(\epsilon)}g.
    \end{align*}
    Using the assumption that $\sigma\leq \Phi^{1/2}$, Cauchy-Schwarz and Young's inequalities, for every $\delta\in(0,1)$ we have the bound
    \begin{align}\label{eq: how to bound new terms involving control}
        \int_0^t\int_U\sigma(\rho^{n,\epsilon,g})\nabla\rho^{n,\epsilon,g}\cdot g&\leq \delta\int_0^t\int_U\Phi(\rho^{n,\epsilon,g})|\nabla \rho^{n,\epsilon,g}|^2+\delta^{-1}\int_0^t\int_U|P_{K(\epsilon)}g|^2\nonumber\\
        &=\delta\int_0^t\int_U|\nabla\Theta_{\Phi,2}(\rho^{n,\epsilon,g})|^2+\delta^{-1}\int_0^t\int_U|P_{K(\epsilon)}g|^2.
    \end{align}
    Picking $\delta\in(0,1)$ sufficiently small so that the first term can be absorbed onto the left hand side of the estimate completes the proof.
    \end{proof}
We note that the presence of the $L^2(U)$-norm of the initial data in the above estimate \eqref{eq: l2 space-time estimate for regularised controlled SPDE} will necessitate the higher order regularity of the initial data for the subsequent estimates.\\
The proof of tightness for laws of  the regularised controlled SPDE \eqref{eq: regularised controlled SPDE} will be via the Aubin-Lions-Simon Lemma, \cite{aubin1963theoreme,lions1969quelques,simon1986compact}. 
We state a sufficient condition below that requires us to prove higher order spatial regularity and higher order time regularity.
\begin{lemma}[Sufficient condition for tightness of controlled SPDE]\label{lemma: Aubin-Lions-Simon tightness result}
  Let $\{\rho^\epsilon\}_{\epsilon\in(0,1)}$ be a sequence that satisfies for every $\beta\in(0,1/2)$ and every $s\geq\frac{d+2}{2}$, that there is a constant constant $c\in(0,\infty)$ independent of $\epsilon$ such that either
    \begin{align}\label{eq: first criteria for tightness}
        \|\rho^\epsilon\|_{L^1([0,T];W^{1,1}(U))}+\|\rho^\epsilon\|_{W^{\beta,1}([0,T];H^{-s}(U))}\leq c,
    \end{align}
    or for every $\gamma\in(0,1)$
     \begin{align}\label{eq: second criteria for tightness}
        \|\rho^\epsilon\|_{L^1([0,T];W^{\gamma,1}(U))}+\|\rho^\epsilon\|_{W^{\beta,1}([0,T];H^{-s}(U))}\leq c.
    \end{align}
    Then $\{\rho^\epsilon\}_{\epsilon\in(0,1)}$ is relatively pre-compact in $L^1(U\times[0,T])$.
\end{lemma}
\begin{proof}
    The proof is a special case of Aubin-Lions-Simon lemma, see specifically Corollary 5 of \cite{simon1986compact}, which follows from the compact embeddings of $W^{1,1}(U),W^{\gamma,1}(U)\hookrightarrow L^1(U)$ and the continuous embedding of $L^1(U)\hookrightarrow H^{-s}(U)$.
\end{proof}
Our goal will therefore be to prove the two regularity estimates for the regularised controlled SPDE.
The two different spatial regularities in equations \eqref{eq: first criteria for tightness} and \eqref{eq: second criteria for tightness} are due to whether we are considering fast diffusion type equations $\Phi(\rho)=\rho^m$ for $m\in(0,1)$ or porus medium type equations $m>1$, which we see through the two conditions in point 3 of Assumption \ref{asm: assumptions for well posedenss of equation} in the Appendix.
\begin{proposition}[Higher order spatial regularity for \eqref{eq: regularised controlled SPDE}]\label{prop: higher order spatial regularity for regularised controlled SPDE}
    Let the coefficients of \eqref{eq: regularised controlled SPDE} satisfy Assumptions \ref{assumption on noise} and \ref{asm: assumptions for well posedenss of equation}.
    Let $g\in L^2(U\times[0,T];\mathbb{R}^d)$ be an $\mathcal{F}_t$-predictable control, $\rho_0\in L^2(\Omega;L^2(U))$ be non-negative and $\mathcal{F}_0$-measurable initial data, and to simplify the below estimate we assume that the boundary data $\bar{f}\in H^1(\partial U)$ satisfies Assumption \ref{asm: constant boundary data}.
     Suppose that $\epsilon\in(0,1)$ is sufficiently small so that the inequality \eqref{eq: scaling so that noise times epsilon is less than 1} is satisfied.
    For fixed regularisation constants $\alpha\in(0,1), n\in\mathbb{N}$, let $\rho^{n,\epsilon,g}$ denote weak solutions of the regularised controlled SPDE \eqref{eq: regularised controlled SPDE}.
    \begin{enumerate}
        \item If $\Phi$ satisfies equation \eqref{eq: first assumption on Phi} in point 3 of Assumption \ref{asm: assumptions for well posedenss of equation}, then there exists a constant $c\in(0,\infty)$ independent of $\alpha,n,\epsilon$ such that
        \[\mathbb{E}\|\rho^{n,\epsilon,g}\|_{L^1([0,T];W^{1,1}(U))}\leq c\left(T+\|\rho_0\|^2_{L^2(U)}+\|P_{K(\epsilon)}g\|^2_{L^2(U\times[0,T];\mathbb{R}^d)}\right).\]
        \item If $\Phi$ satisfies equation \eqref{eq: second assumption on Phi} in point 3 of Assumption \ref{asm: assumptions for well posedenss of equation}, then for every $\gamma\in(0,1)$ there exists a constant $c\in(0,\infty)$ independent of $\alpha,n,\epsilon$ such that
        \[\mathbb{E}\|\rho^{n,\epsilon,g}\|_{L^1([0,T];W^{\gamma,1}(U))}\leq c\left(T+\|\rho_0\|^2_{L^2(U)}+\|P_{K(\epsilon)}g\|^2_{L^2(U\times[0,T];\mathbb{R}^d)}\right).\]
    \end{enumerate}
\end{proposition}
\begin{proof}
    The result is a consequence of Lemma 4.16 and Corollary 4.17 of \cite{popat2025well}.
    The bound for the $L^1(U\times[0,T])$-norm of the solutions follows from the technical Proposition 4.12 of \cite{popat2025well}.
    The bound for the higher order spatial regularity is precisely the technical Lemma 4.16 in \cite{popat2025well} (see Lemma 5.11 of \cite{fehrman2024well} for a proof) combined with the $L^2(U\times[0,T])$-estimate \eqref{eq: l2 space-time estimate for regularised controlled SPDE} from Proposition \ref{prop: l2 estimate from 4.14 of popat} above.
\end{proof}
There is also a technical point that arises when going from the torus to a bounded domain which we hid in the above proof but emphasise below for completeness.
\begin{remark}[Extension operator]\label{rmk: extension operator}
    In order to bound the fractional Sobolev semi-norm in Lemma 4.16 of \cite{popat2025well} (see Lemma 5.11 of \cite{fehrman2024well}), we use the identity for $x,y\in U$,
\begin{equation}\label{eq: issue with bounded domain}
    |\Theta_{\Phi,2}(\rho^{n,\epsilon,g})(x)-\Theta_{\Phi,2}(\rho^{n,\epsilon,g})(y)|=\left|\int_0^1\nabla\Theta_{\Phi,2}(\rho^{n,\epsilon,g})(y+s(x-y))\cdot(x-y)\,ds\right|.
\end{equation}
To ensure that the right hand side makes sense, it appears that we require our domain to be convex, since we require $(y+s(x-y))$, $s\in(0,1)$ to be in the domain $U$ whenever $x,y\in U$.
However, this is not necessary in light of the Sobolev extension theorem, see Chapter 5.4 of \cite{evans2022partial}, and the fact that $\Theta_{\Phi,2}\in H^1(U)$.
The theorem provides the existence of a bounded, linear extension operator
\[E:H^1(U)\to H^1(\mathbb{R}^d),\]
and since our domain $U$ is $C^2$-regular, we have the existence of a constant $c\in (0,\infty)$ depending only on the domain $U$ such that
\[\|E\Theta_{\Phi,2}\|_{H^1(\mathbb{R}^d)}\leq c\|\Theta_{\Phi,2}\|_{H^1(U)}. \]
The analysis in Lemma 5.11 of \cite{fehrman2024well} can therefore be repeated with $\Theta_{\Phi,2}$ replaced with $E\Theta_{\Phi,2}$ without any issue, and in particular there are no issues with interpreting the righ hand side of equation \eqref{eq: issue with bounded domain}.
\end{remark}
We will now state the additional time regularity result. 
Since we know that the It\^o-to-Stratonovich conversion results in a singular term when the singular solution $\rho^{\epsilon,g}$ approaches it's zero set, it follows that we only have additional space-time regularity that is uniform in the $n$-regularisation for a non-linear function of the regularised solution $\rho^{n,\epsilon,g}$ that cuts it off from its zero set.  
\begin{proposition}[Higher order time regularity for \eqref{eq: regularised controlled SPDE}]\label{prop: higher order time regularity for regularised controlled SPDE}
      Let the coefficients of \eqref{eq: regularised controlled SPDE} satisfy Assumptions \ref{assumption on noise} and \ref{asm: assumptions for well posedenss of equation}.
   Let $g\in L^2(U\times[0,T];\mathbb{R}^d)$ be an $\mathcal{F}_t$-predictable control, $\rho_0\in L^2(\Omega;L^2(U))$ be non-negative and $\mathcal{F}_0$-measurable initial data, and to simplify the below estimate we assume that the boundary data $\bar{f}\in H^1(\partial U)$ satisfies Assumption \ref{asm: constant boundary data}.
       Suppose that $\epsilon\in(0,1)$ is sufficiently small so that the inequality \eqref{eq: scaling so that noise times epsilon is less than 1} is satisfied.
    For fixed regularisation constants $\alpha\in(0,1),n\in\mathbb{N}$, let $\rho^{n,\epsilon,g}$ denote weak solutions of the regularised controlled SPDE \eqref{eq: regularised controlled SPDE}.\\
    For every $\beta\in(0,1)$, let $\Psi_\beta$ be defined by
    $\Psi_\beta(\xi):=\xi\psi_\beta(\xi)$,
    where $\psi_\beta$ is a smooth cutoff function around zero of strength $\beta$, see Definition 4.18 of \cite{popat2025well} for precise definition.
    For every $\beta\in(0,1)$, $\delta\in(0,1/2)$ and $s>\frac{d+1}{2}$ there exists a constant $c\in(0,\infty)$ independent of $\alpha,n,\epsilon$ such that
    \[\mathbb{E}\|\Psi_\beta(\rho^{n,\epsilon,g})\|_{W^{\delta,1}([0,T];H^{-s}(U))}\leq c\left(T+\|\rho_0\|^2_{L^2(U)}+\|P_{K(\epsilon)}g\|^2_{L^2(U\times[0,T];\mathbb{R}^d)}\right).\]
\end{proposition}
\begin{proof}
    See Proposition 4.19 of \cite{popat2025well} for a proof for the uncontrolled equation.
    The presence of the additional control term here does not cause additional difficulty. 
The term is
\begin{multline*}
    \int_0^t\Psi'_\beta(\rho^{n,\epsilon,g})\nabla\cdot(\sigma_n(\rho^{n,\epsilon,g})P_{K(\epsilon)}g)\\
    =\int_0^t\nabla\cdot(\Psi'_\beta(\rho^{n,\epsilon,g})\sigma_n(\rho^{n,\epsilon,g})P_{K(\epsilon)}g)-\int_0^t\Psi''_\beta(\rho^{n,\epsilon,g})\sigma_n(\rho^{n,\epsilon,g})\nabla\rho^{n,\epsilon,g}\cdot P_{K(\epsilon)}g.
\end{multline*}
The $W^{1,1}([0,T];H^{-s}(U))$-norm of the above terms can be bounded exactly in the same way as the other finite variations terms of Proposition 4.19 of \cite{popat2025well} and the methods of equation \eqref{eq: how to bound new terms involving control}, which results in the additional factor of the $L^2(U\times[0,T])$-norm of $P_{K(\epsilon)}(g)$ on the right hand side of the estimate.
\end{proof}
\begin{corollary}[Tightness of laws of controlled SPDE]\label{corr: tightness of laws of controlled SPDE}
Let the coefficients of \eqref{eq: regularised controlled SPDE} satisfy Assumptions \ref{assumption on noise} and \ref{asm: assumptions for well posedenss of equation}. 
Let $g\in L^2(U\times[0,T];\mathbb{R}^d)$ be an $\mathcal{F}_t$-predictable control, $\rho_0\in L^2(\Omega;L^2(U))$ be non-negative and $\mathcal{F}_0$-measurable initial data, and $\bar{f}\in H^1(\partial U)$ be boundary data satisfying Assumption \ref{asm: constant boundary data}.
If $\epsilon\in(0,1), K(\epsilon)\in\mathbb{N}$ satisfy the joint scaling \eqref{eq: choice of joint scaling} from Remark \ref{rmk: choice of joint scaling},
 then the laws of the solutions $\{\rho^{n,\epsilon,g}\}_{\alpha,\epsilon\in(0,1),n\in\mathbb{N}}$ of the regularised SPDE \eqref{eq: regularised controlled SPDE} with initial condition $\rho_0$, control $g$ and boundary data $\bar{f}$ are tight\footnote{Recall that the tightness is along appropriate sub-sequences of all three $\alpha\to0, n\to\infty, \epsilon\to0$, which follows from the fact that the right hand side of the energy estimates above are stable in all three $\alpha,n,\epsilon$.} on $L^1(U\times[0,T])$ with respect to the strong norm topology.    
\end{corollary}
\begin{proof}
The proof is a consequence of Lemma \ref{lemma: Aubin-Lions-Simon tightness result} and Propositions \ref{prop: higher order spatial regularity for regularised controlled SPDE} and \ref{prop: higher order time regularity for regularised controlled SPDE}.\\
The only delicate point is that in Proposition \ref{prop: higher order time regularity for regularised controlled SPDE} we only showed higher order spatial regularity for the cutoff equation $\Psi_\beta(\rho^{n,\epsilon,g}), \beta\in(0,1)$.
The tightness of the laws of $\rho^{n,\epsilon,g}$ follows by the introduction of a new metric on $L^1(U\times[0,T])$ that takes into account the non-linear approximation $\Psi_\beta$, see Definitions 5.23 of \cite{fehrman2024well} or Definition 4.33 of \cite{popat2025well}.
The new metric induces the same topology as the usual strong norm topology on $L^1(U\times[0,T])$, and furthermore proving tightness of the laws of $\Psi_\beta(\rho^{n,\epsilon,g})$ with respect to the new metric is equivalent to tightness of the solutions $\rho^{n,\epsilon,g}$ in law.
The precise proof of tightness is identical to Proposition 5.26 of \cite{fehrman2024well}, so we will not repeat it here.
\end{proof}
We end this section with an entropy estimate for the regularised equation.
The estimate is used to control gradient terms when we show convergence of the controlled SPDE \eqref{eq: controlled SPDE} to the skeleton equation \eqref{eq: skeleton equation} in Theorem \ref{thm: convergence of solutions when initial data and controls converge} below. 
Furthermore, the first condition of the definition of weak solutions for the skeleton equation, Definition \ref{def: weak solution of skeleton equation}, requires $H^1(U)$-spatial regularity of $\Phi^{1/2}(\rho)$, which is precisely given by the below estimate.
\begin{proposition}[Entropy estimate]\label{ppn: entropy estimate}
      Let the coefficients of \eqref{eq: regularised controlled SPDE} satisfy Assumptions \ref{assumption on noise} and \ref{asm: assumptions for well posedenss of equation}.
    Let $\Psi_\Phi$ be as in Definition \ref{def: entropy space}, $g\in L^2(U\times[0,T];\mathbb{R}^d)$ be an $\mathcal{F}_t$-predictable control, $\rho_0\in L^2(\Omega;L^2(U))$ be non-negative and $\mathcal{F}_0$-measurable initial data, and to simplify the below estimate we assume that the boundary data $\bar{f}\in H^1(\partial U)$ satisfies Assumption \ref{asm: constant boundary data}.
     Suppose that $\epsilon\in(0,1)$ is sufficiently small so that the inequality \eqref{eq: scaling so that noise times epsilon is less than 1} is satisfied.
    For fixed regularisation constants $\alpha\in(0,1), n\in\mathbb{N}$, let $\rho^{n,\epsilon,g}$ denote weak solutions of the regularised controlled SPDE \eqref{eq: regularised controlled SPDE}.
    Then we have the estimate
    \begin{multline*}
         \mathbb{E}\left[\int_U\Psi_{\Phi}(\rho^{n,\epsilon,g})\,dx\right]+ 4\mathbb{E}\left[\int_0^T\int_U|\nabla\Phi^{1/2}(\rho^{n,\epsilon,g})|^2\right]+\alpha \mathbb{E}\left[\int_0^T\int_U\frac{\Phi'(\rho)}{\Phi(\rho)}|\nabla\rho^{n,\epsilon,g}|^2\right]\\
         \leq \mathbb{E}\left[\int_U\Psi_{\Phi}(x,\rho_0(x))\,dx\right]+c\left(T+\|\rho_0\|^2_{L^2(U)}+\|P_{K(\epsilon)}g\|^2_{L^2(U\times[0,T];\mathbb{R}^d)}\right).
    \end{multline*}
\end{proposition}
\begin{proof}
    The estimate is a consequence of the entropy estimate in Proposition 4.24 of \cite{popat2025well}, which gives the existence of a constant $c\in(0,\infty)$ such that
    \begin{multline}\label{eq: entropy estimate from pop25}
        \mathbb{E}\left[\int_U\Psi_{\Phi}(\rho^{n,\epsilon,g})\,dx\right]+ 4\mathbb{E}\left[\int_0^T\int_U|\nabla\Phi^{1/2}(\rho^{n,\epsilon,g})|^2\right]+\alpha \mathbb{E}\left[\int_0^T\int_U\frac{\Phi'(\rho)}{\Phi(\rho)}|\nabla\rho^{n,\epsilon,g}|^2\right]\\
         \leq cT+\mathbb{E}\left[\int_U\Psi_{\Phi}(x,\rho_0(x))\,dx\right]+\mathbb{E}\left[\int_0^T\int_U|\nabla\Theta_{\Phi,2}(\rho^{n,\epsilon,g})|^2\right]\\
         +\mathbb{E}\left[\lim_{\delta\to0}\int_0^T\int_U\frac{\Phi'(\rho^{n,\epsilon,g})}{\Phi(\rho^{n,\epsilon,g})+\delta}\sigma_n(\rho^{n,\epsilon,g})\nabla\rho^{n,\epsilon,g}\cdot P_{K(\epsilon)} g\right],
    \end{multline}
     where $\Theta_{\Phi,2}$ is as in Definition \ref{def: auxhilary functions} above with the choice $p=2$.
     A bound for the third term on the right hand side of \eqref{eq: entropy estimate from pop25} is obtained using Proposition \ref{prop: l2 estimate from 4.14 of popat} above.
     Furthermore, due to the presence of the control we have the final term of \eqref{eq: entropy estimate from pop25} that needs to be handled.
     An argument similar to \eqref{eq: how to bound new terms involving control} using the assumption that $\sigma\leq c\Phi^{1/2}$, Cauchy-Schwarz and Young's inequalities and the fact that $\frac{x}{x+\delta}\leq 1$ for every $\delta\in(0,1)$ gives for every $\gamma\in(0,1)$
     \begin{multline*}
         \int_0^T\int_U\frac{\Phi'(\rho^{n,\epsilon,g})}{\Phi(\rho^{n,\epsilon,g})+\delta}\sigma_n(\rho^{n,\epsilon,g})\nabla\rho^{n,\epsilon,g}\cdot P_{K(\epsilon)} g\\
         \leq c\gamma\int_0^T\int_U \frac{(\Phi'(\rho^{n,\epsilon,g}))^2}{\left(\Phi(\rho^{n,\epsilon,g})+\delta\right)}|\nabla\rho^{n,\epsilon,g}|^2+c\gamma^{-1}\int_0^T\int_U|P_{K(\epsilon)}g|^2.
     \end{multline*}
     The first term can be absorbed onto the left hand side of \eqref{eq: entropy estimate from pop25} after taking the limit $\delta\to0$ by picking $\gamma>0$ sufficiently small.
     Putting everything together then gives the estimate.
\end{proof}

\subsection{The uniform large deviations principle}\label{sec: checking conditions for BDM}
We are in a position to prove the the remaining conditions of Theorem \ref{thm: LDP theorem from BDM}. 
\begin{theorem}\label{thm: convergence of solutions when initial data and controls converge}
   Let the coefficients of equation \eqref{eq: regularised controlled SPDE} satisfy Assumptions \ref{assumption on noise}, \ref{asm: assumptions for well posedenss of equation} and \ref{asm: assumptions on boundary data for well-posedness}, and furthermore assume that the non-linearity $\Phi$ satisfies for constant $c\in(0,\infty)$ and every $M\in(0,\infty)$
   \begin{equation}\label{eq: assumption on quotient Phi / Phi'}
       \sup_{0\leq \xi\leq M}\frac{\Phi(\xi)}{\Phi'(\xi)}\leq cM.
   \end{equation}
    Let $\{g^\epsilon\}_{\epsilon\in(0,1)}$ and $g$ be $L^2(\Omega; L^2(U\times[0,T];\mathbb{R}^d))$-valued, $\mathcal{F}_t$-predictable processes satisfying \[\sup_{\epsilon\in(0,1)}\|g^\epsilon\|_{L^\infty(\Omega; L^2(U\times[0,T];\mathbb{R}^d))}<\infty\quad \text{and} \quad g^\epsilon\to g \quad \text{weakly in} \quad L^2(U\times[0,T];\mathbb{R}^d) \quad \text{as} \quad \epsilon\to0,\]
    and let $\{\rho_0^\epsilon\}_{\epsilon\in(0,1)},\rho_0\in Ent_\Phi(U)$ be such that 
    \begin{equation}\label{eq: weak convergence of initial data}
    \sup_{\epsilon\in(0,1)}\left(\int_U\Psi_\Phi(\rho_0^\epsilon)\right)<\infty \quad\text{and}\quad \rho_0^\epsilon\to\rho_0 \quad\text{weakly in}\quad L^1(U) \quad \text{as} \quad \epsilon\to0. 
    \end{equation}
    Let $\{\epsilon, K(\epsilon)\}_{\epsilon\in(0,1)}$ be such that the joint scaling \eqref{eq: choice of joint scaling} in Remark \ref{rmk: choice of joint scaling} is satisfied. 
    For every $\epsilon\in(0,1)$ and fixed $\bar{f}\in H^1(\partial U)$, denote by $\rho^{\epsilon,g^\epsilon}(\bar{f},\rho^\epsilon_0)$ the unique stochastic kinetic solutions of the controlled SPDE \eqref{eq: controlled SPDE} with control $g^\epsilon$, boundary data $\bar{f}$, and initial data $\rho_0^\epsilon$.
    Then we have that
    \[\rho^{\epsilon,g^\epsilon}(\bar{f},\rho_0^\epsilon)\to\rho(\bar{f},\rho_0) \quad\text{weakly in}\quad L^1(U\times[0,T]) \quad \text{as} \quad \epsilon\to0,\]
    where as above $\rho$ denotes the solution of the skeleton equation \eqref{eq: skeleton equation} in the sense of Definition \ref{def: weak solution of skeleton equation} and Theorem \ref{thm: well-posedness of skeleton equation assumed}.
\end{theorem}
\begin{proof}
Let us begin by considering non-negative initial data $\{\rho_0^\epsilon\}_{\epsilon\in(0,1)}$ and $\rho_0$ that are uniformly bounded in $L^2(U)$. 
In this case we can use Corollary \ref{corr: tightness of laws of controlled SPDE} that gives us tightness of laws of the regularised controlled SPDEs $\{\rho^{n,\epsilon,g^\epsilon}(\bar{f},\rho_0^\epsilon)\}_{\alpha,\epsilon\in(0,1),n\in\mathbb{N}}$.\\
We refer a reader who is unfamiliar with the kinetic formulation to see Section \ref{sec: kinetic formulation for controlled SPDE} of the Appendix for a derivation of the kinetic equation of the controlled SPDE \eqref{eq: controlled SPDE}.
The sequence of non-negative kinetic measures $\{q^\epsilon\}_{\epsilon\in(0,1)}$ corresponding to the regularised equation \eqref{eq: regularised controlled SPDE} are defined for fixed regularisation constants $\alpha\in(0,1),n\in\mathbb{N}$ by 
\[dq^\epsilon:=\delta_0(\xi-\rho^{n,\epsilon,g^\epsilon})\left(|\nabla\Theta_{\Phi,2}(\rho^{n,\epsilon,g^\epsilon})|^2+\alpha|\nabla \rho^{n,\epsilon,g^\epsilon}|^2\right)\,dx\,dt\,d\xi.\]
Due to Proposition \ref{prop: l2 estimate from 4.14 of popat} and the assumed the $L^2(U)-$integrability of the initial data, the non-negative kinetic measures are uniformly bounded in expectation, in the sense that
\[\sup_{\epsilon\in(0,1)}\mathbb{E}\left[q^\epsilon(U\times(0,\infty)\times[0,T])\right]<\infty.\]
Let the kinetic function $\chi^\epsilon$ corresponding to the regularised controlled equation $\rho^{n,\epsilon,g^\epsilon}$ be defined in Definition \ref{def: kinetic function} in the Appendix.
In Section \ref{sec: kinetic formulation for controlled SPDE} of the Appendix we showed that for every $t\in[0,T]$, $\psi\in C_c^\infty(U\times(0,\infty))$, $\mathbb{P}-a.s.$ the kinetic equation for fixed regularised coefficients $\alpha\in(0,1),n\in\mathbb{N}$ is the equation
\begin{multline}\label{eq: kinetic equation in proof of weak convergence}
\int_{\mathbb{R}}\int_U\chi^\epsilon(x,\xi,t)\psi(x,\xi)\\
=\int_{\mathbb{R}}\int_U\chi^\epsilon(x,\xi,0)\psi(x,\xi)-\int_0^t\int_U\Phi'(\rho^{n,\epsilon, g^\epsilon})\nabla\rho^{n,\epsilon, g^\epsilon}\cdot(\nabla_x\psi)(x,\rho^{n,\epsilon, g^\epsilon})-\int_0^t\int_\mathbb{R}\int_U\partial_\xi\psi(x,\xi)\,dq^\epsilon\\
    -\sqrt{\epsilon}\int_0^t\int_U\psi(x,\rho^{n,\epsilon, g^\epsilon})\nabla\cdot(\sigma_n(\rho^{n,\epsilon, g^\epsilon}) \,d{\xi}^{K(\epsilon)})+\int_0^t\int_U\sigma_n(\rho^{n,\epsilon, g^\epsilon})(\nabla_x\psi)(x,\rho^{n,\epsilon, g^\epsilon})\cdot P_{K(\epsilon)} g^\epsilon\\
    +\int_0^t\int_U(\partial_\xi\psi)(x,\rho^{n,\epsilon, g^\epsilon})\sigma_n(\rho^{n,\epsilon, g^\epsilon})\nabla\rho^{n,\epsilon, g^\epsilon}\cdot P_{K(\epsilon)}g-\int_0^t\int_U\psi(x,\rho^{n,\epsilon, g^\epsilon})\nabla\cdot\nu(\rho^{n,\epsilon, g^\epsilon})\\
    -\frac{\epsilon}{2}\int_0^t\int_U(\nabla_x\psi)(x,\rho^{n,\epsilon, g^\epsilon})\cdot\left(F_1^{K(\epsilon)}(\sigma_n'(\rho^{n,\epsilon, g^\epsilon}))^2\nabla\rho^{n,\epsilon, g^\epsilon} + \sigma_n(\rho^{n,\epsilon, g^\epsilon})\sigma_n'(\rho^{n,\epsilon, g^\epsilon})F_2^{K(\epsilon)}\right)\\
    +\frac{\epsilon}{2}\int_0^t\int_U(\partial_\xi\psi)(x,\rho^{n,\epsilon, g^\epsilon})\left(\sigma_n(\rho^{n,\epsilon, g^\epsilon})\sigma_n'(\rho^{n,\epsilon, g^\epsilon})\nabla\rho^{n,\epsilon, g^\epsilon}\cdot F_2^{K(\epsilon)}+\sigma_n^2(\rho^{n,\epsilon, g^\epsilon})F_3^{K(\epsilon)}\right).
\end{multline}
The terms in the final two lines involving a products of factors of $\epsilon$ and noise coefficients $\{F_i^{K(\epsilon)}\}_{i=1,2,3}$ dictate the joint scaling needed for the large deviations principle, and we note that it is the first term in the penultimate line involving a factor of $(\sigma_n'(\rho^{n,\epsilon, g^\epsilon}))^2$ that necessitates the use of compactly supported test functions in the velocity variable $\xi$.\\
We note that the energy estimate of Proposition \ref{ppn: entropy estimate} only give us weak convergence of the random variables $\nabla\Phi^{1/2}(\rho^{n,\epsilon, g^\epsilon})$ in the space $L^2(U\times[0,T];\mathbb{R}^d)$ rather weak convergence of the gradients $\nabla\rho^{n,\epsilon, g^\epsilon}$. 
Hence using the distributional equality
\[\nabla\Phi^{1/2}(\rho^{n,\epsilon, g^\epsilon})=\frac{1}{2}\Phi^{-1/2}(\rho^{n,\epsilon, g^\epsilon})\Phi'(\rho^{n,\epsilon, g^\epsilon})\nabla\rho^{n,\epsilon, g^\epsilon}\]
we re-write \eqref{eq: kinetic equation in proof of weak convergence} as
\small
\begin{multline}\label{eq: kinetic equation in proof of weak convergence involving gradients of Phi 1/2}
\int_{\mathbb{R}}\int_U\chi^\epsilon(x,\xi,t)\psi(x,\xi)\\
=\int_{\mathbb{R}}\int_U\chi^\epsilon(x,\xi,0)\psi(x,\xi)-2\int_0^t\int_U\Phi^{1/2}(\rho^{n,\epsilon, g^\epsilon})\nabla\Phi^{1/2}(\rho^{n,\epsilon, g^\epsilon})\cdot(\nabla_x\psi)(x,\rho^{n,\epsilon, g^\epsilon})-\int_0^t\int_\mathbb{R}\int_U\partial_\xi\psi(x,\xi)\,dq^\epsilon\\
    -\sqrt{\epsilon}\int_0^t\int_U\psi(x,\rho^{n,\epsilon, g^\epsilon})\nabla\cdot(\sigma_n(\rho^{n,\epsilon, g^\epsilon}) \,d{\xi}^{K(\epsilon)})+\int_0^t\int_U\sigma_n(\rho^{n,\epsilon, g^\epsilon})(\nabla_x\psi)(x,\rho^{n,\epsilon, g^\epsilon})\cdot P_{K(\epsilon)} g^\epsilon\\
    +\int_0^t\int_U(\partial_\xi\psi)(x,\rho^{n,\epsilon, g^\epsilon})\frac{2\Phi^{1/2}(\rho^{n,\epsilon, g^\epsilon})\sigma_n(\rho^{n,\epsilon, g^\epsilon})}{\Phi'(\rho^{n,\epsilon, g^\epsilon})}\nabla\Phi^{1/2}(\rho^{n,\epsilon, g^\epsilon})\cdot P_{K(\epsilon)}g-\int_0^t\int_U\psi(x,\rho^{n,\epsilon, g^\epsilon})\nabla\cdot\nu(\rho^{n,\epsilon, g^\epsilon})\\
    -\frac{\epsilon}{2}\int_0^t\int_U(\nabla_x\psi)(x,\rho^{n,\epsilon, g^\epsilon})\cdot\left(F_1^{K(\epsilon)}\frac{2\Phi^{1/2}(\rho^{n,\epsilon, g^\epsilon})(\sigma_n'(\rho^{n,\epsilon, g^\epsilon}))^2}{\Phi'(\rho^{n,\epsilon, g^\epsilon})}\nabla\Phi^{1/2}(\rho^{n,\epsilon, g^\epsilon}) + \sigma_n(\rho^{n,\epsilon, g^\epsilon})\sigma_n'(\rho^{n,\epsilon, g^\epsilon})F_2^{K(\epsilon)}\right)\\
    +\frac{\epsilon}{2}\int_0^t\int_U(\partial_\xi\psi)(x,\rho^{n,\epsilon, g^\epsilon})\left(\frac{2\Phi^{1/2}(\rho^{n,\epsilon, g^\epsilon})\sigma_n(\rho^{n,\epsilon, g^\epsilon})\sigma_n'(\rho^{n,\epsilon, g^\epsilon})}{\Phi'(\rho^{n,\epsilon, g^\epsilon})}\nabla\Phi^{1/2}(\rho^{n,\epsilon, g^\epsilon})\cdot F_2^{K(\epsilon)}+\sigma_n^2(\rho^{n,\epsilon, g^\epsilon})F_3^{K(\epsilon)}\right),
\end{multline}
\normalsize
where the quotient terms are integrable due to the compact support of the test function $\psi$ in the velocity variable.
Due to the tightness of laws of $\rho^{n,\epsilon,g^\epsilon}(\bar{f},\rho_0^\epsilon)$ given by Corollary \ref{corr: tightness of laws of controlled SPDE}, we want to use Skorokhod representation theorem to pass almost surely to the to the $\alpha,\epsilon\to0, n\to\infty$ limits on an auxiliary probability space.
However, the reason that we can not do this directly is due to the first term on the fourth line of \eqref{eq: kinetic equation in proof of weak convergence involving gradients of Phi 1/2}, where we have the product of two random variables $\nabla\Phi^{1/2}(\rho^{n,\epsilon, g^\epsilon})$ and $P_{K(\epsilon)}g$ which are both only weakly convergent, so it is not possible to characterise convergence of their product.
The key observation is that such a product does not appear in the weak formulation of the skeleton equation, see Definition \ref{def: weak solution of skeleton equation}.
Let us upper bound the contribution of this term by introducing the measures $p^\epsilon$.
For $\epsilon\in(0,1)$ we define $p^{\epsilon}$ to be the non-negative, almost surely finite measure on $U\times(0,\infty)\times[0,T]$ given by
\[d p^\epsilon:=\delta_0(\xi-\rho^{n,\epsilon,g^\epsilon})|\nabla\Phi^{1/2}(\rho^{n,\epsilon,g^\epsilon})|\,|P_{K(\epsilon)}g^\epsilon|\,dx\,dt\,d\xi.\]
The finiteness of the measures follows from the $L^2(U\times[0,T])$-boundedness of the two terms\\ $|\nabla\Phi^{1/2}(\rho^{n,\epsilon,g^\epsilon})|$ and $|P_{K(\epsilon)}g^\epsilon|$. 
By the assumption $\sigma\leq c\Phi^{1/2}$ and by construction, we have the bound
\begin{multline}\label{eq: bound for measure p}
\left|\int_0^t\int_U(\partial_\xi\psi)(x,\rho^{n,\epsilon, g^\epsilon})\frac{2\Phi^{1/2}(\rho^{n,\epsilon, g^\epsilon})\sigma_n(\rho^{n,\epsilon, g^\epsilon})}{\Phi'(\rho^{n,\epsilon, g^\epsilon})}\nabla\Phi^{1/2}(\rho^{n,\epsilon, g^\epsilon})\cdot P_{K(\epsilon)}g\right|\\
\leq c\int_\mathbb{R_+}\int_0^T\int_U |\partial_\xi\psi(x,\xi)|\frac{2\Phi(\xi)}{\Phi'(\xi)}\,dp^\epsilon.
\end{multline}
To understand the limiting behaviour of the stochastic integral in equation \eqref{eq: kinetic equation in proof of weak convergence involving gradients of Phi 1/2}, again omitting the dependence on the regularisation coefficients, define for fixed $\alpha\in(0,1),n\in\mathbb{N}$ and test function $\psi\in C^\infty_c(U\times(0,\infty))$
\[M_t^{\psi}:=\int_0^t\int_U\psi(x,\rho^{n,\epsilon, g^\epsilon})\nabla\cdot(\sigma_n(\rho^{n,\epsilon, g^\epsilon}) \dot{\xi}^{K(\epsilon)}).\]
Following a similar method to Proposition \ref{prop: higher order time regularity for regularised controlled SPDE} above, Proposition 5.27 of \cite{fehrman2024well} proves that for $\gamma\in(0,1/2)$, the laws of the martingales are tight on $C^\gamma([0,T])$.\\
We now follow the methods of the uniqueness proof, see Theorem 5.29 of \cite{fehrman2024well} to characterise the limiting behaviour of the solutions $\rho^{n,\epsilon,g^\epsilon}$.
To do so, we need to characterise the limiting behaviour of each of the components in the kinetic equation \eqref{eq: kinetic equation in proof of weak convergence involving gradients of Phi 1/2} above. 
For $s>\frac{d+2}{2}$ and $\{\psi_j\}_{j\in\mathbb{N}}$ a countable sequence of dense functions in $C_c^\infty(U\times(0,\infty))$ with respect to the $H^s(U\times(0,\infty))$-topology, we want to establish the tightness of the random variables
\begin{align*}
    X^{\alpha,n,\epsilon}:=\left(\rho^{n,\epsilon,g^\epsilon},\rho_0^\epsilon, \nabla\Theta_{\Phi,2}(\rho^{n,\epsilon,g^\epsilon}),\alpha\nabla\rho^{n,\epsilon,g^\epsilon},g^\epsilon, \frac{q^\epsilon}{|q^\epsilon|},\frac{p^\epsilon}{|p^\epsilon|},|q^\epsilon|,|p^\epsilon|,\sqrt{\epsilon} (M_t^{\psi_j})_{j\in\mathbb{N}}\right),
\end{align*}
    where $|q^\epsilon|:=q^\epsilon(U\times(0,\infty)\times[0,T])$ and analogously for $|p^\epsilon|$, on the product metric topology of the state space
    \begin{multline*}
        \bar{X}:=L^1(U\times[0,T])\times L^2(U)\times L^2(U\times[0,T];\mathbb{R}^d)^3\times\mathcal{P}(U\times(0,\infty)\times[0,T])^2
        \times\mathbb{R}^2\times C^\gamma([0,T])^\mathbb{N},
    \end{multline*}
    where $\mathcal{P}(U\times(0,\infty)\times[0,T])$ denotes the space of non-negative probability measures.
    The space $L^1(U\times[0,T])$ is equipped with the strong topology, the spaces $L^2(U)$ and $L^2(U\times[0,T];\mathbb{R}^d)$ is equipped with the weak topology, and the space $C^\gamma([0,T])^{\mathbb{N}}$ is equipped with the topology of component-wise convergence in the strong norm induced by the metric discussed in the proof of Corollary \ref{corr: tightness of laws of controlled SPDE}, see Definition 4.33 of \cite{popat2025well}.
    In much the same way as Theorem 5.29 of \cite{fehrman2024well}, by using Prokhorov's theorem and Skorokhod representation theorem, owing to Corollary \ref{corr: tightness of laws of controlled SPDE} and tightness of the martingale term, we have convergence of $X^{\alpha,n,\epsilon}$ along a subsequence $\alpha_k\to0,n_k\to\infty,\epsilon_k\to0$ in an auxiliary probability space.\\
    To prove the desired result it therefore suffices to prove that if along this subsequence, if $X^{\alpha,n,\epsilon}$ converges $\mathbb{P}-$a.s. to a random variable 
    \begin{align*}
      X:=\left(\rho,\rho_0,\nabla\Theta_{\Phi,2}(\rho),0,g, \tilde{q},\tilde{p},a,b,0\right),
    \end{align*}
    for probability measures $\tilde{q},\tilde{p}\in \mathcal{P}(U\times(0,\infty)\times[0,T])$ and random constants $a,b\in\mathbb{R}$, then $\rho$ is a $\mathbb{P}-$a.s. weak solution to the skeleton equation with control $g$, initial data $\rho_0$ and boundary data $\bar{f}$.\\
The regularisations in $\alpha,n$ were used to obtain the energy estimates that led to tightness, so do not play any further role in the proof and below we denote by $\rho^{\epsilon,g^\epsilon}$ the solution to the limiting singular controlled stochastic PDE \eqref{eq: controlled SPDE}.
   The sufficient condition is proved by passing from the kinetic equation \eqref{eq: kinetic equation in proof of weak convergence involving gradients of Phi 1/2} to the weak formulation of the equation along the subsequence $\epsilon_k\to0$.
    The kinetic equation has the additional velocity component not seen by the weak formulation, so to remove this dependency, for $\delta\in(0,1/2)$, let $\phi_\delta:[0,\infty)\to[0,1]$ be a smooth function satisfying for some constant $c\in(0,\infty)$
    \begin{equation}\label{eq: definition of phi delta velocity cutoff}
   \begin{cases}
    \phi_\delta=0&\text{on}\hspace{5pt} [0,\delta]\cup [2\delta^{-1},\infty)\\
  \phi_\delta=1&\text{on}\hspace{5pt} [2\delta,\delta^{-1}]\\
   |\phi'_\delta|\leq c\delta^{-1}\mathbbm{1}_{\xi\in[\delta,2\delta]}+c\delta\mathbbm{1}_{\xi\in[\delta^{-1},2\delta^{-1}]}.
   \end{cases} 
\end{equation}
 For an arbitrary test function $\psi\in C_c^\infty(U)$, we will consider test functions of the form $\Psi_\delta(x,\xi)=\phi_\delta(\xi)\psi(x)$.
 We first pass to the limit $\epsilon\to0$ and subsequently pass to the velocity limit $\delta\to0$ to recover the weak formulation of the skeleton equation.
 We proceed term by term.\\
 For the initial data $\rho_0^\epsilon$, we only have that it converges weakly to $\rho_0$, which is not compatible with the non-linear convergence of the kinetic function, so using the definition of the kinetic function, we write the relevant term as \say{what we want} plus a correction,
 \begin{equation*}
     \int_{\mathbb{R}}\int_U\chi^\epsilon(x,\xi,0)\Psi_\delta(x,\xi)=\int_U \rho_0^\epsilon(x)\psi(x)+\int_{\mathbb{R}}\int_U\chi^\epsilon(x,\xi,0)\psi(x)\left(\phi_\delta(\xi)-1\right).
 \end{equation*}
 By re-arranging, it follows from the support of $\phi_\delta$, the boundedness of $\psi$ and the weak convergence of $\rho_0^\epsilon$ to $\rho_0$, that there exists a constant $c\in(0,\infty)$ such that
 \begin{align}\label{eq: convergence of initial condition}
     \limsup_{\epsilon\to0}\left|\int_{\mathbb{R}}\int_U\chi^\epsilon(x,\xi,0)\Psi_\delta(x,\xi)-\int_U \rho_0(x)\psi(x)\right|\leq c\delta.
 \end{align}
 The choice of scaling \eqref{eq: choice of joint scaling},  Proposition \ref{prop: l2 estimate from 4.14 of popat} and the compact support in spatial component of $\Psi_\delta$ (which ensures that the first term in \eqref{eq: convergence of mass terms 1} remains integrable for every fixed $\delta\in(0,1)$) then implies that the terms in the penultimate line of \eqref{eq: kinetic equation in proof of weak convergence involving gradients of Phi 1/2} satisfy $\mathbb{P}-a.s.$ for every fixed $\delta\in(0,1)$,
 \small
 \begin{multline}\label{eq: convergence of mass terms 1}
     \lim_{\epsilon\to0}\sup_{t\in[0,T]}\left|\frac{\epsilon}{2}\int_0^t\int_U(\nabla_x\Psi_\delta)(x,\rho^{\epsilon, g^\epsilon})\cdot\left(F_1^{K(\epsilon)}\frac{2\Phi^{1/2}(\rho^{\epsilon, g^\epsilon})(\sigma'(\rho^{\epsilon, g^\epsilon}))^2}{\Phi'(\rho^{\epsilon, g^\epsilon})}\nabla\Phi^{1/2}(\rho^{\epsilon, g^\epsilon}) + \sigma(\rho^{\epsilon, g^\epsilon})\sigma'(\rho^{\epsilon, g^\epsilon})F_2^{K(\epsilon)}\right)\right|\\=0,
 \end{multline}
 \normalsize
 and analogously for terms in the final line of \eqref{eq: kinetic equation in proof of weak convergence involving gradients of Phi 1/2}
 \small
 \begin{multline}\label{eq: convergence of mass terms 2}
     \lim_{\epsilon\to0}\sup_{t\in[0,T]}\left|\frac{\epsilon}{2}\int_0^t\int_U(\partial_\xi\Psi_\delta)(x,\rho^{\epsilon, g^\epsilon})\left(\frac{2\Phi^{1/2}(\rho^{\epsilon, g^\epsilon})\sigma(\rho^{\epsilon, g^\epsilon})\sigma'(\rho^{\epsilon, g^\epsilon})}{\Phi'(\rho^{\epsilon, g^\epsilon})}\nabla\Phi^{1/2}(\rho^{\epsilon, g^\epsilon})\cdot F_2^{K(\epsilon)}+\sigma^2(\rho^{\epsilon, g^\epsilon})F_3^{K(\epsilon)}\right)\right|\\
     =0.
 \end{multline}
 \normalsize
 For the martingale term, we have by Burkholder-Davis-Gundy inequality and the compact support of $\Psi_\delta$ that there exists a constant $c\in(0,\infty)$ such that for every $\epsilon,\delta\in(0,1)$,
 \begin{align*}
     \mathbb{E}&\left[\sup_{t\in[0,T]}\left|\sqrt{\epsilon}\int_0^t\int_U\Psi_\delta(x,\rho^{\epsilon, g^\epsilon})\nabla\cdot(\sigma(\rho^{\epsilon, g^\epsilon}) \,d{\xi}^{K(\epsilon)})\right|\right]\\
     &\leq c\epsilon \int_0^T\int_U\Psi_\delta^2(x,\rho^{\epsilon, g^\epsilon})\left(\sigma^2(\rho^{\epsilon, g^\epsilon})F_1^{K(\epsilon)}+2\sigma(\rho^{\epsilon, g^\epsilon})\nabla\sigma(\rho^{\epsilon, g^\epsilon}) F_2^{K(\epsilon)}+|\nabla\sigma(\rho^{\epsilon, g^\epsilon})|^2 F_3^{K(\epsilon)}\right)\\
&\leq c\epsilon\left(\| F_1^{K(\epsilon)}\|_{L^\infty(U)}+\| F_2^{K(\epsilon)}\|_{L^\infty(U;\mathbb{R}^d)}+\| F_3^{K(\epsilon)}\|_{L^\infty(U)}\right)\\
&\hspace{100pt}\times\int_0^T\int_U \Psi_\delta^2(x,\rho^{\epsilon, g^\epsilon})\left(\sigma^2(\rho^{\epsilon, g^\epsilon})+(\sigma'(\rho^{\epsilon, g^\epsilon}))^2|\nabla\rho^{\epsilon, g^\epsilon}|^2\right).
 \end{align*}
 It therefore follows by the choice of scaling \eqref{eq: choice of joint scaling} that along a further subsequence, for every fixed $\delta\in(0,1)$, as $\epsilon\to0$, $\mathbb{P}-a.s.$,
 \begin{align}\label{eq: convergene of noise term}
     \lim_{\epsilon\to0}\sup_{t\in[0,T]}\left|\sqrt{\epsilon}\int_0^t\int_U\Psi_\delta(x,\rho^{\epsilon, g^\epsilon})\nabla\cdot(\sigma(\rho^{\epsilon, g^\epsilon}) \,d{\xi}^{K(\epsilon)})\right|=0.
 \end{align}
 The remaining terms are handled using the strong convergence of $\rho^{\epsilon,g^\epsilon}$ to $\rho$, the weak convergence of $P_{K(\epsilon)}g$ to $g$, the the weak convergence of $\nabla\Phi^{1/2}(\rho^{\epsilon,g^\epsilon})$ to $\nabla\Phi^{1/2}(\rho)$, the weak convergence of $q^\epsilon$ to $a\tilde{q}$ and the weak convergence of $p^\epsilon$ to $b\tilde{p}$.
 Using this, equation \eqref{eq: kinetic equation in proof of weak convergence involving gradients of Phi 1/2} and the subsequent equations \eqref{eq: bound for measure p}, \eqref{eq: convergence of initial condition}, \eqref{eq: convergence of mass terms 1}, \eqref{eq: convergence of mass terms 2}, and \eqref{eq: convergene of noise term} illustrate that after passing to a subsequence $\epsilon\to0$, $\mathbb{P}-$a.s. for almost every $t\in[0,T]$, if $\chi$ denotes the kinetic function of the skeleton equation $\rho$, then the difference between the kinetic function and the weak formulation satisfies for a constant $c\in(0,\infty)$ independent of $\delta\in(0,1)$,
\begin{multline}\label{eq: working towards convergence to weak solution of skeleton equation}
    \left|\int_{\mathbb{R}}\int_U\chi(x,\xi,t)\Psi_\delta(x,\xi)-\int_U \rho_0(x)\psi(x)+2\int_0^t\int_U\Phi^{1/2}(\rho)\nabla\Phi^{1/2}(\rho)\cdot(\nabla_x\Psi_\delta)(x,\rho)\right.\\
     \left.-\int_0^t\int_U\sigma(\rho)(\nabla_x\Psi_\delta)(x,\rho)\cdot g+\int_0^t\int_U\Psi_\delta(x,\rho)\nabla\cdot\nu(\rho)
    \right|\\
    \leq c\delta+c\left|\int_{\mathbb{R}_+}\int_0^T\int_U|\phi'_\delta(\xi)|\psi(x)\frac{2\Phi(\xi)}{\Phi'(\xi)}b\,d\tilde{p}+|\phi'_\delta(\xi)|\psi(x)a\,d\tilde{q}\right|.
\end{multline}
It follows from the support of $\phi_\delta'$, Proposition 4.29 of \cite{popat2025well} which bounds of the kinetic measure $\tilde{q}$ when the velocity argument approaches zero, the boundedness of $\psi$ and the finiteness of $b\tilde{q}$ that
\[\liminf_{\delta\to0}\max_{t\in[0,T]}\mathbb{E}\left[\int_{\mathbb{R}_+}\int_0^T\int_U|\phi'_\delta(\xi)|\psi(x)a\,d\tilde{q}\right]=0.\]
From the support of $\phi_\delta'$, the boundedness of $\psi$, the finiteness of $\tilde{p}$ and assumption \eqref{eq: assumption on quotient Phi / Phi'} in the statement of the theorem, we have
\begin{align*}
    \liminf_{\delta\to0}\max_{t\in[0,T]}\mathbb{E}&\left[\int_{\mathbb{R}_+}\int_0^T\int_U|\phi'_\delta(\xi)|\psi(x)\frac{2\Phi(\xi)}{\Phi'(\xi)}b\,d\tilde{p}\right]\\
    &\leq c \liminf_{\delta\to0}\mathbb{E}\left[\int_{[\delta,2\delta]}\int_0^T\int_U\delta^{-1}\frac{2\Phi(\xi)}{\Phi'(\xi)}\,d\tilde{p}+\int_{[\delta^{-1},2\delta^{-1}]}\int_0^T\int_U\delta\frac{2\Phi(\xi)}{\Phi'(\xi)}\,d\tilde{p}
    \right]\\
    &\leq c\liminf_{\delta\to0}\mathbb{E}\left[\int_{[\delta,2\delta]\cup [\delta^{-1},2\delta^{-1}]}\int_0^T\int_U\,d\tilde{p}\right]=0,\
\end{align*}
We are left to deal with the terms on the left hand side of \eqref{eq: working towards convergence to weak solution of skeleton equation}.
For the first, we have from the support of $\phi_\delta$ and the definition of the kinetic function that $\mathbb{P}-$a.s. for every $t\in[0,T]$,
\[\lim_{\delta\to0}\int_{\mathbb{R}}\int_U\chi^\epsilon(x,\xi,t)\Psi_\delta(x,\xi)=\int_{\mathbb{R}}\int_U\chi^\epsilon(x,\xi,t)\psi(x)=\int_{U}\rho(x,t)\psi(x),\]
and analogously for the transport term
\[\lim_{\delta\to0}\int_0^t\int_U\Psi_\delta(x,\rho)\nabla\cdot\nu(\rho)=\int_0^t\int_U\psi(x)\nabla\cdot\nu(\rho).\]
For the remaining two terms on the left hand side of \eqref{eq: working towards convergence to weak solution of skeleton equation}, using $\nabla_x\Psi(x,\rho)=\nabla\psi(x)\phi_\beta(\rho)$ it similarly follows that
\[\lim_{\delta\to0}\int_0^t\int_U\Phi'(\rho)\nabla\rho\cdot(\nabla_x\Psi_\delta)(x,\rho)=\int_0^t\int_U\Phi'(\rho)\nabla\rho\cdot\nabla\psi(x)\]
and
\[\lim_{\delta\to0}\int_0^t\int_U\sigma(\rho)(\nabla_x\Psi_\delta)(x,\rho)\cdot g=\int_0^t\int_U\sigma(\rho)\nabla\psi(x)\cdot g.\]
Putting equation \eqref{eq: working towards convergence to weak solution of skeleton equation} and subsequent computations together, after passing to a subsequence $\delta\to0$, $\mathbb{P}-$a.s. for almost every $t\in[0,T]$ we have
\begin{multline*}
    \int_U\rho(x,t)\psi(x)\\
    =\int_U \rho_0(x)\psi(x)-\int_0^t\int_U2\Phi^{1/2}(\rho)\nabla\Phi^{1/2}(\rho)\cdot\nabla\psi(x)+\int_0^t\int_U\sigma(\rho)\nabla\psi(x)\cdot g
     -\int_0^t\int_U\psi(x)\nabla\cdot\nu(\rho),
\end{multline*}
which is a re-writing of the weak formulation of the skeleton equation, see Definition \ref{def: weak solution of skeleton equation}. 
Furthermore, arguments similar to the uniqueness argument of Theorem 5.29 of \cite{fehrman2024well} we can obtain that $\rho$ has an $L^1(U\times[0,T])-$continuous representative.
The additional regularity of $\Phi^{1/2}(\rho)$ is inherited by the entropy estimate Proposition \ref{ppn: entropy estimate}, and the boundary data is constant throughout the analysis, so is inherited by the skeleton equation.
Hence $\rho$ is a $\mathbb{P}-$a.s. weak solution of the skeleton equation \eqref{eq: skeleton equation} with control $g$, initial data $\rho_0$ and boundary data $\bar{f}$.
Uniqueness then follows by Theorem \ref{thm: well-posedness of skeleton equation assumed} which gives us the equivalence of weak and renormalised kinetic solutions; and the uniqueness of kinetic solutions which follows from a simplified version of Theorem 3.6 of \cite{popat2025well}.\\
To extend the result to less regular initial data $\rho_0^\epsilon,\rho_0\in Ent_\Phi(U)$, we just apply an approximation argument. 
Approximate $\rho_0^\epsilon\in Ent_\Phi(U)$ by $(\rho_0^\epsilon\wedge n)\in L^2(U)$ for $n\in\mathbb{N}$ for which the result holds, and then the $L^1(U)$-contraction of kinetic solutions and triangle inequality gives the result.
\end{proof} 
We now check the final remaining condition in Theorem \ref{thm: LDP theorem from BDM}.
\begin{proposition}\label{prop: rate function is lower semicontinuous}
    Let $\bar{f}\in H^1(\partial U)$, $\rho_0\in Ent_\Phi(U)$ and $\rho\in L^1(U\times[0,T])$, and recall from equation \eqref{eq: rate function} that the rate function is given by
    \small
   \[ I_{\bar{f},\rho_0}(\rho)=\frac{1}{2}\inf_{g\in L^2(U\times[0,T];\mathbb{R}^d)}\left\{\|g\|^2_{L^2(U\times[0,T];\mathbb{R}^d)}:\partial_t\rho=\Delta\Phi(\rho)-\nabla\cdot(\sigma(\rho)g+ \nu(\rho)):\,\Phi(\rho)|_{\partial U}=\bar{f}, \rho(\cdot,0)=\rho_0\right\}.\]
   \normalsize
    Under Assumptions \ref{asm: assumptions for well posedenss of equation} and \ref{asm: assumptions on boundary data for well-posedness}, for every $\bar{f}\in H^1(\partial U)$ and $\rho\in L^1(U\times[0,T])$, we have $\rho_0\mapsto I_{\bar{f},\rho_0}(\rho)$ is a lower semi-continuous map from $Ent_\Phi(U)$ to $[0,\infty]$.  
\end{proposition}
\begin{proof}
    We want to show that for any sequence $\{\rho_0^\epsilon\}_{\epsilon\in(0,1)}$ uniformly bounded and weakly converging to $\rho_0$ as in \eqref{eq: weak convergence of initial data}, we have
    \[I_{\bar{f},\rho_0}(\rho)\leq \liminf_{\epsilon\to0} I_{\bar{f},\rho_0^\epsilon}(\rho).\]
For every $\epsilon\in(0,1)$,  there exists a unique minimising element $g^\epsilon$ of the rate function in the sense that
    \[I_{\bar{f}, \rho_0^\epsilon}(\rho)= \frac{1}{2}\int_0^T\int_U |g^\epsilon|^2\,dx\,ds.\]
Since $\|g^\epsilon\|_{L^2(U\times[0,T];\mathbb{R}^d)}<\infty$ uniformly in $\epsilon$, we can extract a weakly convergent subsequence $g^\epsilon\to g$ in $L^2(U\times[0,T];\mathbb{R}^d)$.
A key result is the weak-strong continuity of Theorem 21 of \cite{fehrman2023non}, which then implies that we have the strong convergence in $L^1(U\times[0,T])$ of the corresponding solutions of the skeleton equation, i.e. 
\[I_{\bar{f}, \rho_0^\epsilon}(\rho)\to I_{\bar{f}, \rho_0}(\rho)\]
strongly as $\epsilon\to0$, where the right hand side corresponds to the energy of the solution of the skeleton equation with limiting control $g$.
The weak lower-semicontinuity of the $L^2(U\times[0,T];\mathbb{R}^d)-$norm then gives the result
\begin{align*}
    I_{\bar{f},\rho_0}(\rho)\leq\frac{1}{2}\int_0^T\int_U|g|^2\,dx\,ds\leq \liminf_{\epsilon\to0}\frac{1}{2}\int_0^T\int_U|g^\epsilon|^2\,dx\,ds= \liminf_{\epsilon\to0}I_{\bar{f}, \rho_0^\epsilon}(\rho).
\end{align*}
\end{proof}
Remark \ref{rmk: verifying first point in assumption of ldp}, Propositions \ref{prop: existence of solution map for skeleton equation}, \ref{prop: existence of solution map for controlled SPDE}, Theorem \ref{thm: convergence of solutions when initial data and controls converge} and Proposition \ref{prop: rate function is lower semicontinuous} verify the conditions to prove the large deviations principle for equation \eqref{generalised Dean-Kawasaki Equation Stratonovich with epsilon} from Theorem \ref{thm: LDP theorem from BDM}.
\newpage
\section{Appendix}
\subsection{Properties of eigenfunctions of the Laplacian}\label{subsec: Appendix example of noise}
Recall that the eigenfunctions of the Laplacian with zero Dirichlet boundary conditions are defined in Example \ref{example of noise: eigenvalues of laplacian} by
\begin{equation*}
    \begin{cases}
        -\Delta e_k=\lambda_k e_k,& x\in U\\
        e_k=0,& x\in\partial U.
    \end{cases}
\end{equation*}
We have the following properties, see the useful lecture notes \cite{eigenvalues_of_lap} or Chapter 11 of \cite{strauss2007partial} for further details.
    \begin{itemize}
        \item All eigenvalues $\{\lambda_k\}_{k\in\mathbb{N}}$ are positive.
        All eigenfunctions $\{e_k\}_{k\in\mathbb{N}}$ can be chosen to be real valued, are orthogonal in $L^2(U)$, and by Gram-Schmidt they can be chosen to be orthonormal in $L^2(U)$.
        The eigenfunctions are complete in the $L^2(U)$-sense.
        \item Using the local version of Weyl's law \cite{weyl1911asymptotische}, we have the bounds 
         \begin{align}\label{eq: bounds for F_k in eigenvalue of laplacian case}
             &E^{K}_1:=\sum_{k:\lambda_k\leq\lambda_K}e_k^2\leq C K^{\frac{d}{2}};\quad E^{K}_2(x):=\sum_{k:\lambda_k\leq\lambda_K}e_k\nabla e_k\leq C K^{\frac{d+1}{2}};\nonumber\\
             &E^{K}_3(x):=\sum_{k:\lambda_k\leq\lambda_K}|\nabla e_k|^2\leq C K^{\frac{d+2}{2}}.
         \end{align}
         Furthermore, by the definition of the eigenfunctions $\{e_k\}_{k\in\mathbb{N}}$, we also have for every $K\in\mathbb{N}$,
\begin{align}\label{eq: bound for divergence of F2 in laplacian eigenfunction example}
    \nabla\cdot E^K_2=\nabla\cdot \sum_{k=1}^K e_k\nabla e_k&=\sum_{k=1}^K\left(|\nabla e_k|^2+e_k\Delta e_k\right)= \sum_{k=1}^K\left(|\nabla e_k|^2-\lambda_k e_k^2\right)\leq E_3^K.
\end{align}
    \item In order to obtain probabilistic stationary noise, we consider for $\rho_\delta$ a standard mollifier of scale $\delta$, the noise
\[\xi^{K,\delta}(x,t):=\sum_{k:\lambda_k\leq\lambda_K}(e_k(x)*\rho_\delta) B_t^k:=\sum_{k:\lambda_k\leq\lambda_K}e^\delta_k(x) B_t^k.\]
    \end{itemize}
We also have the below result that illustrates that $\{e_k\}_{k\in\mathbb{N}}$ satisfy point 2 of Assumption \ref{assumption on noise}.
\begin{lemma}\label{lemma sum of norm of e_k}
    Let $\{e_k\}_{k\in\mathbb{N}}$ be defined as in Example \ref{definition truncated noise}. 
    Then we have the equivalence
    \[\sum_{k=1}^\infty \|e_k\|^2_{H^{-s}(U)}<\infty \iff s>\frac{d}{2}. \]
\end{lemma}
\begin{proof}
The $H^{-s}(U)$-norm of the eigenfunctions of the Laplacian is computed using the definition 
\[\|e_k\|_{H^{-s}(U)}:=\sup_{f\in H^s(U):\,\|f\|_{H^s(U)}\leq 1}\left|\langle e_k,f\rangle_s\right|,\]
where the inner product $\langle\cdot,\cdot\rangle_s$ denotes the dual pairing of $H^s(U)$ and $H^{-s}(U)$, which can be interpreted as an $L^2(U)$ inner product by the Riesz Representation Theorem.
    Since $\{e_k\}_{k\in\mathbb{N}}$ is a complete orthonormal basis in $L^2(U)$, we can write a general element $f\in L^2(U)$ as
\[f(x)=\sum_{k=1}^\infty \langle f,e_k\rangle_{L^2(U)}e_k(x)=:\sum_{k=1}^\infty f_ke_k(x),\]
by which it follows that the $H^s(U)$-norm of $f$ is given by
\[\|f\|^2_{H^s(U)}=\sum_{k=1}^\infty (\lambda_k)^s |f_k|^2.\]
The maximisation problem becomes
\[\|e_k\|_{H^{-s}(U)}=\sup_{f\in H^s(U):\sum_{j\in\mathbb{N}} (\lambda_j)^s |f_j|^2\leq 1}\left|f_k\right|.\]
By inspection, the supremum is attained when $f_j=0$ for every $j\neq k$, and for the $k$'th element we have $(\lambda_k)^s|f_k|^2=1$.
Rearranging gives
\[\|e_k\|_{H^{-s}(U)}=(\lambda_k)^{-s/2}.\]
By a corollary of Weyl's law \cite{weyl1911asymptotische}, we have that on a bounded domain $U\subset \mathbb{R}^d$ the eigenvalues grow like 
\[\lambda_k\sim C k^{2/d}.\]
Hence we have 
\[\sum_{k=1}^\infty\|e_k\|^2_{H^{-s}(U)}\approx\sum_{k=1}^\infty\frac{1}{k^{2s/d}}<\infty\iff 2s/d>1 \iff s>\frac{d}{2}. \]
\end{proof}

\subsection{Assumptions for well-posedness}\label{sec: assumptions for well-posedness}
The below assumptions on the non-linear functions $\Phi,\sigma$ and $\nu$ are required for well-posedness of stochastic kinetic solutions to \eqref{generalised Dean-Kawasaki Equation Stratonovich with epsilon}, see Assumptions 3.1, 4.2, 4.20 and 4.28 of \cite{popat2025well}, which we summarise below.
We emphasise that these capture the cases of interest, $\Phi(\xi)=\xi^m$ for the full range $m\in(0,\infty)$, the critical square root $\sigma(\xi)=\sqrt{\xi}$ case, and lipschitz $\nu$.

\begin{assumption}[Assumptions on non-linear functions $\Phi,\sigma,\nu$ for well-posedness of the generalised Dean--Kawasaki equation]\label{asm: assumptions for well posedenss of equation}
       Suppose $\Phi,\sigma\in C([0,\infty))$ and $\nu\in C([0,\infty);\mathbb{R}^d)$ satisfy the following assumptions:
    \begin{enumerate}
        \item We have $\Phi,\sigma\in C^{1,1}_{loc}([0,\infty))$ and $\nu\in C^1_{loc}([0,\infty);\mathbb{R}^d)$. 
        \item The function $\Phi$ is strictly increasing and starts at $0$: $\Phi(0)=0$ with $\Phi'>0$ on $(0,\infty)$. 
        \item Growth conditions of $\Phi$ and $\Phi'$: There exists constants $m\in(0,\infty), c\in(0,\infty)$ such that for every $\xi\in[0,\infty)$
    \begin{equation}\label{eq: growth bound on phi}
        \Phi(\xi)\leq c(1+\xi^m), \quad \text{and}\quad \Phi'(\xi)\leq c(1+\xi+\Phi(\xi)).
    \end{equation}
    Let $\Theta_{\Phi,2}$ be as in Definition \ref{def: auxhilary functions}. 
             Then, either for constants $c\in(0,\infty)$ and $\theta\in[0,1/2]$, we have for every $\xi\in(0,\infty)$,
    \begin{equation}\label{eq: first assumption on Phi}
        (\Theta'_{\Phi,2}(\xi))^{-1}=\Phi'(\xi)^{-1/2}\leq c\xi^\theta,
    \end{equation}
    or there exist constants
    $c\in(0,\infty)$, $q\in[1,\infty)$ such that for every $\xi,\eta\in[0,\infty)$
    \begin{equation}\label{eq: second assumption on Phi}
        |\xi-\eta|^q\leq c|\Theta_{\Phi,2}(\xi)-\Theta_{\Phi,2}(\eta)|^2.
    \end{equation}
        Furthermore, assume that there is a constant $c\in(0,\infty)$ such that for every $\xi\in(0,\infty)$ and $m\in(0,\infty)$ as above, we have
        \[\Theta_{\Phi,2}(\xi)\geq c(\xi^{\frac{m+1}{2}}+1).\]
        \item Growth conditions on $\sigma, \nu$ and other non-linear functions: For a constant $c\in(0,\infty)$  we have for every $\xi\in[0,\infty)$,
        \begin{equation}\label{eq: bound on sigma and so sigma squared}
            |\sigma(\xi)|\leq c\Phi^{1/2}(\xi)
        \end{equation}
    \begin{equation}\label{eq: bound on nu squared and sigma sigma' squared}
        |\nu(\xi)|^2+ |\sigma(\xi)\sigma'(\xi)|^2\leq  c(1+\xi+\Phi(\xi)).
    \end{equation}
       Furthermore, for each $\delta\in(0,1)$ there is a constant $c_\delta\in(0,\infty)$ such that for every $\xi\in(\delta,\infty)$,
    \[\frac{(\sigma'(\xi))^4}{\Phi'(\xi)}\leq c_\delta (1+\xi+\Phi(\xi)).\]
    
        \item At least linear decay of $\sigma^2$ at $0$: There exists a constant $c\in(0,\infty)$ such that 
        \[\limsup_{\xi\to 0^+}\frac{\sigma^2(\xi)}{\xi}\leq c.\]
        In particular this implies that $\sigma(0)=0$.
        \item Regularity of oscillations of of $\sigma^2$ and $\nu$ at infinity: There is a constant $c\in[1,\infty)$ such that for every $\xi\in[0,\infty)$
        \[\sup_{\xi'\in[0,\xi]}\sigma^2(\xi')\leq  c(1+\xi+\sigma^2(\xi)), \quad \text{and} \quad \sup_{\xi'\in[0,\xi]}|\nu(\xi')|\leq  c(1+\xi+|\nu(\xi)|).\]
       \item Entropy assumption: We have $\log(\Phi)$ is locally integrable on $[0,\infty)$, and that $\nu\in L^2([0,\infty);\mathbb{R}^d)$.  
    \end{enumerate}
\end{assumption}
In \cite{popat2025well} there are also several assumptions on the boundary data, see Assumption 4.2 points 9-11, Definition 4.6 and Assumption 4.9, Assumption 4.20, Assumption 4.23 and Assumption 4.28.
The assumptions encompass all non-negative constant functions including zero and all smooth functions bounded away from zero.
All of the assumptions are automatically satisfied in the central limit theorem section due to Assumption \ref{assumption of constant boundary data and constant initial condition}.
However, in the large deviations results of Section \ref{sec: ldp}, this assumption is relaxed, so the below assumptions on the boundary data are required for well-posedness there.
\begin{assumption}[Assumptions on boundary data for well-posedness]\label{asm: assumptions on boundary data for well-posedness} 
Let $\Phi,\sigma,\nu$ denote the non-linear functions in equation \eqref{generalised Dean-Kawasaki Equation Stratonovich with epsilon}, recall $\hat{\eta}=(\hat{\eta}_1,\hdots,\hat{\eta}_d)$ denotes the outward pointing unit normal at the boundary, and recall that $\bar{f}$ denotes the boundary data of various equations.
    \begin{enumerate}
        \item Boundary regularity for non-linear functions: 
        Define the anti-derivative of $\nu$, defined element-wise for $i=1,\hdots,d$, by $\Theta_{\nu,i}(0)=0$, $\Theta'_{\nu,i}(\xi)=\nu_i(\xi)$, and define the unique function $\Psi_\sigma$ given by $\Psi_\sigma(1)=0$, $\Psi'_\sigma(\xi)=F_1[\sigma'(\xi)]^2$.
    We assume that 
        \[\bar{f}\in L^2(\partial U), \quad \Phi^{-1}(\bar{f})\in L^2(\partial U), \quad \sigma^2(\Phi^{-1}(\bar{f}))\in L^1(\partial U), \quad \Psi_\sigma(\Phi^{-1}(\bar{f}))\in L^2(\partial U),\]
        and for $i=1\hdots,d$,
        \[\Theta_{\nu,i}(\Phi^{-1}(\bar{f}))\cdot\hat{\eta}_i\in L^1(\partial U).\]
        \item Higher order spatial regularity of harmonic PDE:
         Define $h:U\to\mathbb{R}$ to be the below harmonic function that captures the regularity of the solution on the boundary
     \begin{equation}\label{eq: pde h}
    \begin{cases}
            -\Delta h = 0 & \text{on} \hspace{5pt}U,\\
    h=\Phi^{-1}(\bar{f})& \text{on} \hspace{5pt}\partial U.
    \end{cases}
\end{equation}
Then we assume that $\Phi^{-1}(\bar{f})\in H^1(\partial U)$.\\
Further assume that either the boundary data is bounded or we have $\Theta_{\Phi,2}(h)\in H^1(U)$ where $\Theta_{\Phi,2}$ is as in Definition \ref{def: auxhilary functions}. 
    \end{enumerate}
The following assumptions are needed for the entropy estimates in \cite{popat2025well}.
\begin{enumerate}[resume]
    \item Entropy assumptions: For the unique function $\Theta_{\Phi,\sigma}$ defined by 
        $\Theta_{\Phi,\sigma}(1)=0$ and $\Theta'_{\Phi,\sigma}(\xi)=\frac{\Phi'(\xi)\sigma'(\xi)\sigma(\xi)}{\Phi(\xi)}$,
        and for $i=1,\hdots,d$,
        the unique functions $\Theta_{\Phi,\nu,i}$ defined by
       $ \Theta_{\Phi,\nu,i}(0)=0$, and $\Theta'_{\Phi,\nu,i}(\xi)=\frac{\Phi'(\xi)\nu_i(\xi)}{\Phi(\xi)}$,
       we have
       \[\Theta_{\Phi,\sigma}(\Phi^{-1}(\bar{f}))\in L^1(\partial U),\quad \Theta_{\Phi,\nu,i}(\Phi^{-1}(\bar{f}))\cdot \eta_i \in L^1(\partial U).\]
    \item Let $v_0:U\to\mathbb{R}$ denote the solution of the PDE
    \begin{equation*}
    \begin{cases}
       -\Delta v_0=0,&\text{in}\hspace{5pt} U,\\
       v_0=\log(\bar{f}),&\text{on}\hspace{5pt} \partial U.    \end{cases}
\end{equation*}
We assume that $\log(\bar{f})\in H^1(\partial U)$ in the sense of page 176 of \cite{fabes1978potential}.
  \item Either $F_2=0$, or the unique function $\Theta_\sigma$ defined by $\Theta_\sigma(1)=0$, $\Theta_\sigma'(\xi)=\frac{\sigma(\xi)\sigma'(\xi)}{\xi}$ satisfies
        $\Theta_\sigma(\Phi^{-1}(\bar{f})\wedge 1)\in L^1(\partial U)$.
\end{enumerate}
\end{assumption}
\begin{remark}
    In the model case, the function $\Psi_\sigma$ in point 1 is 
    \[\Psi_\sigma(\xi)=F_1\log(\xi),\]
    so the assumption requires that if the boundary data is not constant, then it is uniformly bounded away from zero.
    That is to say, we can not handle boundary data that that is zero on some part of the boundary and positive on other parts.\\
    This restriction can also be seen in points 4 and 5.
\end{remark}
\subsection{Space for initial data for the large deviations principle}\label{sec: Correct space for initial data for the large deviations principle}
In this section we produce a formal estimate for the skeleton equation \eqref{eq: skeleton equation} under the assumption that the initial data is non-negative which will suggest that $Ent_\Phi(U)$ is the correct space for the initial data. 
The estimate is a combination of the estimate in Section 2.2 of \cite{fehrman2023non} and the entropy estimate Proposition 4.24 of \cite{popat2025well}.
For simplicity we begin by assuming constant boundary data as in Assumption \ref{asm: constant boundary data}, and towards the end of the section we subsequently we will explain how this can be dispensed of and how one would make the estimate rigorous.\\
Let $\psi:[0,\infty)\to \mathbb{R}$ be an arbitrary function.
Under Assumption \ref{asm: constant boundary data}, define the function $\Psi:[0,\infty)\to\mathbb{R}$ by
\begin{equation}\label{eq: for Psi with constant boundary data}
    \Psi(\xi)=\int_0^\xi\psi(x)-\psi(\Phi^{-1}(\bar{f}))\,dx,
\end{equation}
where the second term in the integral is a constant which ensures $\Psi'$ is compactly supported.
Denoting by $\rho$ the unique solution of the skeleton equation \eqref{eq: skeleton equation} as in Theorem \ref{thm: well-posedness of skeleton equation assumed}, it follows using integration by parts that the composition $\Psi(\rho)$ satisfies
\begin{align}\label{eq: computation for identifying correct space for initial data}
    \partial_t\int_U\Psi(\rho)\,dx&=\int_U \Psi'(\rho)\partial_t\rho\,dx\nonumber\\
    &=\int_U \left(\psi(\rho)-\psi(\Phi^{-1}(\bar{f}))\right)\left(\Delta\Phi(\rho)-\nabla\cdot(\sigma(\rho)g+ \nu(\rho))\right)\,dx\nonumber\\
    &= -\int_U \psi'(\rho)\Phi'(\rho)|\nabla\rho|^2\,dx+ \int_U\psi'(\rho)\sigma(\rho)\nabla\rho\cdot g+ \psi'(\rho)\nabla\rho\cdot\nu(\rho)\,dx.
\end{align}
Our goal is to pick the test function $\psi$ in such a way that we are able to obtain a nice energy estimate.
The first term on the right hand side of \eqref{eq: computation for identifying correct space for initial data} can be moved over to the left hand side of the estimate, and the final term vanishes due Assumption \ref{asm: constant boundary data}. 
The middle term is handled using Cauchy-Schwarz inequality, Young's inequality and point 4 of Assumption \ref{asm: assumptions for well posedenss of equation} that there is a constant $c\in(0,\infty)$ such that $\sigma\leq c\Phi^{1/2}$, which provides for every $\delta\in(0,1)$ the bound
\begin{equation}\label{eq: bounding term involving product of g in skeleton equation proof}
    \int_U\psi'(\rho)\sigma(\rho)\nabla\rho\cdot g\leq\delta^{-1}\int_0^t\int_U|g|^2+c\delta\int_0^t\int_U(\psi'(\rho))^2\Phi(\rho)|\nabla\rho|^2.
\end{equation}
Subsequently putting everything together gives
\small
\begin{align*}
    \left.\int_U\Psi(\rho)\,dx\right|_{s=0}^t+\int_0^t\int_U \psi'(\rho)\Phi'(\rho)|\nabla\rho|^2\,dx\,ds\leq \delta^{-1}\int_0^t\int_U|g|^2\,dx\,ds+c\delta\int_0^t\int_U(\psi'(\rho))^2\Phi(\rho)|\nabla\rho|^2\,dx\,ds.
\end{align*}
\normalsize
We need to estimate the final term on the right hand side, which is difficult for general $\psi$.
However, we realise that this term can be absorbed onto left hand side of the estimate by picking $\delta$ small enough if there exists a constant $c\in(0,\infty)$ such that for every $\xi\in(0,\infty)$ we have
\[(\psi'(\xi))^2\Phi(\xi)\leq c\psi'(\xi)\Phi'(\xi),\]
which is equivalent to 
\[\psi'(\xi)\leq c\frac{\Phi'(\xi)}{\Phi(\xi)}.\]
By observation, a candidate choice is $\psi(\xi)=\log(\Phi(\xi))$, and consequently by the definition \eqref{eq: for Psi with constant boundary data}, we have
\begin{align*}
    \Psi_{\Phi}(\xi)=\int_0^\xi \log(\Phi(x))-\log(\bar{f})\,dx.
\end{align*}
When we assume Assumption \ref{asm: constant boundary data} this is precisely the definition of $\Psi_{\Phi}$ given in Definition \ref{def: entropy space}.
From \eqref{eq: computation for identifying correct space for initial data} we deduce that there is a constant $c\in(0,\infty)$ such that
\begin{align}\label{eq: entropy estimate for initial condition}
   \int_U\Psi_{\Phi}(\rho(\cdot,t))\,dx+\frac{1}{2}\int_0^t\int_U |\nabla\Phi^{1/2}(\rho)|^2\,dx\,ds    \,\leq \int_U\Psi_{\Phi}(\rho_0)\,dx+ c\int_0^t\int_U|g|^2\,dx\,ds.
\end{align}
If we want to dispense of Assumption \ref{asm: constant boundary data}, instead of looking at $\Psi$ of the form \eqref{eq: for Psi with constant boundary data}, we look at the function
    \[\Psi(\xi)=\int_0^\xi(\psi(x)-h_\psi(x))\,dx,\]
    where for each fixed $\psi$, $h_\psi$ is the solution of the harmonic PDE
    \begin{equation}\label{eq: harmonic PDE h}
   \begin{cases}
    -\Delta h_\psi=0, & \text{on} \hspace{5pt} U,\\
h_\psi=\psi(\Phi^{-1}(\bar{f})),&\text{on}\hspace{5pt} \partial U.
   \end{cases} 
\end{equation}
Repeating the energy estimate computation outlined in  \eqref{eq: computation for identifying correct space for initial data} gives the additional boundary dependent terms
\begin{align}\label{eq: estimate for Psi in non constant boundary case}
    \partial_t\int_U\Psi(\rho)\,dx &= -\int_U \psi'(\rho)\Phi'(\rho)|\nabla\rho|^2\,dx+ \int_U\psi'(\rho)\sigma(\rho)\nabla\rho\cdot g+ \psi'(\rho)\nabla\rho\cdot\nu(\rho)\,dx\nonumber\\
    & \quad+\int_U \nabla h_\psi(\rho)\cdot\nabla\Phi(\rho)\,dx- \int_U\sigma(\rho)\nabla h_\psi(\rho)\cdot g+ \nabla h_\psi(\rho)\cdot\nu(\rho)\,dx
\end{align}
Motivated by the estimate with constant boundary data, we wish to obtain an estimate for the choice $\psi(\xi)=\log(\Phi(\xi))$.
For this choice, we have $h_{\log(\Phi(\cdot))}=v_0$ where $v_0$ is as in point 4 of Assumption \ref{asm: assumptions on boundary data for well-posedness}.
We notice that the desired estimate can be obtained from the entropy estimate presented in Proposition 4.24 of \cite{popat2025well}, where instead of the noise term we have to bound the term involving the control $g$ instead. 
Since we only have the regularity $g\in L^2(U\times[0,T];\mathbb{R}^d)$ for the control, for both terms involving the control in \eqref{eq: estimate for Psi in non constant boundary case} we bound them using Cauchy-Schwarz and Young's inequalities as in \eqref{eq: bounding term involving product of g in skeleton equation proof}, and absorb the additional terms onto the left hand side of the estimate.\\
We explain how make the above estimate rigorous.
As we did several times already, in order for solutions to be sufficiently regular to apply It\^o's formula, for fixed $\alpha\in(0,1)$ we add the regularisation $\alpha\Delta\rho$ to the right hand side of the skeleton equation \eqref{eq: skeleton equation}.
It\^o's formula is then applied to the regularised function
\[\Psi_{\Phi,\delta}(\xi):=\int_0^\xi \log(\Phi(x)+\delta)-v_\delta(x)\,dx,\]
where $\delta\in(0,1)$ is fixed and  $v_\delta$ is the harmonic PDE with boundary data $\log(\bar{f}+\delta)$.
The estimate is uniform in both $\alpha$ and $\delta$, so that we can subsequently pass to the limits.
\subsection{Kinetic formulation for the controlled SPDE}\label{sec: kinetic formulation for controlled SPDE}
For fixed $\epsilon\in(0,1), K\in\mathbb{N}$ and control $g\in L^2(U\times[0,T];\mathbb{R}^d)$, we provide the computation of how one derives the kinetic equation for the controlled SPDE \eqref{eq: controlled SPDE}
\begin{equation}\label{eq: controlled SPDE in Appendix}
    \partial_t\rho^{\epsilon,g}=\Delta\Phi(\rho^{\epsilon,g})-\sqrt{\epsilon}\nabla\cdot(\sigma(\rho^{\epsilon,g})\circ\dot{\xi}^K))-\nabla\cdot(\sigma(\rho^{\epsilon,g})P_K g)-\nabla\cdot \nu(\rho^{\epsilon,g}).
\end{equation} 
The kinetic equation is derived in the uncontrolled setting in Section 3 of \cite{fehrman2024well}. 
For ease of notation and because $\epsilon,K,g$ are all fixed, in this section we denote the solutions of the above SPDE as $\rho^{\epsilon,K}\equiv\rho$.
Re-writing equation \eqref{eq: controlled SPDE in Appendix} using It\^o noise gives
\begin{align*}
      \partial_t\rho=\Delta\Phi(\rho)-\sqrt{\epsilon}\nabla\cdot(\sigma(\rho)\dot{\xi}^K)-\nabla\cdot(\sigma(\rho)P_K g)-\nabla\cdot \nu(\rho)+\frac{\epsilon}{2}\nabla\cdot\left(F_1^K(\sigma'(\rho))^2\nabla\rho + \sigma(\rho)\sigma'(\rho)F_2^K\right),
\end{align*}
where we recall that $F_1^K,F_2^K$ are defined in Definition \ref{definition truncated noise}.
Suppose that we were interested in how functions of the solution behaved. 
For smooth function $S:\mathbb{R}\to\mathbb{R}$, applying It\^o's formula gives that formally\footnote{the computation is only formal because the equation does not have enough regularity to apply It\^o's formula. To get around this, we would need to add a regularisation $\alpha\nabla\rho$ to the right hand side of the equation. So in general the identities would be true with equalities replaced by inequalities.}
\begin{align}\label{eq: ito formula for S}
      \partial_tS(\rho)&=S'(\rho)\Delta\Phi(\rho)-S'(\rho)\sqrt{\epsilon}\nabla\cdot(\sigma(\rho)\dot{\xi}^K)-S'(\rho)\nabla\cdot(\sigma(\rho)P_K g)-S'(\rho)\nabla\cdot \nu(\rho)\nonumber\\
      &+\frac{\epsilon}{2}S'(\rho)\nabla\cdot\left(F_1^K(\sigma'(\rho))^2\nabla\rho + \sigma(\rho)\sigma'(\rho)F_2^K\right)+\frac{\epsilon}{2}S''(\rho)\sum_{k=1}^K \left(\nabla\cdot(\sigma(\rho)f_k)\right)^2,
\end{align}
where the final term is the It\^o correction, and can be expanded to get
\[\sum_{k=1}^K \left(\nabla\cdot(\sigma(\rho)f_k)\right)^2=F_1^K|\nabla\sigma(\rho)|^2+2\sigma(\rho)\sigma'(\rho)\nabla\rho\cdot F_2^K+\sigma^2(\rho)F_3^K,\]
where again $F_3^K$ is defined in Definition \ref{definition truncated noise}.
We wish to observe cancellation between the It\^o-to-Stratonovich conversion terms and the It\^o correction term above, particularly in the singular term involving $(\sigma'(\rho))^2$.
However, we notice that the conversion terms involve $S'(\rho)$, whilst the correction terms involve $S''(\rho)$.
To get past this, for the conversion terms we will bring $S'(\rho)$ into the derivative. 
Doing this also to the first and third terms on the right hand side of equation \eqref{eq: ito formula for S} gives
\begin{align*}
      \partial_tS(\rho)&=\nabla\cdot(S'(\rho)\Phi'(\rho)\nabla\rho)-S''(\rho)\nabla\rho\cdot\nabla\Phi(\rho)-S'(\rho)\sqrt{\epsilon}\nabla\cdot(\sigma(\rho)\dot{\xi}^K)-\nabla\cdot\left(S'(\rho)\sigma(\rho)P_K g\right)\\
      &+S''(\rho)\sigma(\rho)\nabla\rho\cdot P_Kg-S'(\rho)\nabla\cdot \nu(\rho)+\frac{\epsilon}{2}\nabla\cdot\left(F_1^K(\sigma'(\rho))^2S'(\rho)\nabla\rho + \sigma(\rho)\sigma'(\rho)S'(\rho)F_2^K\right)\nonumber\\
      &-\frac{\epsilon}{2}\left(S''(\rho)\nabla\rho\cdot F_1^K(\sigma'(\rho))^2\nabla\rho+S''(\rho)\nabla\rho\cdot\sigma(\rho)\sigma'(\rho)F_2^K\right)\\
      &+\frac{\epsilon}{2}S''(\rho)\left(F_1^K|\nabla\sigma(\rho)|^2+2\sigma(\rho)\sigma'(\rho)\nabla\rho\cdot F_2^K+\sigma^2(\rho)F_3^K\right).
\end{align*}
After the cancellation of terms in the final two lines, the equation simplifies to
\begin{align*}
      \partial_tS(\rho)&=\nabla\cdot(S'(\rho)\Phi'(\rho)\nabla\rho)-S''(\rho)\nabla\rho\cdot\nabla\Phi(\rho)-S'(\rho)\sqrt{\epsilon}\nabla\cdot(\sigma(\rho)\dot{\xi}^K)-\nabla\cdot\left(S'(\rho)\sigma(\rho)P_K g\right)\\
      &+S''(\rho)\sigma(\rho)\nabla\rho\cdot P_Kg-S'(\rho)\nabla\cdot \nu(\rho)+\frac{\epsilon}{2}\nabla\cdot\left(F_1^K(\sigma'(\rho))^2S'(\rho)\nabla\rho + \sigma(\rho)\sigma'(\rho)S'(\rho)F_2^K\right)\nonumber\\
      &+\frac{\epsilon}{2}S''(\rho)\left(\sigma(\rho)\sigma'(\rho)\nabla\rho\cdot F_2^K+\sigma^2(\rho)F_3^K\right).
\end{align*}
Rigorously we interpret the above equation by integrating in space and time against a test function $\psi\in C_c^\infty(U)$ and subsequently integrating by parts which gives for every fixed $t\in[0,T]$
\begin{align*}
      \left.\int_U\psi(x) S(\rho)\right|_{s=0}^t&=-\int_0^t\int_U\nabla\psi\cdot(S'(\rho)\Phi'(\rho)\nabla\rho)-\int_0^t\int_U\psi S''(\rho)\nabla\rho\cdot\nabla\Phi(\rho)\\
      &-\sqrt{\epsilon}\int_0^t\int_U\psi S'(\rho)\nabla\cdot(\sigma(\rho)\,d\xi^K)+\int_0^t\int_U\nabla\psi \cdot\left(S'(\rho)\sigma(\rho)P_K g\right)\\
      &+\int_0^t\int_U \psi S''(\rho)\sigma(\rho)\nabla\rho\cdot P_Kg-\int_0^t\int_U 
 \psi S'(\rho)\nabla\cdot \nu(\rho)\\
 &-\frac{\epsilon}{2}\int_0^t\int_U \nabla\psi \cdot\left(F_1^K(\sigma'(\rho))^2S'(\rho)\nabla\rho + \sigma(\rho)\sigma'(\rho)S'(\rho)F_2^K\right)\nonumber\\
      &+\frac{\epsilon}{2}\int_0^t\int_U \psi S''(\rho)\left(\sigma(\rho)\sigma'(\rho)\nabla\rho\cdot F_2^K+\sigma^2(\rho)F_3^K\right).
\end{align*}
We want to write the above equation in terms of the kinetic function.
\begin{definition}[Kinetic function]\label{def: kinetic function}
Given a non-negative solution $\rho$ of the equation \eqref{eq: controlled SPDE in Appendix}, the kinetic function $\chi:U\times(0,\infty)\times[0,T]\to\{0,1\}$ of $\rho$ is defined as
\[\chi(x,\xi,t):=\mathbbm{1}_{0\leq\xi\leq\rho(x,t)}.\]
\end{definition}
The variable $\xi\in(0,\infty)$ is called the velocity variable.
To re-write the equation in terms of the kinetic function, we realise that if $S(0)=0$, we have that 
\[S(\rho)=\int_{\mathbb{R}}S'(\xi)\chi(x,\xi,t)\,d\xi.\]
This can be substituted into the left hand side, and for the right hand side we formally introduce artificial integrals in the velocity variable by adding a Dirac delta term $\delta_0(\xi-\rho)$, which gives
\begin{align*}
      &\left.\int_\mathbb{R}\int_U\psi(x) S'(\xi)\chi(x,\xi,t)\,dx\,d\xi\right|_{s=0}^t=-\int_\mathbb{R}\int_0^t\int_U\delta_0(\xi-\rho)S'(\xi)\Phi'(\xi)\nabla\psi\cdot\nabla\rho\\
      &-\int_\mathbb{R}\int_0^t\int_U\delta_0(\xi-\rho)\Phi'(\xi)\psi S''(\xi)|\nabla\rho|^2-\sqrt{\epsilon}\int_\mathbb{R}\int_0^t\int_U\delta_0(\xi-\rho)\psi S'(\xi)\nabla\cdot(\sigma(\rho)\,d\xi^K)\\
      &+\int_{\mathbb{R}}\int_0^t\int_U\delta_0(\xi-\rho)S'(\xi)\sigma(\xi)\nabla\psi \cdot P_K g+\int_{\mathbb{R}}\int_0^t\int_U\delta_0(\xi-\rho) \psi S''(\xi)\sigma(\xi)\nabla\rho\cdot P_Kg\\
      &-\int_{\mathbb{R}}\int_0^t\int_U \delta_0(\xi-\rho)
 \psi S'(\xi)\nabla\cdot \nu(\rho)\\
 &-\frac{\epsilon}{2}\int_\mathbb{R}\int_0^t\int_U \delta_0(\xi-\rho)\nabla\psi \cdot\left(F_1^K(\sigma'(\xi))^2S'(\xi)\nabla\rho + \sigma(\xi)\sigma'(\xi)S'(\xi)F_2^K\right)\nonumber\\
      &+\frac{\epsilon}{2}\int_{\mathbb{R}}\int_0^t\int_U \delta_0(\xi-\rho)\psi S''(\xi)\left(\sigma(\xi)\sigma'(\xi)\nabla\rho\cdot F_2^K+\sigma^2(\xi)F_3^K\right).
\end{align*}
Finally, we can re-write the equation, factoring terms involving $\psi(x)S'(\xi)$,
\begin{align*}
      &\left.\int_\mathbb{R}\int_U\psi(x) S'(\xi)\chi(x,\xi,t)\,dx\,d\xi\right|_{s=0}^t=-\int_\mathbb{R}\int_0^t\int_U \nabla\left(\psi S'(\xi)\right)\cdot\delta_0(\xi-\rho)\Phi'(\xi)\nabla\rho\\
      &-\int_\mathbb{R}\int_0^t\int_U\partial_\xi\left(\psi S'(\xi)\right)\delta_0(\xi-\rho)\Phi'(\xi)|\nabla\rho|^2-\sqrt{\epsilon}\int_\mathbb{R}\int_0^t\int_U\psi S'(\xi)\delta_0(\xi-\rho)\nabla\cdot(\sigma(\rho)\,d\xi^K)\\
      &+\int_{\mathbb{R}}\int_0^t\int_U\nabla\left(\psi S'(\xi)\right)\cdot\delta_0(\xi-\rho)\sigma(\xi) P_K g+\int_{\mathbb{R}}\int_0^t\int_U\partial_\xi(\psi S'(\xi))\delta_0(\xi-\rho)\sigma(\xi)\nabla\rho\cdot P_Kg\\
      &-\int_{\mathbb{R}}\int_0^t\int_U \psi S'(\xi)\delta_0(\xi-\rho)
 \nabla\cdot \nu(\rho)\\
 &-\frac{\epsilon}{2}\int_\mathbb{R}\int_0^t\int_U \nabla\left(\psi S'(\xi)\right)\cdot\delta_0(\xi-\rho)\left(F_1^K(\sigma'(\xi))^2\nabla\rho + \sigma(\xi)\sigma'(\xi)F_2^K\right)\nonumber\\
      &+\frac{\epsilon}{2}\int_{\mathbb{R}}\int_0^t\int_U\partial_\xi\left(\psi S'(\xi)\right) \delta_0(\xi-\rho)\left(\sigma(\xi)\sigma'(\xi)\nabla\rho\cdot F_2^K+\sigma^2(\xi)F_3^K\right).
\end{align*}
Integrating by parts and using the density of functions of the type $\psi(x)S'(\xi)$ in $C_c^\infty(U\times(0,\infty))$ allows us to conclude that the kinetic function $\chi$ of equation \eqref{eq: controlled SPDE in Appendix} solves
\begin{align}\label{eq: for kinetic function before we introduce kinetic measure}
    \partial_t\chi&=\nabla\cdot(\delta_0(\xi-\rho)\Phi'(\xi)\nabla\rho)+\partial_\xi(\delta_0(\xi-\rho)\Phi'(\xi)|\nabla\rho|^2)-\sqrt{\epsilon}\delta_0(\xi-\rho)\nabla\cdot(\sigma(\rho) \dot{\xi}^K)\nonumber\\
    &-\nabla\cdot(\delta_0(\xi-\rho)\sigma(\xi)P_K g)-\partial_\xi(\delta_0(\xi-\rho)\sigma(\xi)\nabla\rho\cdot P_Kg)-\delta_0(\xi-\rho)\nabla\cdot\nu(\rho)\nonumber\\
    &+\frac{\epsilon}{2}\nabla(\delta_0(\xi-\rho)\left(F_1^K(\sigma'(\xi))^2\nabla\rho + \sigma(\xi)\sigma'(\xi)F_2^K\right))\nonumber\\
    &-\frac{\epsilon}{2}\partial_\xi(\delta_0(\xi-\rho)\left(\sigma(\xi)\sigma'(\xi)\nabla\rho\cdot F_2^K+\sigma^2(\xi)F_3^K\right)).
\end{align}
However, we mentioned previously that the above computation is formal because the solution $\rho$ did not have enough regularity to enable us to apply It\^o's formula to $S(\rho)$.
We mentioned that this is resolved by adding a $\alpha\Delta\rho$ term to the original equation for $\alpha\in(0,1)$.
This results in a non-negative term appearing on the right hand side of \eqref{eq: for kinetic function before we introduce kinetic measure}, so the above equality should actually be replaced with an inequality \say{$\leq$}.
At the end of the analysis, we take $\alpha\to 0$, where we have competing convergence to zero of $\alpha$ and divergence to infinity of $|\nabla\rho|^2$, see for example the definition of the kinetic measures $dq^\epsilon$ in the proof of Theorem \ref{thm: convergence of solutions when initial data and controls converge}.\\
On the kinetic level, the entropy inequality is quantified exactly by the kinetic measure, see \cite{chen2003well}.
That is, a non-negative measure on $U\times(0,\infty)\times[0,T]$ satisfying that in the sense of measures
\begin{equation}\label{eq: inequality for kinetic measure}
    \delta_0(\xi-\rho)\Phi'(\xi)|\nabla\rho|^2\leq q.
\end{equation}
\begin{definition}[Kinetic measure]\label{def: kinetic measure}
Let $(\Omega,\mathcal{F},(\mathcal{F}_t)_{t\geq0},\mathbb{P})$ be the filtered probability space from Definition \ref{definition truncated noise}.
A kinetic measure $q$ is a map from $\Omega$ to the space of non-negative, locally finite measures on $U\times(0,\infty)\times[0,T]$ that satisfies
\[(\omega,t)\in\Omega\times[0,T]\to\int_\mathbb{R}\int_0^t\int_U\psi(x,\xi)dq(\omega)\quad \text{is} \quad \mathcal{F}_t-\text{predictable}\]
for every test function $\psi\in C_c^\infty(U\times(0,\infty))$.
\end{definition}
Using equation \eqref{eq: for kinetic function before we introduce kinetic measure} and the inequality \eqref{eq: inequality for kinetic measure} relating to the kinetic measure, the kinetic equation is then defined for every $\psi\in C^\infty_c(U\times(0,\infty))$ almost surely for every $t\in[0,T]$, by
\begin{align}\label{eq: kinetic equation}
    \int_{\mathbb{R}}&\int_U\chi(x,\xi,t)\psi(x,\xi)=\int_{\mathbb{R}}\int_U\chi(x,\xi,0)\psi(x,\xi)-\int_0^t\int_U\Phi'(\rho)\nabla\rho\cdot(\nabla_x\psi)(x,\rho)-\int_0^t\int_\mathbb{R}\int_U\partial_\xi\psi(x,\xi)\,dq\nonumber\\
    &-\sqrt{\epsilon}\int_0^t\int_U\psi(x,\rho)\nabla\cdot(\sigma(\rho) \,d{\xi}^K)+\int_0^t\int_U\sigma(\rho)(\nabla_x\psi)(x,\rho)\cdot P_K g+\int_0^t\int_U(\partial_\xi\psi)(x,\rho)\sigma(\rho)\nabla\rho\cdot P_Kg\nonumber\\
    &-\int_0^t\int_U\psi(x,\rho)\nabla\cdot\nu(\rho)-\frac{\epsilon}{2}\int_0^t\int_U(\nabla_x\psi)(x,\rho)\cdot\left(F_1^K(\sigma'(\rho))^2\nabla\rho + \sigma(\rho)\sigma'(\rho)F_2^K\right)\nonumber\\
    &+\frac{\epsilon}{2}\int_0^t\int_U(\partial_\xi\psi)(x,\rho)\left(\sigma(\rho)\sigma'(\rho)\nabla\rho\cdot F_2^K+\sigma^2(\rho)F_3^K\right),
\end{align}
where we used the notation $\nabla_x\psi$ to reflect the fact that we are taking the gradient only in the first component rather than the full gradient. 
The compact support of the test function $\psi$ in the velocity variable amounts to a renormalisation of the solution away from its zero set, and so the second term on the third line is well defined.
\newpage
\Large{\textbf{Acknowledgements}}\\
\normalsize
The author is grateful to his PhD advisor Benjamin Fehrman for the constant academic support and guidance. This work would
not have been possible without him.
The author is grateful for support from the EPSRC Centre for Doctoral Training in Mathematics of Random Systems: Analysis, Modelling and Simulation (EP/S023925/1).and 
\bibliographystyle{alpha}  
\bibliography{references.bib}
\end{document}